\documentclass{article}
\usepackage{latexsym}
\usepackage{graphics}
\usepackage{amsmath}
\usepackage{amsfonts}
\begin{document}
\bibliographystyle{plain}
\title{The decomposition of Global Conformal Invariants V.}
\author{Spyros Alexakis\thanks{University of Toronto, alexakis@math.toronto.edu.
\newline
This work has absorbed
the best part of the author's energy over many years. 
This research was partially conducted during 
the period the author served as a Clay Research Fellow, 
an MSRI postdoctoral fellow,
a Clay Liftoff fellow and a Procter Fellow.  
\newline
The author is immensely indebted to Charles
Fefferman for devoting twelve long months to the meticulous
proof-reading of the present paper. He also wishes to express his
gratitude to the Mathematics Department of Princeton University
for its support during his work on this project.}}
\date{}
\maketitle
\newtheorem{proposition}{Proposition}
\newtheorem{theorem}{Theorem}
\newcommand{\Sum}{\sum}
\newtheorem{lemma}{Lemma}
\newtheorem{observation}{Observation}
\newtheorem{formulation}{Formulation}
\newtheorem{definition}{Definition}
\newtheorem{conjecture}{Conjecture}
\newtheorem{corollary}{Corollary}
\numberwithin{equation}{section}
\numberwithin{lemma}{section}
\numberwithin{theorem}{section}
\numberwithin{definition}{section}
\numberwithin{proposition}{section}

\begin{abstract}
This is the fifth in a series of papers where we  prove a conjecture of Deser and Schwimmer
regarding the algebraic structure of ``global conformal
invariants''; these  are defined to
be conformally invariant integrals of geometric scalars.
 The conjecture asserts that the integrand of
any such  integral can be expressed as a linear
combination of a local conformal invariant, a divergence and of
the Chern-Gauss-Bonnet integrand. 

The present paper complements \cite{alexakis4} in 
reducing the purely algebraic results that were 
used in \cite{alexakis1,alexakis2} to certain simpler Lemmas, 
which will be proven in the last paper in this series, \cite{alexakis6}.
\end{abstract}

\tableofcontents
\section{Introduction}

This is the fourth in a series of papers 
\cite{alexakis1,alexakis3,alexakis4,alexakis5,alexakis6} 
where we prove a conjecture of Deser-Schwmimmer \cite{ds:gccaad} regarding the algebraic 
structure of global conformal invariants. We recall that a global 
conformal invariant is an integral of a natural scalar-valued function of 
Riemannian metrics, $\int_{M^n}P(g)dV_g$, with the property that this integral 
remains invariant under conformal re-scalings of the underlying 
metric.\footnote{See the introduction of \cite{alexakis1}
for a detailed discussion of the Deser-Schwimmer 
conjecture, and for background on scalar Riemannian invariants.} 
More precisely, $P(g)$ is assumed to be a linear combination, $P(g)=\sum_{l\in L} a_l C^l(g)$, 
where each $C^l(g)$ is a complete contraction in the form:

\begin{equation}
\label{contraction} 
contr^l(\nabla^{(m_1)}R\otimes\dots\otimes\nabla^{(m_s)}R);
\end{equation}
 here each factor $\nabla^{(m)}R$ stands for the $m^{th}$ iterated 
covariant derivative of the curvature tensor $R$. $\nabla$ is the Levi-Civita 
connection of the metric $g$ and $R$ is the curvature associated to this connection. 
The contractions are taken with respect to the quadratic form $g^{ij}$.
In this series of papers we prove:

\begin{theorem}
\label{thetheorem} 
Assume that $P(g)=\sum_{l\in L} a_l C^l(g)$, where each $C^l(g)$ is a 
complete contraction in the form (\ref{contraction}), with weight $-n$. 
Assume that for every closed Riemannian manifold $(M^n,g)$ and every $\phi\in C^\infty (M^n)$:
$$\int_{M^n}P(e^{2\phi}g)dV_{e^{2\phi}g}=\int_{M^n}P(g)dV_g.$$

We claim that $P(g)$ can then be expressed in the form:
$$P(g)=W(g)+div_iT^i(g)+\operatorname{Pfaff}(R_{ijkl}).$$
Here $W(g)$ stands for a local conformal invariant of weight $-n$  (meaning 
that $W(e^{2\phi}g)=e^{-n\phi}W(g)$ for every $\phi\in C^\infty (M^n)$), 
$div_iT^i(g)$ is the divergence of a Riemannian vector field of 
weight $-n+1$, and   $\operatorname{Pfaff}(R_{ijkl})$ is the Pfaffian of the curvature tensor. 
\end{theorem}

Before we discuss the position of the present paper in the 
series \cite{alexakis1}--\cite{alexakis6}, we digress to 
describe the relation between the present series 
of papers with classical and recent work on 
scalar local invariants in various geometries.

{\bf Broad Discussion:} The theory of {\it local} invariants of Riemannian structures 
(and indeed, of more general geometries,
e.g.~conformal, projective, or CR)  has a long history. 
As discussed in \cite{alexakis1}, the original foundations of this 
field were laid in the work of Hermann Weyl and \'Elie Cartan, see \cite{w:cg, cartan}. 
The task of writing out local invariants of a given geometry is intimately connected
with understanding polynomials in a space of tensors with  given symmetries;  
these polynomials are required to remain invariant under the action of a Lie group
on the components of the tensors. 
In particular, the problem of writing down all 
 local Riemannian invariants reduces to understanding 
the invariants of the orthogonal group. 

 In more recent times, a major program was laid out by C.~Fefferman in \cite{f:ma}
aimed at finding all scalar local invariants in CR geometry. This was motivated 
by the problem of understanding the  
local invariants which appear in the asymptotic expansion of the 
Bergman and Szeg\"o kernels of strictly pseudo-convex CR manifolds,
 in a similar way to which Riemannian invariants appear in the asymptotic expansion  
of the heat kernel; the study of the local invariants
in the singularities of these kernels led to important breakthroughs 
in \cite{beg:itccg} and more recently by Hirachi in \cite{hirachi1}.
 This program was later extended  to conformal geometry in \cite{fg:ci}. 
Both these geometries belong to a 
broader class of structures, the
{\it parabolic geometries}; these admit a principal bundle with 
structure group a parabolic subgroup $P$ of a semi-simple 
Lie group $G$, and a Cartan connection on that principle bundle 
(see the introduction in \cite{cg1}). 
An important question in the study of these structures 
is the problem of constructing all their local invariants, which 
can be thought of as the {\it natural, intrinsic} scalars of these structures.

  In the context of conformal geometry, the first (modern) landmark 
in understanding {\it local conformal invariants} was the work of Fefferman 
and Graham in 1985 \cite{fg:ci},
where they introduced the {\it ambient metric}. This allows one to 
construct local conformal invariants of any order in odd 
dimensions, and up to order $\frac{n}{2}$ in even dimensions. 
The question is then whether {\it all} invariants arise via this construction. 

The subsequent work of Bailey-Eastwood-Graham \cite{beg:itccg} proved that 
this is indeed true in odd dimensions; in even dimensions, 
they proved that the result holds  
when the weight (in absolute value) is bounded by the dimension. The ambient metric construction 
in even dimensions was recently extended by Graham-Hirachi, \cite{grhir}; this enables them to 
indentify in a satisfactory way {\it all} local conformal invariants, 
even when the weight (in absolute value) exceeds the dimension.  

 An alternative 
construction of local conformal invariants can be obtained via the {\it tractor calculus} 
introduced by Bailey-Eastwood-Gover in \cite{bego}. This construction bears a strong 
resemblance to the Cartan conformal connection, and to 
the work of T.Y.~Thomas in 1934, \cite{thomas}. The tractor 
calculus has proven to be very universal; 
tractor buncles have been constructed \cite{cg1} for an entire class of parabolic geometries. 
The relation betweeen the conformal tractor calculus and the Fefferman-Graham 
ambient metric  has been elucidated in \cite{cg2}.

The present work \cite{alexakis1}--\cite{alexakis6}, while pertaining to the question above
(given that it ultimately deals with the algebraic form of local 
{\it Riemannian} and {\it conformal} invariants), nonetheless addresses a different 
{\it type} of problem:  We here consider Riemannian invariants $P(g)$ for 
which the {\it integral} $\int_{M^n}P(g)dV_g$ remains invariant 
under conformal changes of the underlying metric; we then seek to understand 
the possible algebraic form of the {\it integrand} $P(g)$, 
ultimately proving that it can be de-composed 
in the way that Deser and Schwimmer asserted. 
It is thus not surprising that the prior work on 
 the construction and understanding of local {\it conformal} 
invariants plays a central role in this endeavor, in 
the papers \cite{alexakis2,alexakis3}.

On the other hand, a central element of our
proof are the main algebraic 
Propositions 5.1, 3.1, 3.2 in \cite{alexakis1,alexakis2}; 
these deal {\it exclusively} 
with algebraic properties of the {\it classical} scalar Riemannian 
invariants.\footnote{These ``main algebraic propositions'' are 
discussed in brief below. A {\it generalization} of these Propositions 
is the Proposition \ref{giade} below.}
The ``fundamental Proposition \ref{giade}'' makes no reference to 
integration; it is purely a statement concerning 
{\it local Riemannian invariants}. Thus, 
while the author was led to led to the main algebraic 
Propositions in \cite{alexakis1,alexakis2} 
out of the strategy that he felt was necessary to 
solve the Deser-Schwimmer conjecture, they can 
be thought of as results of an independent interest. 
The {\it proof} of these Propositions, presented
 in the present paper and in \cite{alexakis6} is in fact 
not particularily intuitive. It is the author's 
sincere hope that deeper insight (and hopefuly a more intuitive proof) 
will be obtained in the future as to  {\it why} these algebraic 
 Propositions hold. 
\newline

We now discuss the position of the present paper in this series of papers.

The purpose of the present paper is to complete the part IIA in this series:
 In \cite{alexakis1,alexakis2,alexakis3}
 we proved that the Deser-Schwimmer conjecture 
holds, {\it provided} one can show certain ``Main algebraic propositions'', namely 5.2 in 
\cite{alexakis1} and 3.1, 3.2  in \cite{alexakis2}.  
In \cite{alexakis4} we claimed a more general Proposition which {\it implies} 
Proposition 5.2 in \cite{alexakis1} and 
Propositions 3.1, 3.2 in \cite{alexakis2}; 
this new ``fundamental Proposition'' 2.1 in \cite{alexakis4} 
is to be proven by an induction of four parameters.  
In \cite{alexakis4} we also reduced the inductive step of
 Proposition 2.1 to three 
Lemmas (in particular we distinguished cases I,II,III on 
Proposition 2.1 by examiniming 
the tensor fields appearing in its hypothesis, see (\ref{hypothese2}) below; Lemmas 
 3.1, 3.2, 3.5 in \cite{alexakis4}
correspond to these three cases). We proved that these three Lemmas 
{\it imply} the inductive step of the fundamental Proposition in cases I,II,III respectively,  apart from certain 
{\it special cases} which were deferred to the present paper. In these 
special cases we will derive Proposition 2.1 in \cite{alexakis4} directly,\footnote{By this 
we mean {\it without} recourse to the 
Lemmas 3.1, 3.2, 3.5 in \cite{alexakis4}.} 
in section \ref{specialcases}. Now, in {\it proving} 
that the inductive step of Proposition \ref{giade}
follows from Lemmas 3.1, 3.2, 3.5 in \cite{alexakis4} we 
asserted certain {\it technical Lemmas}, whose 
proof was deferred to the present paper. 
These were Lemmas 4.6, 4.8, 
and 4.7, 4.9 in \cite{alexakis4}; also, the proof of Lemma 5.1 in \cite{alexakis4}
was deferred to the present paper. We prove all these Lemmas from  
\cite{alexakis4} in section \ref{technical.lemmas}. 
\newline

For reference purposes, and for the reader's convenience, 
we recall the precise formulation of the ``fundamental Proposition'' 2.1 in 
\cite{alexakis4}, referring the reader to \cite{alexakis4} 
for a definition of many of the terms appearing in the formulation. First however, 
we will recall (schematically) the  
``main algebraic Proposition'' 5.2 in \cite{alexakis1}; 
this is a {\it special case}  of Proposition 2.1 in \cite{alexakis4}, and 
provides a simpler version of it. 
\newline

{\it A simpler version of Proposition 2.1 in \cite{alexakis4}:} 
Given a Riemannian metric $g$ over an $n$-dimensional 
 manifold $M^n$ and auxilliary $C^\infty$ scalar-valued functions 
$\Omega_1,\dots,\Omega_p$ defined over $M^n$, the 
objects of study are linear combinations of tensor fields 
$\sum_{l\in L} a_l C^{l,i_1\dots i_\alpha}_g$, where each 
$C^{l,i_1\dots i_\alpha}_g$ is a {\it partial contraction} 
with $\alpha$ free indices, in the form: 
\begin{equation}
\label{gen.form1}
 pcontr(\nabla^{(m)}R\otimes\dots\otimes\nabla^{(m_s)}R\otimes
\nabla^{(b_1)}\Omega_1\otimes\dots\otimes\nabla^{(b_m)}\Omega_p);
\end{equation}
here $\nabla^{(m)}R$ stands for the $m^{th}$ covariant derivative of 
the curvature tesnor $R$,\footnote{In particular it is a tensor of rank $m+4$; 
 if we write out its free 
indices it would be in the form 
$\nabla^{(m)}_{r_1\dots r_m}R_{ijkl}$.} and 
$\nabla^{(b)}\Omega_h$ stands for the $b^{th}$ covariant 
derivative of the function $\Omega_h$. A {\it partial contraction} 
means that we have list of pairs of indices $({}_a,{}_b),\dots, ({}_c,{}_d)$ 
in (\ref{gen.form1}), 
which are contracted against each other using the  metric $g^{ij}$. The remaining 
 indices (which are not contracted against another index in (\ref{gen.form1})) are the 
{\it free indices} ${}_{i_1},\dots, {}_{i_\alpha}$. 

The ``main algebraic Proposition'' 5.2 in \cite{alexakis1} (roughly) asserts the following: Let 
$\sum_{l\in L_{\mu}} a_l C^{l,i_1\dots i_\mu}_g$ stand for a 
linear combination of partial contractions in the form (\ref{gen.form1}), 
where each $C^{l,i_1\dots i_\mu}_g$ has a given number $\sigma_1$ of factors and a given number 
$p$ of factor $\nabla^{(b)}\Omega_h$. Assume also that $\sigma_1+p\ge 3$, 
each $b_i\ge 2$,\footnote{This means that 
each function $\Omega_h$ is differentiated at least twice.} and that for each contracting 
 pair of 
indices $({}_a,{}_b)$ in any given $C^{l,i_1\dots i_\alpha}_g$, 
the indices ${}_a,{}_b$ do not belong to the 
same factor. Assume also the rank 
$\mu>0$ is fixed and each partial contraction $C^{l,i_1\dots i_\mu}_g, l\in L_\mu$ 
has a given {\it weight} $-n+\mu$.\footnote{See \cite{alexakis1} 
for a precise definition of weight.} 
Let also $\sum_{l\in L_{>\mu}} a_l C^{l,i_1\dots i_{y_l}}_g$ stand for 
a (formal) linear combination of partial contractions of 
weight $-n+y_l$, with all the properties of the 
terms indexed in $L_\mu$, {\it except} that now all the partial 
contractions have a different rank $y_l$, and each $y_l>\mu$. 

The assumption of the ``main algebraic Proposition'' 5.1 in \cite{alexakis1} is 
a local equation in the form:

\begin{equation}
\label{assumption.1} 
\sum_{l\in L_\mu} a_l Xdiv_{i_1}\dots Xdiv_{i_\mu}C^{l,i_1\dots i_\mu}_g+
\sum_{l\in L_{>\mu}} a_l Xdiv_{i_1}\dots Xdiv_{i_{y_l}}C^{l,i_1\dots i_{y_l}}_g=0,
\end{equation}
which
is assumed to hold {\it modulo} complete contractions with $\sigma+1$ factors.
Here given a partial contraction $C^{l,i_1\dots i_\alpha}_g$ 
in the form (\ref{gen.form1}) $Xdiv_{i_s}[C^{l,i_1\dots i_\alpha}_g]$
stands for sum of $\sigma-1$ terms in $div_{i_s}[C^{l,i_1\dots i_\alpha}_g]$
where the derivative $\nabla^{i_s}$ is {\it not} 
allowed to hit the factor to which the free index ${}_{i_s}$ belongs.\footnote{Recall that given a partial contraction
$C^{l,i_1\dots i_\alpha}_g$ in the form (\ref{gen.form1}) 
with $\sigma$ factors, $div_{i_s}C^{l,i_1\dots i_\alpha}_g$ 
is a sum of $\sigma$ partial contractions of rank $\alpha-1$. 
The first summand arises by adding a derivative $\nabla^{i_s}$ 
onto the first factor $T_1$ and then contracting the upper index ${}^{i_s}$
against the free index ${}_{i_s}$; the second summand 
arises by adding a derivative $\nabla^{i_s}$ 
onto the second factor $T_2$ and then contracting the upper index ${}^{i_s}$
against the free index ${}_{i_s}$ etc.}

Proposition 5.2 in \cite{alexakis1} 
then asserts that there will exist a linear combination of partial 
contactions in the form (\ref{gen.form1}), 
$\sum_{h\in H} a_h C^{h,i_1\dots i_{\mu+1}}_g$ with all the properties 
of the terms indexed in $L_{>\mu}$, and all with rank $(\mu+1)$, so that:

\begin{equation}
\label{assumption.1} 
\sum_{l\in L_1} a_l C^{l,(i_1\dots i_\mu)}_g+
\sum_{h\in H} a_h  Xdiv_{i_{\mu+1}}
C^{l,(i_1\dots i_\mu)i_{\mu+1}}_g=0;
\end{equation}
the above holds modulo terms of length $\sigma+1$. The symbol $(\dots)$ 
means that we are {\it symmetrizing} over the indices between prentheses. 
\newline

 In \cite{alexakis4} we set up a multiple induction by which we 
 will prove Proposition 5.2 in \cite{alexakis1} (outlined above) 
 and also the main algebraic Propositions 3.1, 3.2 in 
 \cite{alexakis2}. The generalized proposition \ref{giade}
 which we formulated in \cite{alexakis4} 
 deals with tensor fields in the forms:
 \begin{equation}
\label{form1}
\begin{split}
&pcontr(\nabla^{(m_1)}R_{ijkl}\otimes\dots\otimes\nabla^{(m_s)}R_{ijkl}
\otimes
\\& \nabla^{(b_1)}\Omega_1\otimes\dots\otimes \nabla^{(b_p)}\Omega_p
\otimes\nabla\phi_1\otimes\dots \otimes\nabla\phi_u),
\end{split}
\end{equation}

\begin{equation}
\label{form2}
\begin{split}
&pcontr(\nabla^{(m_1)}R_{ijkl}\otimes\dots\otimes\nabla^{(m_{\sigma_1})}
R_{ijkl}\otimes \\&S_{*}\nabla^{(\nu_1)}R_{ijkl}\otimes\dots\otimes
S_{*}\nabla^{(\nu_t)} R_{ijkl}\otimes
\\& \nabla^{(b_1)}\Omega_1\otimes\dots\otimes
\nabla^{(b_p)}\Omega_p\otimes
\\& \nabla\phi_{z_1}\dots \otimes\nabla\phi_{z_w}\otimes\nabla
\phi'_{z_{w+1}}\otimes
\dots\otimes\nabla\phi'_{z_{w+d}}\otimes\dots \otimes
\nabla\tilde{\phi}_{z_{w+d+1}}\otimes\dots\otimes\nabla\tilde{\phi}_{z_{w+d+y}}).
\end{split}
\end{equation}
(See the introduction in \cite{alexakis4} for a detailed 
description of the above form). We remark that a complete or partial contraction 
in the above form will be called ``acceptable'' if each $b_i\ge 2$, 
for $1\le i\le p$.\footnote{In other words, 
we are requiring each function $\Omega_i$ is differentiated at least twice.} 
This convention was introduced in \cite{alexakis4}.

 The claim of Proposition is a generalization of the 
``main algebraic Proposition'' in \cite{alexakis1}:

\begin{proposition}
\label{giade}
Consider two linear combinations of acceptable tensor fields in
the form (\ref{form2}):

$$\Sum_{l\in L_\mu} a_l
C^{l,i_1\dots i_{\mu}}_{g} (\Omega_1,\dots
,\Omega_p,\phi_1,\dots ,\phi_u),$$

$$\Sum_{l\in L_{>\mu}} a_l
C^{l,i_1\dots i_{\beta_l}}_{g} (\Omega_1,\dots
,\Omega_p,\phi_1,\dots ,\phi_u),$$
 where each tensor field above has real length $\sigma\ge 3$ and a given
simple character $\vec{\kappa}_{simp}$. We assume  that for each
$l\in L_{>\mu}$,  $\beta_l\ge \mu+1$. We also assume 
 that none of the tensor fields of maximal refined double
character in $L_\mu$ are ``forbidden'' (see Definition (2.12)).

 We denote by
$$\Sum_{j\in J} a_j C^j_{g}(\Omega_1,\dots ,\Omega_p,
\phi_1,\dots ,\phi_u)$$ a generic linear combination of complete
contractions (not necessarily acceptable) in the form
(\ref{form1}) that are simply subsequent to
$\vec{\kappa}_{simp}$.\footnote{Of course if $Def(\vec{\kappa}_{simp})=\emptyset$ then  by
definition $\Sum_{j\in J} \dots=0$.} We assume that:

\begin{equation}
\label{hypothese2}
\begin{split}
&\Sum_{l\in L_1} a_l Xdiv_{i_1}\dots Xdiv_{i_{\alpha}}
C^{l,i_1\dots i_{\alpha}}_{g} (\Omega_1,\dots
,\Omega_p,\phi_1,\dots ,\phi_u)+
\\&\Sum_{l\in L_2} a_l Xdiv_{i_1}\dots Xdiv_{i_{\beta_l}}
C^{l,i_1\dots i_{\beta_l}}_{g} (\Omega_1,\dots
,\Omega_p,\phi_1,\dots ,\phi_u)+
\\& \Sum_{j\in J} a_j
C^j_{g}(\Omega_1,\dots ,\Omega_p,\phi_1,\dots ,\phi_u)=0.
\end{split}
\end{equation}

\par We draw our conclusion with a little more notation: We break the index set
$L_\mu$ into subsets $L^z, z\in Z$, ($Z$ is finite)
with the rule that each $L^z$  indexes tensor fields
with the same refined double character, and conversely two tensor
fields with the same refined double character must be indexed in
the same $L^z$. For each index set $L^z$, we denote the
refined double character in question by $\vec{L}^z$. Consider
the subsets $L^z$ that index the tensor fields of {\it maximal} refined double
 character.\footnote{Note that in any set $S$ of $\mu$-refined double characters
 with the same simple character there is going to be a subset $S'$
 consisting of the maximal refined double characters.}
 We assume that the index set of those $z$'s
is $Z_{Max}\subset Z$.

We claim that for each $z\in Z_{Max}$ there is some linear
combination of acceptable $(\mu +1)$-tensor fields,

$$\Sum_{r\in R^z} a_r C^{r,i_1\dots i_{\alpha +1}}_{g}(\Omega_1,
\dots ,\Omega_p,\phi_1,\dots ,\phi_u),$$ where 
 each $C^{r,i_1\dots i_{\mu +1}}_{g}(\Omega_1, \dots
,\Omega_p,\phi_1,\dots ,\phi_u)$ has a $\mu$-double
character $\vec{L^z_1}$ and also the same set of factors 
$S_{*}\nabla^{(\nu)}R_{ijkl}$ as in $\vec{L}^z$ contain 
special free indices, so that:

\begin{equation}
\label{bengreen}
\begin{split}
& \Sum_{l\in L^z} a_l C^{l,i_1\dots i_\mu}_{g}
(\Omega_1,\dots ,\Omega_p,\phi_1,\dots
,\phi_u)\nabla_{i_1}\upsilon\dots\nabla_{i_\mu}\upsilon -
\\&\Sum_{r\in R^z} a_r X div_{i_{\mu +1}}
C^{r,i_1\dots i_{\mu +1}}_{g}(\Omega_1,\dots
,\Omega_p,\phi_1,\dots ,\phi_u)\nabla_{i_1}\upsilon\dots
\nabla_{i_\mu}\upsilon=
\\& \Sum_{t\in T_1} a_t
C^{t,i_1\dots i_\mu}_{g}(\Omega_1,\dots ,\Omega_p,,\phi_1,\dots
,\phi_u)\nabla_{i_1}\upsilon\dots \nabla_{i_\mu}\upsilon,
\end{split}
\end{equation}
modulo complete contractions of length $\ge\sigma +u+\mu +1$.
Here each
$$C^{t,i_1\dots i_\mu}_{g}(\Omega_1,\dots
,\Omega_p,\phi_1,\dots ,\phi_u)$$ is acceptable and is either
simply or doubly subsequent to $\vec{L}^z$.\footnote{Recall that
``simply subsequent'' means that the simple character of
$C^{t,i_1\dots i_\mu}_{g}$ is subsequent to $Simp(\vec{L}^z)$.}
\end{proposition}
 
 (See the first section in \cite{alexakis4} for a description 
 of the notions of {\it real length, acceptable tensor fields, 
simple character, 
refined double character, maximal refined double character, 
simply subsequent, strongly doubly subsequent}). 
The Proposition \ref{giade} is proven by a multiple induction on the parameters
$-n$ (the {\it weight} of the complete contractions appearing in (\ref{hypothese2})), 
$\sigma$ (the total number of factors in the form 
  $\nabla^{(m)}R_{ijkl}, S_{*}\nabla^{(\nu)}R_{ijkl}, \nabla^{(A)}\Omega_h$ among the partial 
  contractions in (\ref{hypothese2})),\footnote{The partial contractions in (\ref{hypothese2})
 are assumed to all have
  the same simple character--this implies that they all have the same number
   of factors $\nabla^{(m)}R_{ijkl}, S_{*}\nabla^{(\nu)}R_{ijkl}, \nabla^{(A)}\Omega_h$ respectively.}  
$\Phi$ (the number of factors $\nabla\phi_1,\dots,\nabla\phi_u$ appearing in (\ref{hypothese2})),
 $\Phi$ (the number of factors $\nabla\phi,\nabla\phi',\nabla\tilde{\phi}$
appearing appearing in (\ref{hypothese2})), 
  and $\sigma_1+\sigma_2$ (the total number of factors 
$\nabla^{(m)}R_{ijkl}, S_{*}\nabla^{(\nu)}R_{ijkl}$). Proposition \ref{giade} 
when $\Phi=0$ {\it coincides} with the ``Main algebraic Proposition'' in \cite{alexakis1} 
outlined above.\footnote{Similarly, the ``Main algebraic Propositions'' 3.1, 3.2 
in \cite{alexakis2} {\it coincide} with Proposition \ref{giade} above when $\Phi=1$.}

\section{Proof of the technical Lemmas from \cite{alexakis4}.}
\label{technical.lemmas}
\subsection{Re-statement of the technical Lemmas 4.6--4.9 from \cite{alexakis4}.}

We start by recalling 
a definition from \cite{alexakis4} that will be used frequently in the present paper:

\begin{definition}
\label{proextremovable} Consider any tensor field in the form
(\ref{form2}). We consider any set of indices,
$\{{}_{x_1},\dots,{}_{x_s}\}$ belonging to a factor $T$ (here $T$ is not in the form $\nabla\phi$). We
assume that these indices are neither free nor are contracting
against a factor $\nabla\phi_h$.

If the indices belong to a factor $T$ in the form
$\nabla^{(B)}\Omega_1$ then $\{{}_{x_1},\dots,{}_{x_s}\}$ are
removable provided $B\ge s+2$.

\par Now, we consider indices that belong to a factor
$\nabla^{(m)}R_{ijkl}$ (and are neither free nor are contracting
against a factor $\nabla\phi_h$). Any such index ${}_x$ which is a
derivative index will be removable. Furthermore, if $T$ has at
least two free derivative indices, then if neither of the indices
${}_i,{}_j$ are free then we will say one of ${}_i,{}_j$ is removable;
accordingly, if neither of ${}_k,{}_l$ is free then we will say
that one of ${}_k,{}_l$ is removable. Moreover, if $T$ has one free
derivative index then: if none of the indices ${}_i,{}_j$ are free
then we will say that one of the indices ${}_i,{}_j$ is removable;
on the other hand if one of the indices ${}_i,{}_j$ is also free
and none of the indices ${}_k,{}_l$ are free then we will say that
one of the indices ${}_k,{}_l$ is removable.

\par Now, we consider a set of indices $\{{}_{x_1},\dots,{}_{x_s}\}$ that belong to
a factor $T=S_{*}\nabla^{(\nu)}R_{ijkl}$ and are not special, and
are not free and are not contracting against any $\nabla\phi$. We
will say this set of indices is removable if $s\le\nu$.
Furthermore, if none of the indices ${}_k,{}_l$ are free and
$\nu>0$ and at least one of the other indices in $T$ is free, we
will say that one of the indices ${}_k,{}_l$ is removable.
\end{definition}
For the first two Lemmas, \ref{obote}, \ref{obote3} 
we will consider tensor fields in the form:

\begin{equation}
\label{form2'}
\begin{split}
&pcontr(\nabla^{(m_1)}R_{ijkl}\otimes\dots\otimes\nabla^{m_{\sigma_1}}
R_{ijkl}\otimes
\\&S_{*}\nabla^{(\nu_1)}R_{ijkl}\otimes\dots\otimes
S_{*}\nabla^{(\nu_t)} R_{ijkl}\otimes \nabla Y\otimes
\\& \nabla^{(b_1)}\Omega_1\otimes\dots\otimes
\nabla^{(b_p)}\Omega_p\otimes
\\& \nabla\phi_{z_1}\dots \otimes\nabla\phi_{z_w}\otimes\nabla
\phi'_{z_{w+1}}\otimes
\dots\otimes\nabla\phi'_{z_{w+d}}\otimes\dots \otimes
\nabla\tilde{\phi}_{z_{w+d+1}}\otimes\dots\otimes\nabla\tilde{\phi}_{z_{w+d+y}}).
\end{split}
\end{equation}
(Notice this is the same as the form (\ref{form2}), but for the
fact that we have inserted a factor $\nabla Y$ in the second line). Our claims are then the following:

\begin{lemma}
\label{obote} Assume an equation:

\begin{equation}
\label{guillemin3}
\begin{split}
& \Sum_{h\in H_2} a_h X_{*} div_{i_{1}}\dots X_{*} div_{i_{a_h}}
C^{h,i_{1}\dots i_{a_h}}_{g}(\Omega_1,\dots , \Omega_p,Y,
\phi_{1},\dots,\phi_{u'})=
\\& \Sum_{j\in J} a_j
C^j_{g}(\Omega_1,\dots ,\Omega_p, \phi_{1},\dots,\phi_{u'}),
\end{split}
\end{equation}
where all tensor fields have rank $a_h\ge \alpha$. All tensor
fields have a given $u$-simple character $\vec{\kappa}'_{simp}$,
for which $\sigma\ge 4$. Moreover, we assume that if we formally
treat the factor $\nabla Y$ as a factor $\nabla\phi_{u'+1}$ in the
above equation, then the inductive assumption of Proposition
\ref{giade} can be applied. (See subsection 3.1 in \cite{alexakis4} 
 for a strict discussion of the multi-parameter 
induction by which we prove Proposition \ref{giade}.)

The conclusion (under various assumptions which we will explain
below):  Denote by $H_{2,\alpha}$ the index set of tensor fields
with rank $\alpha$.

 We claim that there
is a linear combination of acceptable\footnote{``Acceptable'' in
the sense that each factor $\Omega_i$ is differentiated at least twice).}
 tensor fields, $\Sum_{d\in D} a_d C^{d,i_{1}\dots i_{\alpha +1}}_{g}(\Omega_1,\dots ,
\Omega_p,Y,\phi_{1},\dots, \phi_u)$, each with a simple character
$\vec{\kappa}'_{simp}$ so that:

\begin{equation}
\begin{split}
\label{vaskonik} &\Sum_{h\in H_{2,\alpha}} a_h C^{h,i_{1}\dots
i_{\alpha }}_{g}(\Omega_1,\dots ,\Omega_p, Y,\phi_{1},\dots,
\phi_{u'})\nabla_{i_{1}}\upsilon\dots \nabla_{i_{\alpha}}\upsilon-
\\& X_{*}div_{i_{\alpha+1}}\Sum_{d\in D} a_d
C^{d,i_{1}\dots i_{\alpha +1}}_{g}(\Omega_1,\dots ,\Omega_p,
Y,\phi_{1},\dots, \phi_{u'})\nabla_{i_{1}}\upsilon\dots
\nabla_{i_{\alpha}}\upsilon=
 \\&+\Sum_{t\in T}
a_t C^t_{g}(\Omega_1,\dots , \Omega_p
,Y,\phi_{1},\dots,\phi_{u'},\upsilon^{\alpha}).
\end{split}
\end{equation}
 The linear combination on the right hand side stands for a generic
linear combination of complete contractions in the form (\ref{form2'}) with a factor $\nabla
Y$ and with a simple character that is
 subsequent to $\vec{\kappa}'_{simp}$.
\newline

{\it The assumptions under which (\ref{vaskonik}) will hold:} The
assumption under which (\ref{vaskonik}) holds is that there should be no tensor fields 
 of rank $\alpha$ in (\ref{guillemin3}) which are ``bad''. Here ``bad'' means the following: 

If $\sigma_2=0$ in $\vec{\kappa}'_{simp}$ then a tensor field in the form (\ref{form2'}) 
is ``bad'' provided:
\begin{enumerate}
 \item{The factor $\nabla Y$ contains a free index.}
\item{If we formally erase the factor $\nabla Y$ (which contains 
a free index), then the resulting tensor field should have no removable indices,\footnote{Thus,
the tensor field should consist of factors $S_*R_{ijkl},\nabla^{(2)}\Omega_h$,
 and factors $\nabla^{(m)}_{r_1\dots r_m}R_{ijkl}$
 with all the indices ${}_{r_1},\dots ,{}_{r_m}$ 
contracting against factors $\nabla\phi_h$.} and no free 
indices.\footnote{I.e.~$\alpha=1$ in (\ref{guillemin3}).}
Moreover, any factors $S_{*}R_{ijkl}$ should be {\it simple}. }
\end{enumerate}

If $\sigma_2>0$ in $\vec{\kappa}'_{simp}$ then a tensor field in the form (\ref{form2'}) 
is ``bad'' provided:
\begin{enumerate}
 \item{The factor $\nabla Y$ should contain a free index.}
\item{If we formally erase the factor $\nabla Y$ (which contains 
a free index), then the resulting tensor field should have no removable indices, 
any factors $S_{*}R_{ijkl}$ should be {\it simple}, any
 factor $\nabla^{(2)}_{ab}\Omega_h$ should 
have at most one of the indices 
${}_a,{}_b$ free or contracting 
against  a factor $\nabla\phi_s$.}
\item{Any factor $\nabla^{(m)}R_{ijkl}$ can contain at most 
one (necessarily special, by virtue of 2.) free index. }
\end{enumerate}

  Furthermore, we claim that the proof of
 this Lemma will only rely on the inductive assumption of
 Proposition \ref{giade}. Moreover, we claim that if all the tensor
 fields indexed in $H_2$ (in (\ref{guillemin3})) do not have a free index in $\nabla Y$
 then we may assume that the tensor fields indexed in $D$ in (\ref{vaskonik}) have
 the same property.
\end{lemma}

\begin{lemma}
\label{obote3} We assume (\ref{guillemin3}), where $\sigma=3$. We
also assume that for each of the tensor fields in
$H_2^{\alpha,*}$\footnote{Recall from \cite{alexakis4} that $H_2^{\alpha,*}$
 is the index set of tensor fields of rank $\alpha$ in (\ref{guillemin3})
 with a free index in the factor $\nabla Y$.} there is at least
 one removable index. We then have two
claims:

Firstly, the conclusion of Lemma \ref{obote}
 holds in this setting also. Secondly, we  can write:
\begin{equation}
\begin{split}
\label{vaskonik} &\Sum_{h\in H_2} a_h Xdiv_{i_{1}}\dots
Xdiv_{i_{\alpha}}C^{h,i_{1}\dots i_{\alpha}}_{g}(\Omega_1,\dots
,\Omega_p, Y,\phi_{1},\dots, \phi_{u'})=
\\&\Sum_{q\in Q} a_q Xdiv_{i_{1}}\dots
Xdiv_{i_{a'}} C^{q,i_{1}\dots i_{a'}}_{g} (\Omega_1,\dots
,\Omega_p, Y,\phi_{1},\dots,\phi_{u'})
\\&+\Sum_{t\in T} a_t C^t_{g}(\Omega_1,\dots , \Omega_p
,Y,\phi_{1},\dots,\phi_{u'}),
\end{split}
\end{equation}
where the linear combination 
$\Sum_{q\in Q} a_q C^{q,i_{1}\dots
i_{a'}}_{g}$ stands
 for a generic linear combination of tensor fields
  in the form:
\begin{equation}
\label{form2''}
\begin{split}
&pcontr(\nabla^{(m_1)}R_{ijkl}\otimes\dots\otimes\nabla^{(m_{\sigma_1})}
R_{ijkl}\otimes
\\&S_{*}\nabla^{(\nu_1)}R_{ijkl}\otimes\dots\otimes
S_{*}\nabla^{(\nu_t)} R_{ijkl}\otimes \nabla^{(B)} Y\otimes
\\& \nabla^{(b_1)}\Omega_1\otimes\dots\otimes
\nabla^{(b_p)}\Omega_p\otimes
\\& \nabla\phi_{z_1}\dots \otimes\nabla\phi_{z_w}\otimes\nabla
\phi'_{z_{w+1}}\otimes
\dots\otimes\nabla\phi'_{z_{w+d}}\otimes\dots \otimes
\nabla\tilde{\phi}_{z_{w+d+1}}\otimes\dots\otimes\nabla\tilde{\phi}_{z_{w+d+y}}),
\end{split}
\end{equation}
 with $B\ge 2$, with a
 simple character $\vec{\kappa}'_{simp}$ and with each $a'\ge \alpha$. The acceptable complete
  contractions $C^t_{g}(\Omega_1,\dots , \Omega_p
 ,Y,\phi_{1},\dots,\phi_{u'})$ are simply subsequent to
 $\vec{\kappa}'_{simp}$. $Xdiv_i$ here means that $\nabla_i$ is
 not allowed to hit the factors $\nabla\phi_h$ (but it is
 allowed to hit $\nabla^{(B)}Y$).
\end{lemma}

\par For our next two Lemmas, we will be considering 
tensor fields in the general form:

\begin{equation}
\label{tavuk}
\begin{split}
&contr(\nabla^{(m_1)}R_{ijkl}\otimes\dots \otimes
\nabla^{(m_s)}R_{ijkl}\otimes
\\&S_{*}\nabla^{(\nu_1)}R_{ijkl}\otimes\dots\otimes
S_{*}\nabla^{(\nu_b)}R_{ijkl}\otimes \nabla^{(B,+)}_{r_1\dots
r_B}(\nabla_a\omega_1\nabla_b\omega_2-\nabla_b\omega_1\nabla_a\omega_1)
\\&\otimes \nabla^{(d_1)}\Omega_p\otimes\dots\otimes\nabla^{(d_p)}
\Omega_p \otimes \nabla\phi_1\otimes\dots\otimes
\nabla\phi_u);
\end{split}
\end{equation}
here $\nabla^{(B,+)}_{r_1\dots r_B}(\dots )$ stands for the
sublinear combination in $\nabla^{(B)}_{r_1\dots r_B}(\dots )$
where each $\nabla$ is not allowed to hit the factor
$\nabla\omega_2$.

\begin{lemma}
\label{vanderbi} Consider a linear combination of partial
contractions,
$$\Sum_{x\in X} a_x C^{x,i_1\dots
i_a}_{g}(\Omega_1,\dots
,\Omega_p,[\omega_1,\omega_2],\phi_1,\dots,\phi_{u'}),$$ where each
of the tensor fields $C^{x,i_1\dots
i_a}_{g}$ is in the form (\ref{tavuk}) with $B=0$ (and
is antisymmetric in the factors
$\nabla_a\omega_1,\nabla_b\omega_2$ by definition), with rank
$a\ge\alpha$ and real length $\sigma\ge 4$.\footnote{ Recall
that in the definition of ``real length'' in this setting,
we count each factor $\nabla^{(m)}R,\S_{*}\nabla^{(\nu)}R, \nabla^{(B)}\Omega_x$ once,
the two factors $\omega_1,\omega_2$ for one, and the 
factors $\nabla\phi,\nabla\phi',\nabla\tilde{\phi}$ nor nothing.} We assume
that all these tensor fields have a given simple character which we
denote by $\vec{\kappa}'_{simp}$ (we use $u'$ instead of $u$ to
stress that this Lemma holds in generality). We assume an
equation:

\begin{equation}
\label{para3enh}
\begin{split}
& \Sum_{x\in X} a_x X_{*}div_{i_1}\dots
X_{*}div_{i_a}C^{x,i_1\dots i_a}_{g}(\Omega_1,\dots
,\Omega_p,[\omega_1,\omega_2],\phi_1,\dots,\phi_u)+
\\& \Sum_{j\in J} a_j C^j_{g}(\Omega_1,\dots
,\Omega_p,[\omega_1,\omega_2],\phi_1,\dots,\phi_u)=0,
\end{split}
\end{equation}
where $X_{*}div_i$ stands for the sublinear combination in
$Xdiv_i$ where $\nabla_i$ is in addition not allowed to hit the
factors $\nabla\omega_1,\nabla\omega_2$. The contractions $C^j$
here are simply subsequent to $\vec{\kappa}'_{simp}$. We assume
that if we formally treat the factors
$\nabla\omega_1,\nabla\omega_2$ as factors
$\nabla\phi_{u+1},\nabla\phi_{u+2}$ (disregarding whether they are
contracting against special indices) in the above, then the
inductive assumption of Proposition \ref{giade} applies.

\par The conclusion we will draw (under various hypotheses that
we will explain below) is that we can write:

\begin{equation}
\label{para3enh2} \begin{split} &\Sum_{x\in X} a_x
X_{+}div_{i_1}\dots X_{+}div_{i_a}C^{x,i_1\dots
i_a}_{g}(\Omega_1,\dots
,\Omega_p,[\omega_1,\omega_2],\phi_1,\dots,\phi_u)=
\\&\Sum_{x\in X'} a_x
X_{+}div_{i_1}\dots X_{+}div_{i_a}C^{x,i_1\dots
i_a}_{g}(\Omega_1,\dots
,\Omega_p,[\omega_1,\omega_2],\phi_1,\dots,\phi_u)+
\\& \Sum_{j\in
J} a_j C^j_{g}(\Omega_1,\dots
,\Omega_p,[\omega_1,\omega_2],\phi_1,\dots,\phi_u)=0,
\end{split}
\end{equation}
where the tensor fields indexed in $X'$ on the right hand side are
in the form (\ref{tavuk}) with $B>0$. All the other sublinear
combinations are as above. We recall from \cite{alexakis4} that $X_{+}div_i$ stands for
the sublinear combination in $Xdiv_i$ where $\nabla_i$ is in
addition not allowed to hit the factor $\nabla\omega_2$ (it is
allowed to hit the factor $\nabla^{(B)}\omega_1$).
\newline

{\it Assumptions needed for (\ref{para3enh2}):} We claim
(\ref{para3enh2}) under certain assumptions on the $\alpha$-tensor
fields in (\ref{para3enh}) which have rank $\alpha$ and have a
free index in one of the factors $\nabla\omega_1,\nabla\omega_2$
(say to $\nabla\omega_1$ wlog)--we denote the index set of those
tensor fields by $X^{\alpha,*}\subset X$.

 The assumption we need in order for the claim to hold is that no
tensor field indexed in $X^{\alpha,*}$ should be ``bad''. A  tensor field 
is ``bad'' if it has the property that 
when we erase the expression 
$\nabla_{[a}\omega_1\nabla_{b]}\omega_2$ (and make 
the index that contrated against ${}_b$ into a free index) then the resulting tensor field 
will have no removable indices, and all factors $S_{*}R_{ijkl}$ will be simple.
\end{lemma}

\begin{lemma}
\label{vanderbi3} We assume (\ref{para3enh}), where now the tensor 
fields have length $\sigma=3$.
We also assume that for each of the tensor fields indexed in  $X$, there
is a removable index in each of the real factors. We then claim that the conclusion of
Lemma \ref{vanderbi} is still true in this setting.
\end{lemma}

For the most part, the remainder of this 
paper is devoted to proving the above Lemmas. However, 
we first state and prove some further technical claims, one of which 
appeared as Lemma 5.1 in \cite{alexakis4}.\footnote{Its 
proof was also deferred to the present paper.}

\subsection{Two more technical Lemmas.}

We claim an analogue of Lemma 4.10 in \cite{alexakis4} can
be derived when we have tensor fields with a given simple
character $\vec{\kappa}_{simp}$, and where rather than having {\it
one} additional factor $\nabla \phi_{u+1}$ (which is not encoded
in the simple character $\vec{\kappa}_{simp}$), we have {\it two
additional factors}  $\nabla_{a}\phi_{u+1},\nabla_{b}\phi_{u+2}$.

\begin{lemma}
\label{addition2}
 Consider a linear combination of acceptable tensor fields in the form
(\ref{form2}) with a given $u$-simple character
 $\vec{\kappa}_{simp}$: \\$\sum_{l\in L} a_l
C^{l,i_1\dots
i_\beta}_g(\Omega_1,\dots,\Omega_p,\phi_1,\dots,\phi_u)$.
 Assume that the minimum rank
among those tensor fields above is $\alpha\ge 2$.  Assume an
equation:

\begin{equation}
\begin{split}
\label{alexandraoik} & \sum_{l\in L} a_l X_{*}div_{i_3}\dots
X_{*}div_{i_\beta} C^{l,i_1\dots
i_\beta}_g(\Omega_1,\dots,\Omega_p,\phi_1,\dots,\phi_u)
\nabla_{i_1}\phi_{u+1}\nabla_{i_2}\phi_{u+2}+
\\&\sum_{j\in J} a_j C^j_g(\Omega_1,\dots,\Omega_p,\phi_1,\dots,\phi_u)=0
\end{split}
\end{equation}
(here $X_{*}div_i$ means that $\nabla^i$ is in addition not
allowed to hit the factors $\nabla\phi_{u+1}$, $\nabla\phi_{u+2}$).
We also assume that if we formally treat the factors
 $\nabla\phi_{u+1}$, $\nabla\phi_{u+2}$ as factors
 $\nabla\phi_{u+1}$, $\nabla\phi_{u+2}$ then (\ref{alexandraoik}) falls
 under the inductive assumption of Proposition \ref{giade} (with respect to the parameters
$(n,\sigma,\Phi,u)$). Denote by $L^\alpha\subset L$ the index set
of terms with rank $\alpha$. We additionally assume that none of the tensor fields 
$C^{l,i_1\dots
i_\beta}_g(\Omega_1,\dots,\Omega_p,\phi_1,\dots,\phi_u)$ are ``forbidden'', 
in the sense defined above Proposition 2.1 in \cite{alexakis4}.

\par We then claim that there exists a linear combination
 of $(\alpha+1)$-tensor fields with a
$u$-simple character $\vec{\kappa}_{simp}$ (indexed in $Y$ below)
so that:

\begin{equation}
\label{palataki}
\begin{split}
 & \sum_{l\in L^\alpha} a_l
C^{l,i_1\dots
i_\alpha}_g(\Omega_1,\dots,\Omega_p,\phi_1,\dots,\phi_u)
\nabla_{i_1}\phi_{u+1}\nabla_{i_2}\phi_{u+2}
\nabla_{i_3}\upsilon\dots\nabla_{i_\alpha}\upsilon
\\+&X_{*}div_{i_{\alpha+1}}\sum_{y\in Y} a_y
C^{l,i_1\dots
i_{\alpha+1}}_g(\Omega_1,\dots,\Omega_p,\phi_1,\dots,\phi_u)
\nabla_{i_1}\phi_{u+1}\nabla_{i_2}\phi_{u+2}
\nabla_{i_3}\upsilon\dots\nabla_{i_\alpha}\upsilon
\\&+\sum_{j\in J} a_j C^{j,i_1\dots i_\alpha}_g
(\Omega_1,\dots,\Omega_p,\phi_1,\dots,\phi_u)
\nabla_{i_1}\phi_{u+1}\nabla_{i_2}\phi_{u+2}
\nabla_{i_3}\upsilon\dots\nabla_{i_\alpha}\upsilon.
\end{split}
\end{equation}

\par Furthermore, we also claim that we can write:

\begin{equation}
\begin{split}
\label{alexandraoik'} & \sum_{l\in L} a_l Xdiv_{i_3}\dots
Xdiv_{i_\beta} C^{l,i_1\dots
i_\beta}_g(\Omega_1,\dots,\Omega_p,\phi_1,\dots,\phi_u)
\nabla_{i_1}\phi_{u+1}\nabla_{i_2}\phi_{u+2}=
\\&\sum_{j\in J} a_j C^j_g(\Omega_1,\dots,\Omega_p,\phi_1,\dots,\phi_u)+
\\&\sum_{q\in Q_1} a_q Xdiv_{i_3}\dots
Xdiv_{i_\alpha} C^{q,i_1\dots
i_\alpha}_g(\Omega_1,\dots,\Omega_p,\phi_1,\dots,\phi_{u+2})
\nabla_{i_1}\phi_{u+1}\nabla_{i_2}\phi_{u+2}+
\\&\sum_{q\in Q_2} a_q Xdiv_{i_3}\dots
Xdiv_{i_\alpha} C^{q,i_1\dots
i_\alpha}_g(\Omega_1,\dots,\Omega_p,\phi_1,\dots,\phi_{u+2}),
\end{split}
\end{equation}
where the tensor fields indexed in $Q_1$ are acceptable
with a u-simple character $\vec{\kappa}_{simp}$ and with a
factor $\nabla^{(2)}\phi_{u+1}$ and a factor $\nabla\phi_{u+2}$.
The tensor fields indexed in $Q_2$ are acceptable
with a u-simple character $\vec{\kappa}_{simp}$ and with a
factor $\nabla^{(2)}\phi_{u+2}$ and a factor $\nabla\phi_{u+1}$.
\end{lemma}

{\it Proof of Lemma \ref{addition2}:} We may divide the index set
$L^\alpha$ into subsets $L^\alpha_I,L^\alpha_{II}$ according to
whether the two factors $\nabla\phi_{u+1},\nabla\phi_{u+2}$
 are contracting against the same factor
or not--we will then prove our claim for those two index sets
separately. Our claim for the index set $L^\alpha_{II}$ follows by
a straightforward adaptation of the proof of  Lemma
4.10 in \cite{alexakis1}. (Notice that the forbidden cases of the present
 Lemma are exactly in correspondence  with the forbidden cases of that Lemma).
Therefore, we now prove our claim for the index
set $L^\alpha_{I}$:

We denote by $L_I\subset L, J_I\subset J$ the index set of terms
for which the two factors $\nabla\phi_{u+1},\nabla\phi_{u+2}$ are
contracting against the same factor. It then follows that
(\ref{alexandraoik}) holds with the index sets $L,J$ replaced by
$L_I,J_I$--denote the resulting new equation by
New[(\ref{alexandraoik})]. Now, for each tensor field $C^{l,i_1\dots i_\beta}_g$ and each
complete contraction $C^j_g$,
 we let $Sym[C^{l,i_1\dots i_\beta}_g]$, $Sym[C^{l,i_1\dots i_\beta}_g]$, 
$AntSym[C^j_g], AntSym[C^j_g]$
stand for the tensor field/complete contraction  that arises from 
$C^{l,i_1\dots i_\beta}_g,C^j_g$ by
 symmetrizing (resp.~anti-symmetrizing) the indices ${}_a,{}_b$ in the two factors
$\nabla_a\phi_{u+1},\nabla_b\phi_{u+2}$. 
We accordingly derive two new equations from
New[(\ref{alexandraoik})], which we denote by
New[(\ref{alexandraoik})$]_{Sym}$ and
New[(\ref{alexandraoik})$]_{AntSym}$.

\par We will then prove the claim separately for the tensor fields
in the sublinear combination $\sum_{l\in L_I^\alpha} a_l
Sym[C]^{l,i_1\dots i_\alpha}_g$ and the tensor fields in the
sublinear combination
 $\sum_{l\in L_I^\alpha} a_l AntSym[C]^{l,i_1\dots i_\alpha}_g$.

\par The claim (\ref{palataki}) for the sublinear combination
 $\sum_{l\in L_I^\alpha} a_l AntSym[C]^{l,i_1\dots i_\alpha}_g$
follows directly from the arguments in the proof of Lemma
\ref{vanderbi}. Therefore it suffices to show our claim for the
sublinear combination $\sum_{l\in L_I^\alpha} a_l
Sym[C]^{l,i_1\dots i_\alpha}_g$.

\par We prove this claim as follows: We divide the index
set $L_I^\alpha$ according to the {\it form} of the factor against
which the two factors $\nabla\phi_{u+1},\nabla\phi_{u+2}$ are
contracting: List out the non-generic factors in
$\vec{\kappa}_{simp}$,\footnote{Recall from the introduction in \cite{alexakis4}
that the non-generic factors in $\vec{\kappa}_{simp}$ are all 
the factors in the form $\nabla^{(A)}\Omega_h, S_{*}\nabla^{(\nu)}R_{ijkl}$, 
and also all the factors $\nabla^{(m)}R_{ijkl}$ that contract 
against at least one factor $\nabla\phi_s$.}
 $\{T_1,\dots, T_a\}$. Then, for each $k\le a$ we let $L^\alpha_{I,k}$ stand for the
index set of terms for which the factors
$\nabla\phi_{u+1},\nabla\phi_{u+2}$ are contracting against the
factor $T_k$. We also let $L^\alpha_{I,a+1}$ stand for the index
set of terms for which the factors 
$\nabla\phi_{u+1},\nabla\phi_{u+2}$ are contracting against a
generic factor $\nabla^{(m)}R_{ijkl}$. We will prove our claim for
each of the sublinear combinations $\sum_{l\in L^\alpha_{I,a+1}}
a_l Sym[C]^{l,i_1\dots i_\alpha}_g$ separately.

\par We firstly observe that for each $k\le a+1$,
 we may obtain a new true equation from
(\ref{alexandraoik}) by replacing $L$ by $L_{I,a+1}$--denote the
resulting equation by (\ref{alexandraoik}$)_{I,Sym,k}$. Therefore,
for each $k\le a+1$ for which $T_k$ is in the form
$\nabla^{(p)}\Omega_h$, our claim follows straightforwardly by
applying Corollary 1 from \cite{alexakis1}.\footnote{There is no danger of
falling under a ``forbidden case'', since 
we started with tensor fields which were not forbidden.}

\par Now, we consider the case where the factor
 $T_k$ is in the form $S_{*}\nabla^{(\nu)}R_{abcd}$:
In that case we denote by $L_{I,k,\sharp}$
 the index set of terms for
which one of the factors $\nabla\phi_{u+1},\nabla\phi_{u+2}$ is
contracting against a special index in $T_k$. In particular, we will
let $L^\alpha_{I,k,\sharp}\subset L_{I,k,\sharp}$ stand for the index
set of terms with rank $\alpha$.
 We will then show
two equations:

Firstly, that there exists a linear combination of tensor fields
as claimed in (\ref{palataki}) so that:
\begin{equation}
\label{bestdays}
\begin{split}
&\sum_{l\in L^\alpha_{I,k,\sharp}} a_l Sym[C]^{l,i_1\dots
i_\alpha}_g(\Omega_1,\dots,\Omega_p,\phi_1,\dots,\phi_u)
 \nabla_{i_1}\phi_{u+1}\nabla_{i_2}\omega\nabla_{i_3}\upsilon\dots\nabla_{i_\alpha}\upsilon-
\\&\sum_{y\in Y} a_y Xdiv_{i_{\alpha+1}}
C^{y,i_1\dots
i_{\alpha+1}}_g(\Omega_1,\dots,\Omega_p,\phi_1,\dots,\phi_u)
 \nabla_{i_1}\phi_{u+1}\nabla_{i_2}\omega\nabla_{i_3}
\upsilon\dots\nabla_{i_\alpha}\upsilon
\\&=\sum_{l\in L^\alpha_{OK}} a_l Xdiv_{i_{\alpha+1}}  
C^{l,i_1\dots i_\alpha i_{\alpha+1}}_g\nabla_{i_1}
\phi_{u+1}\nabla_{i_2}\omega\nabla_{i_3}
\upsilon\dots\nabla_{i_\alpha}\upsilon+
\\&\sum_{j\in J} a_j C^{j,i_1\dots i_\alpha}_g(\Omega_1,\dots,\Omega_p,\phi_1,\dots,\phi_u)
 \nabla_{i_1}\phi_{u+1}\nabla_{i_2}\omega\nabla_{i_3}\upsilon\dots\nabla_{i_\alpha}\upsilon,
\end{split}
\end{equation}
where the tensor fields in $L^\alpha_{OK}$  have all the properties of the
terms in $L_{I,k}$, rank  $\alpha$
and furthermore none of the factors $\nabla\phi_{u+1},
\nabla\phi_{u+2}$ are contracting against a special index.

Then (under the assumption that $L^\alpha_{I,k,\sharp}=\emptyset$)
we claim that we can write:
\begin{equation}
\label{bestdays2}
\begin{split}
&\sum_{l\in L_{I,k,\sharp}} a_l Xdiv_{i_3}\dots Xdiv_{i_\beta}
Sym[C]^{l,i_1\dots
i_\beta}_g(\Omega_1,\dots,\Omega_p,\phi_1,\dots,\phi_u)
 \nabla_{i_1}\phi_{u+1}\nabla_{i_2}\phi_{u+2}
\\&=\sum_{l\in L_{I,k,OK}} a_l Xdiv_{i_3}\dots Xdiv_{i_\beta}
Sym[C]^{l,i_1\dots
i_\beta}_g(\Omega_1,\dots,\Omega_p,\phi_1,\dots,\phi_u)
 \\&\nabla_{i_1}\phi_{u+1}\nabla_{i_2}\phi_{u+2}
+\sum_{j\in J} a_j Sym[C]^{j,i_1 i_2}_g
(\Omega_1,\dots,\Omega_p,\phi_1,\dots,\phi_u)
 \nabla_{i_1}\phi_{u+1}\nabla_{i_2}\phi_{u+2},
\end{split}
\end{equation}
where the tensor fields in $L_{I,k,OK}$ have all the properties of the
terms in $L_{I,k}$, but they additionally have rank $\ge \alpha+1$
and furthermore none of the factors $\nabla\phi_{u+1},
\nabla\phi_{u+2}$ are contracting against a special index.

\par If we can show the above two equations, then we are
 reduced to showing our claim under the additional assumption
that no tensor field indexed in $L$ in Sym(\ref{alexandraoik}) has
any factor $\nabla\phi_{u+1}$, $\nabla\phi_{u+1}$ contracting
against a special index in $T_k$. Under that assumption, we may
additionally assume that none of the complete contractions indexed
in $J$ in (\ref{alexandraoik}) have that property.\footnote{This 
can be derived by repeating the {\it proof} of 
(\ref{bestdays}), (\ref{bestdays2}).} Therefore,
we may then {\it erase} the factor $\nabla\phi_{u+1}$ from all the
complete contractions and tensor fields in
(\ref{alexandraoik}$)_k$ by virtue of the operation $Erase$, 
introduced in the Appendix of \cite{alexakis1}--our claim then follows by applying
Corollary 1 from \cite{alexakis1} to the resulting equation and then
re-introducing the erased factor
 $\nabla\phi_{u+1}$.

{\it Outline of the proof of (\ref{bestdays}), (\ref{bestdays2}):}
Firstly we prove (\ref{bestdays}): Suppose wlog $T_k$ is
contracting against $\nabla\tilde{\phi}_1$ and
$\nabla\phi'_2,\dots,\nabla\phi'_h$; then replace the two factors
$\nabla_a\phi_1,\nabla_b\phi_{u+1}$ by $g_{ab}$ and then apply
$Ricto \Omega_{p+1}$,\footnote{See the relevant Lemma in the Appendix of \cite{alexakis1}.}
 (obtaining a new true equation) an then apply the eraser to  the resulting true equation.
We then apply Corollary 1 from \cite{alexakis1} to the resulting
equation,\footnote{Since the factor $\nabla\phi_{u+2}$ survives this operation, 
and since we started out with terms that were not ``forbidden'', there is
no danger of falling under a ``forbidden
case'' of Corollary 1 from \cite{alexakis1}.}
and finally we replace the factor $\nabla^{(b)}_{r_1\dots r_b}
\Omega_{p+1}$
 by an expression
$$S_{*}\nabla^{(b+h-1)}_{y_2\dots y_hr_1\dots r_{b-1}}
R_{ijkr_b}\nabla^i\tilde{\phi}_1\nabla^j\phi_{u+2}\nabla^k\phi_{u+1}
\nabla^{y_2}\phi'_2\dots \nabla^{y_h}\phi'_h.$$
 As in the proof of Lemma
 4.10 in \cite{alexakis4}, we derive our claim. Then, (\ref{bestdays2})
 is proven by iteratively applying this
step and making $\nabla\upsilon$'s into $Xdiv$'s at each stage.
\newline

\par We analogously show our claim when the factor $T_k$
is in the form $\nabla^{(m)}R_{ijkl}$: In that case we denote by
$L_{I,k,\sharp}$
 the index set of terms for
which {\it both} the factors $\nabla\phi_{u+1},\nabla\phi_{u+2}$
are contracting against a special index in $T_k$. We will then show
two equations:

Firstly, that there exists a linear combination of tensor fields
as claimed in (\ref{palataki}) so that:
\begin{equation}
\label{bestdays'}
\begin{split}
&\sum_{l\in L^\alpha_{I,k,\sharp}} a_l Sym[C]^{l,i_1\dots
i_\alpha}_g(\Omega_1,\dots,\Omega_p,\phi_1,\dots,\phi_u)
 \nabla_{i_1}\phi_{u+1}\nabla_{i_2}\omega\nabla_{i_3}\upsilon\dots\nabla_{i_\alpha}\upsilon-
\\&\sum_{y\in Y} a_y Xdiv_{i_{\alpha+1}}
C^{y,i_1\dots
i_{\alpha+1}}_g(\Omega_1,\dots,\Omega_p,\phi_1,\dots,\phi_u)
 \nabla_{i_1}\phi_{u+1}\nabla_{i_2}\omega\nabla_{i_3}\upsilon\dots\nabla_{i_\alpha}\upsilon=
\\&\sum_{l\in L^\alpha_{OK}} a_l Xdiv_{i_{\alpha+1}}  
C^{l,i_1\dots i_\alpha i_{\alpha+1}}_g\nabla_{i_1}\phi_{u+1}\nabla_{i_2}\omega\nabla_{i_3}
\upsilon\dots\nabla_{i_\alpha}\upsilon+
\\&\sum_{j\in J} a_j C^{j,i_1\dots i_\alpha}_g(\Omega_1,\dots,\Omega_p,\phi_1,\dots,\phi_u)
 \nabla_{i_1}\phi_{u+1}\nabla_{i_2}\omega\nabla_{i_3}\upsilon\dots\nabla_{i_\alpha}\upsilon,
\end{split}
\end{equation}
where the tensor fields in $L^\alpha_{OK}$ have all the properties of the
terms in $L_{I,k}$, but they additionally have rank $\alpha$
and furthermore one
 of the factors $\nabla\phi_{u+1}, \nabla\phi_{u+2}$
does not contract against a  special index.
Then (under the assumption that $L^\alpha_{I,k,\sharp}=\emptyset$)
we denote by $L_{I,k,\sharp}$ the sublinear combination of terms
in $L_{I,k}$ with {\it both} factors $\nabla\phi_{u+1}$ or
$\nabla\phi_{u+1}$ contracting against a special index in $T_k$.
We claim that we can write:

\begin{equation}
\label{bestdays2'}
\begin{split}
&\sum_{l\in L_{I,k,\sharp}} a_l Xdiv_{i_3}\dots Xdiv_{i_\beta}
Sym[C]^{l,i_1\dots
i_\beta}_g(\Omega_1,\dots,\Omega_p,\phi_1,\dots,\phi_u)
 \nabla_{i_1}\phi_{u+1}\\&\nabla_{i_2}\phi_{u+2}=
\sum_{l\in L_{I,k,OK}} a_l Xdiv_{i_3}\dots Xdiv_{i_\beta}
Sym[C]^{l,i_1\dots
i_\beta}_g(\Omega_1,\dots,\Omega_p,\phi_1,\dots,\phi_u)
\\& \nabla_{i_1}\phi_{u+1}\nabla_{i_2}\phi_{u+2}+
\sum_{j\in J} a_j Sym[C]^{j,i_1 i_2}_g(\Omega_1,\dots,\Omega_p,\phi_1,\dots,\phi_u)
 \nabla_{i_1}\phi_{u+1}\nabla_{i_2}\phi_{u+2},
\end{split}
\end{equation}
where the tensor fields in $L_{I,k,OK}$ have all the properties of the
terms in $L_{I,k}$, but they additionally have rank $\ge \alpha+1$
and furthermore one
 of the factors $\nabla\phi_{u+1}, \nabla\phi_{u+2}$
does not contract against a  special index.

\par If we can show the above two equations, then we are
 reduced to showing our claim under the additional assumption
that no tensor field indexed in $L$ in Sym(\ref{alexandraoik}) has
the two factors $\nabla\phi_{u+1}, \nabla\phi_{u+2}$,
 contracting against a special index in $T_k$.
Under that assumption, we may additionally assume that none of the
complete contractions indexed in $J$ in
 (\ref{alexandraoik}) have that property. Therefore,
we may then {\it erase} the factor $\nabla\phi_{u+1}$ from all the
complete contractions and tensor fields in
(\ref{alexandraoik}$)_k$--our claim then follows by applying Lemma
4.10 in \cite{alexakis4} to the resulting equation\footnote{Notice that
there is no danger of falling under a ``forbidden case'' of that Lemma,
 since there will be a non-simple factor
$S_{*}\nabla^{(\nu)}R_{ijkl}$, by virtue of the factor
$\nabla\phi_{u+2}$.} and then re-introducing the erased factor
$\nabla\phi_{u+1}$.

{\it Outline of the proof of (\ref{bestdays'}),
(\ref{bestdays2'}):} Firstly we
 prove  (\ref{bestdays'}). Suppose wlog $T_k$ is contracting
against $\nabla\phi_1,\dots,\nabla\phi_h$ (possibly with $h=0$);
 then replace the two factors
$\nabla_a\phi_1,\nabla_b\phi_{u+1}$ by $g_{ab}$ and then apply
$Ricto \Omega_{p+1}$
 (obtaining a new true equation) an then apply the eraser to
the factors $\nabla\phi_1,\dots,\nabla\phi_h$ in the resulting
true equation. Then  (apart from the cases, discussed below, 
where the above operation may
lead to a ``forbidden case'' of Corollary 1 in \cite{alexakis4}), 
we apply Corollary 1 frm \cite{alexakis4} to the
resulting equation,
and finally we replace the factor
$\nabla^{(b)}_{r_1\dots r_b}\Omega_{p+1}$ by an expression
$$\nabla^{(b+h)}_{s_1\dots s_hr_1\dots
r_{b-2}}R_{ir_{b-1}kr_b}\nabla^i\phi_{u+1}\nabla^k\phi_{u+2}\nabla^{s_1}\phi_1\dots\nabla^{s_h}\phi_h.$$
 As in the proof of Lemma 4.10 in \cite{alexakis4}, we derive our claim. Then, (\ref{bestdays'})
 is proven by iteratively applying this
step and making $\nabla\upsilon$'s into $Xdiv$'s at each stage
(again,  provided we never encounter  ``forbidden cases'').
If we do encounter forbidden cases, then our claims
follow by just making the factors $\nabla\phi_{u+1},\nabla\phi_{u+2}$
into $Xdiv$'s  and then applying Corollary
1 in \cite{alexakis1} to the resulting equation (the 
resulting equation is not forbidden, since it will
contain a factor $\nabla^{(m)}R_{ijkl}$ with two free indices), and in the end
re-naming two factors $\nabla\upsilon$ into
$\nabla\phi_{u+1},\nabla\phi_{u+2}$. $\Box$
\newline

{\bf A Further Generalization: Proof of Lemma 5.1 from \cite{alexakis4}.}
\newline

We remark that on a few occasions later in this series of papers we will be using a 
generalized version of the Lemma \ref{addition2}. 
The generalized version asserts that the claim of Lemma \ref{addition2} remains true, 
for the general case where rather than one or two ``additional''
 factors $\nabla\phi_{u+1},\nabla\phi_{u+2}$ we have $\beta\ge 3$ ``additional'' 
factors $\nabla\phi_{u+1},\dots,\nabla\phi_{u+\beta}$. Moreover, 
in that case there are no ``forbidden cases''.

\begin{lemma}
\label{additiongen}
Let $\sum_{l\in L_1} a_l C^{l,i_1\dots i_\mu,i_{\mu+1}\dots i_{\mu+\beta}}_g
(\Omega_1,\dots,\Omega_p,\phi_1,\dots,\phi_u)$,

$\sum_{l\in L_2} a_l C^{l,i_1\dots i_{b_l},i_{b_l+1}\dots i_{b_l+\beta}}_g
(\Omega_1,\dots,\Omega_p,\phi_1,\dots,\phi_u)$ 
 stand for two linear combinations of acceptable 
tensor fields in the form (\ref{form2}), with a $u$-simple 
character $\vec{\kappa}_{simp}$. We assume that the terms 
indexed in $L_1$ have rank $\mu+\beta$, while 
the ones indexed in $L_2$ have rank greater than $\mu+\beta$. 

Assume an equation:

\begin{equation}
\begin{split}
\label{paragon} & \sum_{l\in L_1} a_l Xdiv_{i_{\beta+1}}\dots
Xdiv_{i_{\mu+\beta}} C^{l,i_1\dots
i_{\mu+\beta}}_g(\Omega_1,\dots,\Omega_p,\phi_1,\dots,\phi_u)
\nabla_{i_1}\phi_{u+1}\dots \nabla_{i_\beta}\phi_{u+\beta}
\\&+\sum_{l\in L_2} a_l Xdiv_{i_{\beta+1}}\dots
Xdiv_{i_{b_l}} C^{l,i_1\dots
i_{b_l+\beta}}_g(\Omega_1,\dots,\Omega_p,\phi_1,\dots,\phi_u)
\nabla_{i_1}\phi_{u+1}\dots \nabla_{i_\beta}\phi_{u+\beta}
\\&+\sum_{j\in J} a_j C^{j}_g(\Omega_1,\dots,\Omega_p,\phi_1,\dots,\phi_{u+\beta})=0,
\end{split}
\end{equation}
modulo terms of length $\ge\sigma+u+\beta+1$.
Furthermore, we assume that the above equation 
falls under the inductive assumption of Proposition 2.1 in \cite{alexakis4} 
(with regard to the parameters weight, $\sigma,\Phi,p$). 
 We are not excluding any  ``forbidden cases''.

We claim that there exists a linear combination of $(\mu+\beta+1)$-tensor fields 
in the form (\ref{form2}) with $u$-simple character $\vec{\kappa}_{simp}$ and 
length $\sigma+u$ (indexed in $H$ below) such that:

\begin{equation}
\begin{split}
\label{paragon.conc} & \sum_{l\in L_1} a_l C^{l,i_1\dots
i_{\mu+\beta}}_g(\Omega_1,\dots,\Omega_p,\phi_1,\dots,\phi_u)
\nabla_{i_{1}}\phi_{u+1}\dots \nabla_{i_\beta}
\phi_{u+\beta}\nabla_{i_{\beta+1}}\upsilon\dots\nabla_{i_{\beta+\mu}}\upsilon
\\&+\sum_{h\in H} a_h Xdiv_{i_{\mu+\beta+1}} C^{l,i_1\dots
i_{\mu+\beta+1}}_g(\Omega_1,\dots,\Omega_p,\phi_1,\dots,\phi_u)
\nabla_{i_1}\phi_{u+1}\dots \nabla_{i_\beta}
\phi_{\beta+1}\\&\nabla_{i_1}\upsilon\dots\nabla_{i_{\beta+\mu}}\upsilon+
\sum_{j\in J} a_j C^j_g(\Omega_1,\dots,
\Omega_p,\phi_1,\dots,\phi_{u+\beta},\upsilon^\mu)=0,
\end{split}
\end{equation}
modulo terms of length $\ge\sigma+u+\beta+\mu+1$. The terms indexed in 
$J$ here are $u$-simply subsequent to $\vec{\kappa}_{simp}$.
\end{lemma}

{\bf Proof of Lemma \ref{additiongen}:} The proof of the above is 
a straightforward adaptation of the proof of Lemma \ref{addition2}, 
{\it except} for the cases where 
 the tensor fields 
$C^{l,i_1\dots i_\mu,i_{\mu+1}\dots i_{\mu+\beta}}_g$ are ``bad'',
 where ``bad'' in this case would mean that all
 factors are in the form $R_{ijkl}$, $S_*R_{ijkl}$, 
$\nabla^{(2)}\Omega_h$,\footnote{Notice that by 
weight considerations, if this property holds for one of the 
terms $C^{l,i_1\dots i_\mu,i_{\mu+1}\dots i_{\mu+\beta}}_g$, 
then it will hold for all of them.} 
and in addition each factor $\nabla^{(2)}\Omega_h$  contracts 
against at most one factor $\nabla\phi_h, 1\le h\le u+\beta$. 
So we now focus on that case:

Let us observe that by weight considerations,
 {\it all} tensor fields 
in (\ref{alexandraoik}) must now have rank $\mu$.

We recall that this special proof applies only in the case where 
there are special free indices in factors $S_{*}R_{ijkl}$ 
among the tensor fields of minimum rank in (\ref{alexandraoik}). 
(If there were no such terms, then the 
regular proof of Lemma \ref{addition2} applies).
We distinguish three cases: Either $p>0$,
 or $p=0, \sigma_1>0$ or $p=\sigma_1=0$ 
and $\sigma_2>0$. We will prove the above by an induction on 
the parameters $(weight)$, $\sigma$: Suppose that the weight of the terms in (\ref{paragon}) 
is $-K$ and the real length is $\sigma\ge 3$. We assume that 
the Lemma holds when the equation (\ref{paragon}) 
consists of terms with weight $-K', K'<K$, 
or of terms with weight $-K$ and real length $\sigma', 3\le \sigma'<\sigma$. 
\newline

{\bf The case $p>0$:} We first consider the $\mu$-tensor fields in (\ref{alexandraoik})
with the extra factor $\nabla\phi_{u+1}$ contracting against a 
 factor $\nabla^{(2)}\Omega_h$. Denote the 
index set of those terms by $\overline{L}_\mu$. 
We will firstly prove that:

\begin{equation}
\label{basilissa}
\sum_{l\in L_\mu} a_l C^{l,i_1\dots
i_{\mu+\beta}}_g(\Omega_1,\dots,\Omega_p,\phi_1,\dots,\phi_u)
\nabla_{i_1}\phi_{u+1}\dots \nabla_{i_\beta}\phi_{u+\beta}
\nabla_{i_{\beta+1}}\upsilon\dots\nabla_{i_{\beta+\mu}}\upsilon=0.
\end{equation}

It suffices to prove the above for the sublinear 
combination of $\mu$-tensor fields 
where $\nabla\phi_{u+1}$ contracts against $\nabla^{(2)}\Omega_1$. 
(\ref{basilissa}) will then follows by repeating this step $p$ times.

We start by a preparatory claim: Let us denote by $\overline{L}_{\mu,\sharp}\subset \overline{L}_\mu$
the index set of $\mu$-tensor fields for which the factor $\nabla^{(2)}\Omega_1$
contains a free index, say the index ${}_{i_1}$ wlog. We will firstly prove that:

\begin{equation}
\label{sponsor} 
\sum_{l\in L_{\mu,\sharp}} a_l C^{l,i_1\dots i_\mu}_g\nabla_{i_1}\phi_{u+1}\dots\nabla_{i_\beta}\phi_{u+\beta}
\nabla_{i_{\beta+1}}\upsilon\dots\nabla_{i_{\beta+\mu}}\upsilon=0.
\end{equation}

{\it Proof of (\ref{sponsor}):} We will use the technique (introduced
in subsection 3.1 of \cite{alexakis5}) of ``inverse integration by parts'' followed by 
the silly divergence formula. 

Let us denote by $\hat{C}^l_g$ the complete contraction that 
arises from each $C^{l,i_1\dots i_{\mu+\beta}}_g$ by formally erasing 
the expression $\nabla^{(2)}_{si_1}\Omega_1\nabla^s\phi_{u+1}$ and then making all free 
indices ${}_{i_{\beta+1}},\dots ,{}_{i_{\beta+\mu}}$ 
into internal contractions.\footnote{We recall that 
to ``make a free index ${}_{i_y}$ into an internal contraction'' means 
that we add a derivative $\nabla_{i_y}$ onto 
the factor $T_{i_y}$  to which the free 
index ${}_{i_y}$ belongs. The new derivative 
index $\nabla^{i_y}$ is then contracted against the
index ${}_{i_y}$ in $T_{i_y}$.} Then, the 
``inverse integration by parts'' implies 
a new integral equation:

\begin{equation}
\label{oloklhrwma} 
\int_{M^n} \sum_{l\in L_\mu} a_l \hat{C}^l_g
+\sum_{j\in J} a_j C^j_g+\sum_{z\in Z} a_z C^z_gdV_g=0.
\end{equation}
Here the complete contractions indexed in $J$ have length $\sigma+u$, $u$ 
factors $\nabla\phi_u$ but are simply subsequent to the simple 
character $\vec{\kappa}_{simp}$. The terms indexed in $Z$ either
 have length $\ge\sigma+u+1$ or have length $\sigma+u$, but 
also at least one factor $\nabla^{(B)}\phi_h$ with $B\ge 2$. 
 
Now, in the above, we consider the complete contractions indexed in 
$\overline{L}_{\mu,\sharp}\subset \overline{L}_\mu$ and we ``pull out'' 
the expression $\Delta\nabla_t\Omega_1\nabla^t\phi_{u+1}$ to write: 
$$\sum_{l\in L_{\mu,\sharp}} a_l \hat{C}^l_g=
\sum_{l\in L_{\mu,\sharp}} a_l \overline{C}^l_g\cdot (\Delta
\nabla_t\Omega_1\nabla^t\phi_1).$$

Now, we consider the silly divergence formula applied to (\ref{oloklhrwma}) 
obtained by integrating by parts with respect to the 
function $\Omega_1$. If we denote the integrand in (\ref{oloklhrwma}) by $F_g$,
 we denote the resulting (local) equation by $silly[F_g]=0$. 
We consider the sublinear combination $silly^*[F_g]$ which consists of terms with length $\sigma+u$, 
$\mu$ internal contractions and $u-1+\beta$ factors $\nabla\phi_h, h\ge 2$ 
and a factor $\Delta\phi_{u+1}$. Clearly, this sublinear combination 
must vanish separately modulo longer terms:
$$silly^*[F_g]=0.$$ 
The above equation can be expressed as:
\begin{equation}
\label{papadopoulos}
Spread^{\nabla^s,\nabla_s}
[\sum_{l\in \overline{L}_{\mu,\sharp}} a_l \overline{C}^l_g]\cdot \Omega_1\cdot \Delta\phi_{u+1}=0.
\end{equation}
(Here $Spread^{\nabla^s,\nabla_s}$ is a formal operation 
that acts on complete contractions 
in the form (\ref{form1}) by hitting a factor $T$ in the form 
$\nabla^{(m)}R_{ijkl}$ or $\nabla^{(p)}\Omega_h$ with
a derivative $\nabla^s$ and then hitting another 
factor $T'\ne T$ in the form 
$\nabla^{(m)}R_{ijkl}$ or $\nabla^{(p)}\Omega_h$
 by a derivative $\nabla_s$ which contracts
against $\nabla^s$ and then adding over all the terms we can thus obtain.) 
Now, using the fact that (\ref{papadopoulos}) holds formally, 
 we derive:\footnote{This can be proven by using 
 the operation $Erase[\dots]$, see the 
 Appendix in \cite{alexakis4}.} 
\begin{equation}
\label{papadopoulos'}
\sum_{l\in L_{\mu,\sharp}} a_l \overline{C}^l_g=0.
\end{equation}
Thus, applying the operation $Sub_\upsilon$ 
$\mu-1$ times to the above and then multiplying 
by $\nabla_{i_1i_2}\Omega_1\nabla^{i_1}\upsilon\nabla^{i_2}\phi_{u+1}$ 
we derive (\ref{sponsor}). So for the rest of this proof we may assume that 
$\overline{L}_{\mu,\sharp}=\emptyset$. 
\newline

 Now we prove our claim under the additional assumption that
for the tensor fields indexed in $\overline{L}_\mu$, the 
factor $\nabla^{(2)}\Omega_1$ contains no free index.  
 
    We again refer to the equation (\ref{oloklhrwma})
  and perform integrations by parts with respect to the
   factor $\nabla^{(B)}\Omega_1$. We denote the resulting 
   local equation by $silly[L_g]=0$. We pick out the sublinear
    combination $silly^*[L_g]$ of terms with $\sigma+u$ factors, 
    $u+\beta$ factors $\nabla\phi_h$, $\mu$ internal 
    contractions, with $u+\beta-1$ factors 
$\nabla\phi_h,h\ge 2$ and a factor $\Delta\phi_1$. 
 This sublinear combination must vanish separately, $silly^*[L_g]=0$; the
 resulting new true equation can be described easily: 
Let us denote by $\hat{C}^{l,j_1}_g$ the 1-vector field that
 arises from $C^{l,i_1\dots i_\mu}_g, l\in L_{\mu,*}$ by formally erasing 
 the factor $\nabla^{(2)}_{js}\Omega_1\nabla^s\phi_1$, making the index ${}^j$ that
 contracted against ${}_j$  into a free index
  ${}_{j_1}$, and making all 
 the free indices ${}_{i_1},\dots, {}_{i_{\mu}}$ into internal 
 contractions.
 (Denote by $\vec{\kappa}'_{simp}$ the simple character of 
 these vector fields).
  Then the  
  equation $silly^*[L_g]=0$ can be expressed in the form:
  
  \begin{equation}
\label{dialogos}
\sum_{l\in L_{\mu,*}} a_l\{ Xdiv_{j_1}\hat{C}^{l,j_1}_g\}\Delta\phi_1
+\sum_{\in J} a_j C^j_g\Delta\phi_1=0;
\end{equation}
 here the complete contractions $C^j_g$ are simply subsequent
  to $\vec{\kappa}'_{simp}$. The above holds 
  modulo terms of length $\ge\sigma+u+1$. Now, we apply the 
 operation $Sub_\omega$ $\mu$ times (see the 
   Appendix in \cite{alexakis1}). In the case 
$\sigma>3$, applying the inductive
assumption of our Lemma \ref{additiongen} 
to the resulting equation (notice that 
the above falls under the inductive assumption of this Lemma
since we have lowered the weight in absolute value;
we ensure that Lemma \ref{additiongen}  
can be applied by just labeling one of the factors 
$\nabla\omega$ into $\nabla\phi_{u+1}$. 
 We derive (due to weight considerations) that there {\it can not be} 
 tensor fields of higher rank, thus:
 \begin{equation}
 \label{dikh}
 \sum_{l\in L_{\mu}} a_lSub_\omega^{\mu-1}[\hat{C}^{l,j_1}_g]\nabla_{i_1}\upsilon\Delta\phi_1=0.
 \end{equation}
Now, formally replacing the factor $\nabla_{i_1}\upsilon$ by 
$\nabla^{(2)}_{j_1t}\Omega_1\nabla^{t} \phi_1$, and then 
setting $\omega=\upsilon$, we derive the claim of 
our Lemma. In the case $\sigma=3$ (\ref{dikh}) follows by inspection, since 
the only two possible cases are $\sigma_2=2$ and $\sigma_1=2$; in the first case 
there are only two possible tensors field in 
$\overline{L}_\mu$ while in the second there are four. 
The equation (\ref{dialogos}) (by inspection) implies that the coefficients 
of all these tensor fields must vanish, which is equivalent to (\ref{dikh}). 
 \newline
 
 \par Now, we will prove our claim under the additional 
 assumption $\overline{L}_{\mu}=\emptyset$ (still for $p>0$). We again 
 refer to (\ref{oloklhrwma}) and again consider 
 the same equation $silly[L_g]=0$ as above. We now pick out the 
 sublinear combination of terms with $\sigma+u$ factors,
 $u+\beta$ factors $\nabla\phi_h$, and $\mu$ internal contractions. We 
 derive:

  \begin{equation}
\label{dialogos3}
\sum_{l\in L_{\mu}} a_l Xdiv_{j_1}Xdiv_{j_2}\hat{C}^{l,j_1j_2}_g
+\sum_{\in J} a_j C^j_g=0;
\end{equation} 
 here the terms $\hat{C}^{l,j_1j_2}_g$ arise from the $\mu$-tensor 
 fields $C^{li_1\dots i_\mu}_g$ by replacing all $\mu$ free indices
  by internal contractions, erasing the factor $\nabla_{jk}^{(2)}\Omega_1$
and making the indices ${}^j,{}^k$ into free indices ${}^{j_1},{}^{j_2}$.  Now,
 applying $Sub_\omega$ $\mu$ times, and then applying
  the inductive assumption of 
Lemma 4.10 \ref{additiongen} (this applies 
by length considerations as above for $\sigma>3$; while
  if $\sigma=3$ the claim (\ref{dialogos4}) will again follow by inspection) we derive: 
  
  \begin{equation}
\label{dialogos4}
\sum_{l\in L_{\mu}} a_l\hat{C}^{l,j_1j_2}_g\nabla_{j_1}\upsilon\nabla_{j_2}\upsilon=0;
\end{equation} 
 Replacing the expression $\nabla_{j_1}\upsilon\nabla_{j_2}\upsilon$
 by a factor $\nabla_{j_1j_2}^{(2)}\Omega_2$ and then setting $\omega=\upsilon$,
 we derive our claim in this case $p>0$. 
 \newline

  {\bf The case $p=0,\sigma_1>0$:} We will reduce ourselves to the previous case: We 
let $L_\mu^1$ the index set of $\mu$-tensor fields where the factor $T_1=S_{*}R_{ijkl}\nabla^i\tilde{\phi}_1$ 
contains a special free index (say the index ${}_k$ is the free index ${}_{i_{\beta+1}}$ wlog). 
We will prove our claim for the index set $L_\mu^1$; 
if we can prove this, then clearly 
our Lemma will follow by induction. 

 To prove this claim, we consider the first conformal 
variation of our hypothesis, $Image^1_Y[L_g]=0$, and we pick out the sublinear combination
of terms with length $\sigma+u+\beta$, with the factor 
$\nabla^{(\nu)}S_{*}R_{ijkl}\nabla^i\tilde{\phi}_1$ 
has been replaced by a factor $\nabla^{(\nu+2)}Y$, and the factor $\nabla\phi_1$ 
now contracts against a factor $T_2=R_{ijkl}$. 
This sublinear combination vanishes separately, 
thus we derive a new local equation. To describe 
the  resulting equation, we denote by 
\\$\hat{C}^{l,i_1\dots \hat{i}_{\beta+1}\dots i_{\mu+\beta}}_g
(Y,\phi_1,\dots,\phi_u)\nabla_{i_1}\phi_{u+1}
\dots\nabla_{i_{\beta}}\phi_{u+\beta}$ 
the $(\mu-1)$-tensor field 
that arises from by formally replacing the factor by and also 
adding a derivative index $\nabla_{i_*}$ onto the factor 
$T_2=R_{ijkl}$ and then contracting 
that index against a factor $\nabla\phi_1$. 
Denote the $(u-1)$-simple character of the above 
(the one defined by $\nabla\phi_2,\dots,\nabla\phi_u$) by $\vec{\kappa}'_{simp}$. 
 We then have an equation:

\begin{equation}
\label{papatsakalos} 
\begin{split}
&\sum_{l\in L_\mu^1} a_l Xdiv_{i_{\beta+2}}\dots Xdiv_{i_{\beta+\mu}} 
\hat{C}^{l,i_1\dots \hat{i}_{\beta+1}\dots i_{\mu+\beta}}_g(Y,\phi_1,\dots,\phi_u)\nabla_{i_1}\phi_{u+1}
\dots\nabla_{i_{\beta}}\phi_{u+\beta}
\\&+\sum_{h\in H} a_h Xdiv_{i_{\beta+1}}\dots Xdiv_{i_{\beta+\mu}}
C^{l,i_2\dots i_{\mu+\beta}}_g(Y,\phi_1,\dots,\phi_u)\nabla_{i_{1}}\phi_{u+1}
\dots\nabla_{i_{\beta+\mu}}\phi_{u+\beta}
\\&+\sum_{j\in J} a_j C^{j,i_{\mu+1}\dots i_{\mu+\beta}}_g
(Y,\phi_1,\dots,\phi_u)\nabla_{i_{\mu+1}}\phi_{u+1}
\dots\nabla_{i_{\mu+\beta}}\phi_{u+\beta}.
\end{split}
\end{equation}
The terms indexed in $H$ are acceptable, have a $(u-1)$-simple character 
$\vec{\kappa}'_{simp}$ and the factor $\nabla\phi_1$ contracts against the index ${}_i$ in 
the factor $T_2=R_{ijkl}$; writing that factor as $S_*R_{ijkl}\nabla^i\tilde{\phi}_1$, 
we denote the resulting $u$-simple factor by $\tilde{\kappa}_{simp}$. 
The terms indexed in $J$ are simply subsequent to 
$\vec{\kappa}'_{simp}$.  Now, applying the inductive assumption of 
Lemma \ref{additiongen},\footnote{The terms indexed in $L_\mu^1$ are
 now simply subsequent to $\tilde{\kappa}_{simp}$.} we derive that:

\begin{equation}
\label{de3iwsh}
 \sum_{h\in H} a_h 
C^{l,i_2\dots i_{\mu+\beta}}_g(Y,\phi_1,\dots,\phi_u)\nabla_{i_{1}}\phi_{u+1}
\dots\nabla_{i_{\beta+1}}\phi_{u+\beta}\nabla_{i_{\beta+1}}\upsilon\dots\nabla_{i_{\beta+\mu}}\upsilon=0.
\end{equation}
Thus, we may assume wlog that $H=\emptyset$ in (\ref{papatsakalos}). 
Now, we again apply Lemma \ref{additiongen} to (\ref{papatsakalos})
(under that additional assumption), and we derive that:

\begin{equation}
\label{de3iwsh2}
\sum_{l\in L_\mu^1} a_l  
\hat{C}^{l,i_1\dots \hat{i}_{\beta+1}\dots i_{\mu+\beta}}_g(Y,\phi_1,\dots,\phi_u)\nabla_{1}\phi_{u+1}
\dots\nabla_{i_{\beta+1}}\phi_{u+\beta}\nabla_{i_{\beta+2}}\upsilon\dots\nabla_{i_{\beta+\mu}}\upsilon=0. 
\end{equation}
Now, erasing the factor $\nabla\phi_1$ from the above, and then formally 
replacing the factor $\nabla^{(2)}_{ab}Y$ by 
$S_*R_{i(ab)l}\nabla^i\tilde{\phi}_1\nabla^l\upsilon$, we derive our claim.

 {\bf The case $p=0,\sigma_1=0$:} In this case $\sigma=\sigma_2$. In other 
 words, all factors  in $\vec{\kappa}_{simp}$ are simple factors
  in the form $S_{*}R_{ijkl}\nabla^i\tilde{\phi}_h$. 
We recall that in this case all 
$\mu$-tensor fields in (\ref{alexandraoik}) must have at 
most one free index in any factor $S_{*}R_{ijkl}$. 
  In that case, we will prove our claim in a more convoluted
   manner, again reducing ourselves to the inductive
    assumption of Proposition 2.1 in \cite{alexakis4}.
    
A key observation is the following: By the definition of the special cases, 
$\mu+\beta\le \sigma_2$. In the case of strict inequality, 
we see (by a counting argument) that at least one of
 the special indices in one of the factors $S_*R_{ijkl}$ 
 must contract against a special index in another
factor $S_{*}R_{abcd}$. In the case $\mu+\beta=\sigma_2$ 
this remains true, {\it except} for the 
terms for which the $\beta$ factors $\nabla\phi_{u+h}$ 
contract against special indices, say the indices ${}_k$, 
in $\beta$ factors $T_y=S_*R_{ikl}\nabla^i\tilde{\phi}_y$, 
 and moreover these factors must not contain a  free 
 index, and all 
 other factors $S_*R_{ikl}$ contain exactly one free index,
which must be special. 
 In this subcase, we will prove our claim 
 for all $\mu$-tensor fields {\it excluding} this particular ``bad'' 
 sublinear combination; we will prove our claim 
 for this sublinear combination  in the end. 
\newline

 We will now proceed to {\it normalize} the 
 different $(\mu+\beta)$-tensor fields in (\ref{alexandraoik}). 
 A {\it normalized} tensor field will be in the form (\ref{form2}), 
 with possibly certain pairs of indices in 
 certain of the factors $S_{*}R_{ijkl}$ being symmetrized. 
 
  Let us first introduce some definitions:
  Given each $C^{l,i_1\dots i_\mu}_g$, we 
  list out the factors $T_1,\dots T_{\sigma_2}$ in the form
   $S_*R_{ikl}$. Here $T_a$ is the factor for which the 
   index ${}_i$ is contracting against the factor $\nabla\tilde{\phi}_a$. 
  We say that factors $S_{*}R_{ikl}$ are type I if they 
  contain no free index. We say they are of type II if they contain 
  a special free index. We say they are of type
   III if they contain a non-special free index.

Given any tensor field  $C^{l,i_1\dots i_\mu}_g$ in the form (\ref{form2}), 
pick out the pairs of factors, $T_\alpha,T_\beta$ in the form $S_*R_{ijkl}$ 
for which a special index in $T_\alpha$ contracts 
against a a special index in $T_\beta$. (Call
such particular contractions ``special-to-special'' 
particular contractions). Now, in any $C^{l,i_1\dots i_\mu}_g$ 
we define an ordering among all its factors $S_*R_{ijkl}$: 
The factor $T_a=S_*R_{ikl}\nabla^i\tilde{\phi}_a$ is more important 
than $T_b=S_*R_{i'j'k'l'}\nabla^{i'}\tilde{\phi}_b$ if $a<b$.

 Now, consider a tensor field $C^{l,i_1\dots i_\mu}_g$ and list out all the 
 pairs of factors $T_a,T_b$ with a 
 special-to-special particular contraction. We say that $(T_a,T_b)$ 
 is the {\it most important} pair of factors with a special-to-special particular 
 contraction\footnote{Asume wlog that $T_a$ is
 more important than $T_b$ } if any other such 
 pair $(T_c,T_d)$,\footnote{Again assume wlog
  that $T_c$ is more important than $T_d$.} 
  has either $T_c$ being less important than $T_a$ or 
  $T_a=T_c$ and $T_d$ less important than $T_b$.

Now, consider a tensor field $C^{l,i_1\dots i_\mu}_g$ and
consider the most important pair of factors $(T_a,T_b)$ with a 
 special-to-special particular contraction.
Assume wlog that the index ${}_l$ in 
$T_a=S_{*}R_{ijkl}\nabla^i\tilde{\phi}_a$
   contracts against the index ${}_{l'}$ in 
   $T_b=S_*R_{i'j'k'l'}\nabla^{i'}\tilde{\phi}_b$.
    We say that  $C^{l,i_1\dots i_\mu}_g$ is 
  normalized if both factors $T_a, T_b$ are normalized. 
The factor $T_a=S_*R_{ikl}\nabla^i\tilde{\phi}_a$ is normalized  if:
 Either the index ${}_j$ contracts against 
a factor $T_c$
which is more important than $T_b$, {\it or} 
 if the indices ${}_j,{}_k$ are symmetrized. 
If $T_a$ is of type II, then we require that the index ${}_j$ in $T_b=S_*R_{ij(free)l}$ 
{\it must} be contracting against  a special index of some other factor $T_c$,
 and moreover $T_c$ must be more important than $T_b$.
If $T_a$ is of type $III$, then it is automatically normalized.
The same definition applies to $T_b$, where any reference to $T_b$ 
must be replaced by a reference to $T_a$.

Let us now prove that we may assume wlog that all $\mu$-tensor 
fields in (\ref{alexandraoik}) are normalized:
 Consider a $C^{l,i_1\dots i_\mu}_g$ in (\ref{alexandraoik}) 
 for which the most important pair of factors with a  
 special-to-special particular contraction is the pair $(T_a,T_b)$.    
 We will prove that we can write:
   
\begin{equation}
\label{grammos}C^{l,i_1\dots i_{\mu+\beta}}_g=\tilde{C}^{l,i_1\dots i_{\mu+\beta}}_g+
\sum_{t\in T} a_t C^{t,i_1\dots i_{\mu+\beta}}_g;
\end{equation}
here the term $\tilde{C}^{l,i_1\dots i_{\mu+\beta}}_g$
 is normalized,  the most important pair 
of factors with a  special-to-special particular 
contraction is the pair $(T_a,T_b)$, and moreover 
its refined double character is the same as for $C^{l,i_1\dots i_{\mu+\beta}}_g$. 
Each $C^{t,i_1\dots i_{\mu+\beta}}_g$ has either the same,  
or a doubly subsequent refined double character   to 
$C^{l,i_1\dots i_{\mu+\beta}}_g$; moreover in the first case its most important pair 
of factors with a special-to-special particular contraction
will be less important than the pair $(T_a,T_b)$. In the second case 
the most important pair will either be $(T_a,T_b)$ or a less important pair. 

Clearly, if we can prove the above, then by iterative repetition we 
may assume wlog that all $(\mu+\beta)$-tensor fields in (\ref{alexandraoik}) 
are normalized. 

{\it Proof of (\ref{grammos}):} Pick out the most important pair 
of factors with a  
special-to-special particular contraction is the 
pair $(T_a,T_b)$ in $C^{l,i_1 \dots i_{\mu+\beta}}_g$. 
Let us first normalize $T_a$. If $T_a$ is 
of type III, there is nothing to do. If it is of type II 
and already normalized, there is again nothing to do. 
If it is of type II and not normalized, then we {\it interchange} the 
indices ${}_j,{}_k$. The resulting factor {\it is} normalized. The correction term we obtain
by virtue of the first Bianchi identity is also normalized (it is of type III). 
Moreover, the resulting tensor field is doubly subsequent 
to $C^{l,i_1 \dots i_{\mu+\beta}}_g$.  Finally, if the factor $T_a$ is of type I, we 
inquire on the factor $T_c$ against which ${}_j$ in $T_a=S_*R_{ijkl}$ contracts:
If it is more important than $T_b$, then we leave $T_a$ as it is; it is already normalized. 
If not, we symmetrize ${}_j,{}_k$. The resulting tensor field is 
normalized. The correction term we obtain by virtue 
of the first Bianchi identity will then have the same refined double
 character as $C^{l,i_1\dots i_{\mu+\beta}}_g$,  and moreover its most important pair 
of factors with a special-to-special particular contraction
 is less important than that pair $(T_a, T_b)$.
\newline

We now prove the claim of Lemma \ref{additiongen} in this special case,
 under the additional assumption that all tensor fields in (\ref{alexandraoik}) 
 are normalized. We list out the most important 
 pair of special-to-special particular contractions 
in each $C^{l,i_1\dots i_{\mu+\beta}}_g$, 
 and denote it by $(a,b)_l$. We let $(\alpha,\beta)$ stand for the lexicographicaly minimal 
 pair among the list $(a,b)_l, l\in L_\mu$. 
 We denote by $L_\mu^{(\alpha,\beta)}\subset L_\mu$
 the index set of terms with a special-to-special particular contraction
 among the terms $T_\alpha,T_\beta$. We will prove that: 
 
 \begin{equation}
 \label{sweden}
 \sum_{l\in L_\mu^{(\alpha,\beta)}} a_l
  C^{l,i_1\dots i_{\mu+\beta}}_g\nabla_{i_1}\upsilon\dots\nabla_{i_\mu}\upsilon=0.
 \end{equation}
Clearly, the above will imply our claim,
 by iterative repetition.\footnote{In 
the subcase ${\mu+\beta}=\sigma_2$
 it will only imply it for the ``excluded'' sublinear combination 
 defined above.} 
\newline

{\it Proof of (\ref{sweden}):} Consider $Image^2_{Y_1,Y_2}[L_g]=0$ and 
pick out the sublinear combination where the factors $T_\alpha,T_\beta$ are
 replaced by $\nabla^{(A)}Y_1\otimes g, \nabla^{(B)}Y_2\otimes g$, and the
 two factors $\nabla\tilde{\phi}_\alpha,\nabla\tilde{\phi}_\beta$
 contract against each other. The resulting sublinear combination must 
 vanish separately. We erase the expression 
 $\nabla_t\tilde{\phi}_\alpha\nabla^t\tilde{\phi}_\beta$,\footnote{Denote 
 the resulting $(u-2)$-simple character by $\vec{\kappa}'''_{simp}$.} 
 and derive a new true eqation which will be in the form:
  
  \begin{equation}
  \label{write}
  \sum_{l\in L_\mu^{(\alpha,\beta)}} a_l
Xdiv_{i_1}\dots Xdiv_{i_\mu} 
 \tilde{C}^{l,i_1\dots i_{\mu+\beta}}_g(\Omega_1,Y_1,Y_2)
 +\sum_{j\in J} a_j C^j_g(\Omega_1,Y_1,Y_2)=0;
  \end{equation}
here the tensor fields $\tilde{C}^{l,i_1\dots i_{\mu+\beta}}_g(\Omega_1,Y_1,Y_2)$
arise from the tensor fields $C^{l,i_1\dots i_{\mu+\beta}}_g$ by 
replacing the expression \\$\nabla^i\tilde{\phi}_{\alpha}S_*R_{ijkl}
\otimes S_*{R_{i'jk}}^l\nabla^{i'}\tilde{\phi}_\beta$
by an expression $\nabla_{jk}Y_1\otimes \nabla_{j'k'}Y_2$ (notice we
have lowered the weight in absolute value).
   
 Now, applying the inductive assumption of Lemma \ref{additiongen}
 to the above,\footnote{We have lowered the weight in absolute value.} we derive:  
  
\begin{equation}
  \label{write'}
  \sum_{l\in L_\mu^{(\alpha,\beta)}} a_l
 \tilde{C}^{l,i_1\dots i_{\mu+\beta}}_g(\Omega_1,Y_1,Y_2)
 \nabla_{i_1}\upsilon\dots\nabla_{i_\mu}\upsilon=0;
  \end{equation}  
  The proof of (\ref{sweden}) is only one step away. 
Let us start with an important observation: 
For each given complete contraction above, examine the factor 
$\nabla^{(2)}_{zx}Y_1$;
it either contracts against no factor $\nabla\upsilon$ or one factor
$\nabla\upsilon$.\footnote{The two corresponding sublinear 
combinations vanish separately, of course.} In the first case, 
the factor $\nabla^{(2)}_{zx}Y_1$ must have
arisen from a factor $S_*R_{ijkl}$ of type I. 
In fact, the indices ${}_z,{}_x$ correspond to the indices ${}_j,{}_k$
in the original factor, and we can even determine their position:
Since the pair $(\alpha,\beta)$ is the most important pair in (\ref{alexandraoik}), 
at most one of the indices ${}_z,{}_x$ can contract against a special index in a 
more important factor than $T_\beta$. If one of them does (say ${}_z$), 
then that index must have been the index ${}_j$
in $T_\alpha=S_*R_{ikl}$. If none of them does, then the two 
indices ${}_z,{}_x$ must be symmetrized over, since the two indices 
${}_j,{}_k$ in $T_\alpha$ to which they correspond were symmetrized over.  
Now, these two separate sublinear combinations in (\ref{write'}) must vanish 
separately (this can be proven using the eraser fro the
Appendix in \cite{alexakis1}), and  furthermore in the first case, 
we may assume that the index ${}_z$ (which contracts against a 
special index in a more important factor than $T_\beta$) 
occupies the leftmost position in $\nabla^{(2)}_{zx}Y_1$ 
and is {\it not permuted} in the formal permutations of 
indices that make (\ref{write'}) hold {\it formally}).  
    
On the other hand, consider the terms in (\ref{write'}) with the factor 
$\nabla^{(2)}Y_1$  contracting against a factor $\nabla\upsilon$. 
By examining the index ${}_y$ in the factor $\nabla^{(2)}_{yt}Y_1\nabla^t\upsilon$,
 we can determine the {\it type} of factor in $C^{l,i_1\dots i_{\mu+\beta}}_g$ 
 from which the factor $\nabla^{(2)}Y_1$ arose: If the index ${}_y$ is contracting against 
 a special index in a factor $S_{*}R_{ijkl}$ which is {\it more important} than $T_\beta$,
 then $\nabla^{(2)}Y_1$ can only have arisen from a factor 
 of type II in $C^{l,i_1\dots i_{\mu+\beta}}_g$. In fact, the index ${}_y$ in $\nabla^{(2)}Y_1$ 
 must correspond to the index ${}_j$ in $S_* R_{ij(free)l}$ in $T_\alpha$. If the index 
 ${}_y$ in $\nabla^{(2)}_{yt}Y_1\nabla^t\upsilon$
 does not contract against a special index in 
 a factor $T_c$ which is more important than $T_\beta$, then 
 the factor $\nabla^{(2)}Y_1$ 
can only have arisen from a factor 
of type III in $C^{l,i_1\dots i_{\mu+\beta}}_g$. 
In fact, the index ${}_y$ in $\nabla^{(2)}Y_1$ 
must correspond to the index ${}_k$ in $S_* R_{i(free)kl}$ in $T_\alpha$. 
   
 The same analysis can be repeated for the factor $\nabla^{(2)}Y_2$, with
any reference to the factor $T_\beta$ now replaced by the factor $T_\alpha$. 
    
In view    of the above analysis, we can break the LHS of (\ref{write'}) into
four sublinear combinations which vanish separately (depending 
on whether $\nabla^{(2)}Y_1,\nabla^{(2)}Y_2$ contract against a  
factor $\nabla\upsilon$ or not). Then in each of the four sublinear combinations,
 we can arrange that in the formal permutations that make 
the LHS of (\ref{write'}) formally zero, the two indices in the factors 
$\nabla^{(2)}Y_1,\nabla^{(2)}Y_2$ are {\it not} permuted (by virtue of the remarks above).
In view of this and the analysis in the previous paragraph, we can then  
{\it replace}  the two factors $\nabla^{(2)}_{zx}Y_1,\nabla^{(2)}_{qw}Y_2$ by an expression 
 $\nabla^i\tilde{\phi}_\alpha S_*R_{izxl}\otimes
 S_*{R_{i'qw}}^l\nabla^{i'}\tilde{\phi}_\beta$,
 in such a way that the resulting linear combination {\it vanishes formally
  without permuting the two indices ${}_q,{}_w,{}_{q'}, {}_{w'}$}. This proves our claim, except 
  for   the subcase $\mu+\beta=\sigma_2$
where we only derive our claim  for all terms except for the ``bad sublinear combination''.
We now prove our claim for that. 
\newline

 We then break up the LHS of (\ref{paragon})
 according to which factor $T_s$ the factor $\nabla\phi_{u+1}$
 contracts--denote the index set of those terms by $L_\mu^K$.
  Denote the resulting sublinear 
 combinations by $L^K_g, K=1,\dots ,\sigma_2$. Given any $K$,  
 we consider the eqation $Image^1_Y[L_g]=0$, and we pick out the 
 sublinear combination where the term
  $\nabla^{(B)}S_*R_{ijkl}\nabla^i\tilde{\phi}_K$,
 is replaced by $\nabla^{(B+2)}Y$, and the factor $\nabla\phi_K$ 
 now contracts against the factor $\nabla\phi_{u+1}$. This 
 sublinear combination must vanish separately. We then 
again perform the ``inverse integration by parts'' to this true
 equation (deriving an integral equation), and then we 
 consider the silly divergence formula for this integral equation, 
 obtained by integrating by parts with respect to $\nabla^{(B)}Y$. 
 We pick out the sublinear combination with $\sigma+u+\beta$ factors, 
 $\mu$ internal contractions and $u+\beta$
factors $\nabla\phi_h$, and a expression $\nabla_s\phi_{u+1}\nabla^s\tilde{\phi}_K$
This gives us a new true local equation:
 
 \begin{equation}
 \label{marcus}
 \sum_{l\in L_\mu^K} a_l X_*div_{j_1}X_*div_{j_2}\tilde{C}^{l,j_1j_2}_g+\sum_{j\in J} a_j C^j_g=0.
 \end{equation}
 Here the tensor fields $\tilde{C}^{l,j_1j_2}_g$ arise from $C^{l,i_1\dots i_\mu}_g$
by formally replacing all $\mu$ free indices by
internal contractions, and also replacing the expression 
$\nabla_{x}\phi_{u+1}\otimes S_* {R_{i(jk)}}^x\nabla^i\tilde{\phi}_K$ by 
$\nabla_x\phi_{u+1}\nabla^s\tilde{\phi}_K\otimes Y$, 
and then making the 
indices ${}^j,{}^k$ that contracted against ${}_j,{}_k$ 
into free indices ${}^{j_1},{}^{j_2}$. $X_{*}div_j$ stands for the
sublinear combination in $Xdiv_j$ where $\nabla^j$ 
is not allowed to hit the factor $Y$. Now, applying 
the inductive assumption of Lemma \ref{additiongen}  
to the above,\footnote{We have lowered the weight in absolute value.} 
we derive that: 
    
    $$ \sum_{l\in L_\mu^K} a_l\tilde{C}^{l,_1j_2}_g\nabla_{j_1}\omega\nabla_{j_2}\omega=0.$$
 Now, we replace the expression $\nabla^x\phi_K\nabla_x\phi_{u+1}\nabla_{j_1}\omega\nabla_{j_2}\omega Y$
by \\$\nabla^l\phi_{u+1}S_*R_{i(j_1j_2)l}\nabla^i\tilde{\phi}_K$ 
and then replacing all internal contractions by
factors $\nabla\upsilon$ (applying
the operation $Sub_\upsilon$ from the
Appendix in \cite{alexakis1}). The resulting (true) equation  
is precisely our remaining claim for the ``bad'' sublinear
combination. $\Box$

\subsection{Proof of Lemmas 4.6, 4.8 in \cite{alexakis4}: The main part.}

\par We first write down the form of the complete and partial 
 contractions that we are dealing with in
Lemmas \ref{obote} and \ref{vanderbi}. In the setting of Lemma
\ref{obote} we recall that the tensor fields $C^{h,i_1\dots i_\alpha}$
 indexed in $H_2$ (in the hypothesis of Lemma \ref{obote}) are all 
partial contractions in the form:

\begin{equation}
\label{form2obote}
\begin{split}
&pcontr(\nabla^{(m_1)}R_{ijkl}\otimes\dots\otimes
\nabla^{(m_{\sigma_1})} R_{ijkl}\otimes
S_{*}\nabla^{(\nu_1)}R_{ijkl}\otimes\dots\otimes
S_{*}\nabla^{(\nu_t)} R_{ijkl}\otimes
\\& \nabla^{(b_1)}\Omega_1\otimes\dots\otimes \nabla^{(b_p)}\Omega_p\otimes\nabla Y\otimes
\\& \nabla\phi_{z_1}\dots \otimes\nabla\phi_{z_f}\otimes\nabla
\phi'_{z_{f+1}}\otimes
\dots\otimes\nabla\phi'_{z_{f+d}}\otimes\dots
\otimes\nabla\tilde{\phi}_{z_{f+d+1}}\otimes\dots\otimes
\nabla\tilde{\phi}_{z_{f+d+y}}),
\end{split}
\end{equation}
where we let $f+d+y=u'$.  The main assumption here is that
all tensor fields have the same
 $u'$-simple character (the one defined by
 $\nabla\phi_1,\dots,\nabla\phi_{u'}$), which we denote by $\vec{\kappa}^{+}_{simp}$. The
 other main assumption is that if we formally treat the factor $\nabla Y$
 as a function $\nabla\phi_{u+1}$, then the hypothesis of Lemma
 \ref{obote} falls under the inductive assumptions of Proposition
 \ref{giade} (i.e. the weight, real length, $\Phi$
  and $p$ are as in our inductive assumption of Proposition \ref{giade}).

\par In the setting of Lemma \ref{vanderbi} we recall that we are dealing
with complete and partial contractions in the form:

\begin{equation}
\label{form2vanderbi}
\begin{split}
&contr(\nabla^{(m_1)}R_{ijkl}\otimes\dots\otimes
\nabla^{(m_{\sigma_1})} R_{ijkl}\otimes
S_{*}\nabla^{(\nu_1)}R_{ijkl}\otimes\dots\otimes
S_{*}\nabla^{(\nu_t)} R_{ijkl}\otimes
\\& \nabla^{(b_1)}\Omega_1\otimes\dots\otimes \nabla^{(b_p)}\Omega_p\otimes
[\nabla\omega_1\otimes\nabla\omega_2]\otimes
\\& \nabla\phi_{z_1}\dots \otimes\nabla\phi_{z_f}
\otimes\nabla\phi'_{z_{f+1}}\otimes \dots\otimes\nabla
\phi'_{z_{f+d}}\otimes\dots\otimes\tilde{\phi}_{z_{f+d+1}}
\otimes\dots\otimes\tilde{\phi}_{z_{f+d+y}}),
\end{split}
\end{equation}
where we let $f+d+y=u'$. The main assumption here is that
all partial contractions have the same
 $u'$-simple character (the one defined by
 $\nabla\phi_1,\dots,\nabla\phi_{u'}$), which we denote by $\vec{\kappa}^{+}_{simp}$. The
 other main assumption is that if we formally treat the factors $\nabla \omega_1,\nabla\omega_2$
 as factors $\nabla\phi_{u+1},\nabla\phi_{u+2}$, then the hypothesis of Lemma
 \ref{vanderbi} falls under the inductive assumptions of Proposition
 \ref{giade} (i.e. the weight, real length, $\Phi$
  and $p$ are as in our inductive assumption of Proposition \ref{giade}).

 {\it Note:} From now on, we will be writing $u'=u$, for
 simplicity.  We will also be writing $\vec{\kappa}_{simp}^{+}=\vec{\kappa}_{simp}$, for simplicity.
 We will also be labelling the indices
 ${}_{i_1},\dots ,{}_{i_\alpha}$ as ${}_{i_{\pi+1}},\dots ,{}_{i_{\alpha+1}}$.
\newline

{\it New induction:} We will now prove the two Lemmas \ref{obote}, \ref{vanderbi}
by a new induction on the weight of the complete
contractions in the hypotheses of those Lemmas.
We will assume that these two Lemmas are true when the weight of the
 complete contractions in their hypotheses is $-W$, for any $W<K\le n$.
  We will then show our Lemmas for weight $-K$.
\newline

{\it Reduce Lemma \ref{obote} to two Lemmas:} In order to show
Lemma \ref{obote}, we further break up
 $H_2$ into
subsets: We say that $h\in H_2^a$ if and only if
 $C^{h,i_{\pi+1}\dots i_{\alpha+1}}$ has a free
 index ( say the free index ${}_{i_{\alpha +1}}$ wlog)
  belonging to the factor
$\nabla Y$. On the other hand, we say that $h\in H^b_2$ if the
index in the factor $\nabla Y$ is not free.
Lemma \ref{obote} will then follow from Lemmas
 \ref{firstclaim}, \ref{secondclaim} below:

\begin{lemma}
\label{firstclaim}
 There exists a linear combination of acceptable
$(\alpha-\pi+1)$-tensor fields, $\Sum_{v\in V} a_v
C^{v,i_{\pi+1}\dots i_{\alpha +2}}_{g}(\Omega_1,\dots ,\Omega_p,
Y, \phi_1,\dots,\phi_u)$, where the
 index ${}_{i_{\alpha +1}}$ belongs to the factor $\nabla Y$, with
 a simple character $\vec{\kappa}_{simp}$, so that:

\begin{equation}
\label{episkepsh}
\begin{split}
&\Sum_{h\in H_2^a} a_h C^{h,i_{\pi+1}\dots i_{\alpha
+1}}_{g}(\Omega_1,\dots ,\Omega_p, Y,\phi_1,\dots,\phi_u)
\nabla_{i_{\pi+1}}\upsilon\dots\nabla_{i_{\alpha+1}}\upsilon=
\\&\Sum_{v\in V} a_v X_{*} div_{i_{\alpha +2}}
C^{v,i_{\pi+1}\dots i_{\alpha+2}}_{g}(\Omega_1,\dots ,
\Omega_p,Y,\phi_1,\dots,\phi_u)\nabla_{i_{\pi+1}}\upsilon\dots
\nabla_{i_{\alpha+1}}\upsilon+
\\&\Sum_{j\in J} a_j C^{j,i_{\pi+1}\dots
i_{\alpha+1}}_{g}(\Omega_1,\dots,\Omega_p, Y,\phi_1,\dots,
\phi_u)\nabla_{i_{\pi+1}}\upsilon\dots\nabla_{i_{\alpha+1}}
\upsilon.
\end{split}
\end{equation}
 Each $C^j$ is simply subsequent to
$\vec{\kappa}_{simp}$.
\end{lemma}

\par We observe that if we can show our first claim, then we
 can assume, with
 no loss of generality, that $H_2^a=\emptyset$, since
it immediately follows from the above that:

\begin{equation}
\label{stefania}
\begin{split}
&\Sum_{h\in H_2^a} a_h X_{*}div_{i_{\pi+1}}\dots X_{*}div_{i_{\alpha +1}}
C^{h,i_{\pi+1}\dots i_{\alpha +1}}_{g} (\Omega_1,\dots ,\Omega_p,
Y,\phi_1,\dots,\phi_u)=
\\&\Sum_{v\in V} a_v X_{*}div_{i_{\pi+1}}\dots
X_{*}div_{i_{\alpha +1}} X_{*}div_{i_{\alpha +2}} C^{v,i_{\pi+1}\dots
i_{\alpha+2}}_{g}(\Omega_1,\dots , \Omega_p,Y,
\phi_1,\dots,\phi_u)
\\&+\Sum_{j\in J} a_j C^j_{g}
(\Omega_1,\dots,\Omega_p, Y,\phi_1,\dots,
\phi_u),
\end{split}
\end{equation}
where each complete contraction $C^j$ is subsequent to
$\vec{\kappa}_{simp}$. (Note that one of the free indices in the tensor fields
$C^{v,i_{\pi+1}\dots i_{\alpha+2}}_{g}$ will belong to the factor $\nabla Y$).
\newline

\par Second claim, in the setting of Lemma \ref{obote}:
\begin{lemma}
\label{secondclaim}
 We assume $H^a_2=\emptyset$. We then claim that modulo
complete contractions of length $\ge\sigma +u +1$:

\begin{equation}
\label{episkepsh4}
\begin{split}
&\Sum_{h\in H_2} a_h C^{h,i_{\pi+1}\dots i_{\alpha
+1}}_{g}(\Omega_1,\dots ,\Omega_p, Y,\phi_1,\dots,\phi_u)
\nabla_{i_{\pi+1}}\upsilon\dots\nabla_{i_{\alpha+1}}\upsilon=
\\&\Sum_{t\in T} a_t X_{*}div_{i_{\alpha +2}}
C^{t,i_{\pi+1}\dots i_{\alpha +2}}_{g} (\Omega_1,\dots
,\Omega_p,Y,\phi_1,\dots,\phi_u)
\nabla_{i_{\pi+1}}\upsilon\dots\nabla_{i_{\alpha+1}}\upsilon+
\\&\Sum_{j\in J} a_j C^j_{g}
(\Omega_1,\dots,\Omega_p, Y,\phi_1,\dots,\phi_u),
\end{split}
\end{equation}
where each $C^j$ is acceptable and subsequent to
$\vec{\kappa}_{simp}$.
\end{lemma}

\par We observe that if we can show the above two
Lemmas then Lemma \ref{obote} will follow. (Notice
 that replacing by the RHS of (\ref{stefania}) into the hypothesis of Lemma \ref{obote},
 we do not introduce 1-forbidden terms).
\newline

\par We make two analogous claims for Lemma \ref{vanderbi}:
\newline

{\it Reduce Lemma \ref{vanderbi} to two Lemmas:} We say that
$h\in H^a_2$ if and only if $C^{h,i_{\pi+1}\dots i_{\alpha+1}}$
has a free
 index belonging to one of the factors
 $\nabla\omega_1,\nabla\omega_2$.
On the other hand, we say that $h\in H^b_2$ if in none of the
factors $\nabla\omega_1,\nabla\omega_2$ in $C^{h,i_{\pi+1}\dots
i_{\alpha+1}}$ contains a free index. (Observe that we may assume
with no loss of generality that there are no tensor fields
$C^{h,i_{\pi+1}\dots i_{\alpha+1}}$ with free indices in both
factors $\nabla\omega_1,\nabla\omega_2$-this is by virtue of the
anti-symmetry of the factors $\nabla\omega_1,\nabla\omega_2$).  We make two claims.
Firstly:

\begin{lemma}
\label{firstclaimb}
 There is a linear combination of acceptable
$(\alpha-\pi+1)$-tensor fields, $\Sum_{v\in V} a_v
C^{v,i_{\pi+1}\dots i_{\alpha +2}}_{g}(\Omega_1,\dots
,\Omega_p,[\omega_1,\omega_2], \phi_1,\dots,\phi_u)$, in the
 form (\ref{form2obote}) with a simple character $\vec{\kappa}_{simp}$, so that:

\begin{equation}
\label{episkepshb}
\begin{split}
&\Sum_{h\in H^a_2} a_h X_{+}div_{i_{\pi+1}}\dots
X_{+}div_{i_{\alpha+1}} C^{h,i_{\pi+1}\dots
i_{\alpha+1}}_{g}(\Omega_1,\dots ,\Omega_p,
[\omega_1,\omega_2],\phi_1,\dots,\phi_u)=
\\&\Sum_{v\in V} a_v X_{+}div_{i_{\pi+1}}\dots
X_{+}div_{i_{\alpha+1}} X_{+}div_{i_{\alpha+2}}
C^{v,i_{\pi+1}\dots i_{\alpha+2}}_{g}(\Omega_1,\dots
,\Omega_p,[\omega_1,\omega_2],\phi_1,\dots,\phi_u)
\\&+\Sum_{q\in Q} a_q X_{+}div_{i_{\pi+1}}\dots
X_{+}div_{i_{\alpha +1}} C^{q,i_{\pi+1}\dots i_{\alpha+1}}_{g}
(\Omega_1,\dots ,\Omega_p,
\nabla_{+}[\omega_1,\omega_2],\phi_1,\dots, \phi_u)
\\&+\Sum_{j\in J} a_j C^j_{g}
(\Omega_1,\dots,\Omega_p, [\omega_1,\omega_2],\phi_1,\dots,\phi_u).
\end{split}
\end{equation}
(Recall that by definition the complete contractions indexed in $Q$ have a
factor $\nabla^{(2)}\omega_1$).
\end{lemma}

\par We observe that if we can show our first claim, then we
 can, with
 no loss of generality, assume that $H^a_2=\emptyset$.

\par Second claim:
\begin{lemma}
\label{secondclaimb}
 We assume $H^a_2=\emptyset$, and that for some $k\ge 1$,
we can write:

\begin{equation}
\label{episkepsh3b}
\begin{split}
&\Sum_{h\in H^b_2} a_x X_{+}div_{i_{\pi+1}}\dots
X_{+}div_{i_{\alpha+1}} C^{h,i_{\pi+1}\dots i_{\alpha+1}}_{g}
(\Omega_1,\dots ,\Omega_p,
[\omega_1,\omega_2],\phi_1,\dots,\phi_u)
\\&+\Sum_{t\in T_k} a_t X_{+}div_{i_{\pi+1}}\dots
X_{+}div_{i_{\alpha +k}} C^{t,i_{\pi+1}\dots
i_{\alpha+k}}_{g}(\Omega_1,\dots
,\Omega_p,[\omega_1,\omega_2],\phi_1,\dots,\phi_u) 
\\& +\Sum_{q\in Q} a_q X_{+}div_{\pi+1}\dots
X_{+}div_{i_{\alpha +1}} C^{h,i_{\pi+1}\dots i_{\alpha
+1}}_{g}(\Omega_1,\dots ,\Omega_p,
\nabla_{+}[\omega_1,\omega_2],\phi_1,\dots,\phi_u)
\\&+\Sum_{j\in J} a_j C^j_{g}
(\Omega_1,\dots,\Omega_p, [\omega_1,\omega_2],\phi_1,\dots,\phi_u),
\end{split}
\end{equation}
where the last two linear combinations on
 the right hand side of the above are generic linear combinations in
 the form described in the claim of Lemma
 \ref{vanderbi}.\footnote{In Lemma \ref{vanderbi}, $Q$ is called $V$.} On the other hand,
$\Sum_{t\in T_k} a_t C^{t,i_{\pi+1}\dots i_{\alpha+k}}_{g}
(\Omega_1,\dots ,\Omega_p, [\omega_1,\omega_2],
\phi_1,\dots,\phi_u)$ is a linear combination of acceptable
$(\alpha-\pi+k)$-tensor fields in the form (\ref{form2vanderbi})
with a simple character $\vec{\kappa}_{simp}$, and with two
anti-symmetric factors $\nabla \omega_1,\nabla\omega_2$ which do
not contain a free index. We then claim that modulo complete
contractions of length $\ge\sigma +u +1$ we can write:

\begin{equation}
\label{episkepsh4'}
\begin{split}
&\Sum_{t\in T_k} a_t X_{+}div_{i_{1}}\dots X_{+}div_{i_{a +k}}
C^{t,i_{1}\dots i_{a +k}}_{g}
(\Omega_1,\dots,\Omega_p,[\omega_1,\omega_2],\phi_1,\dots,\phi_u)=
\\&\Sum_{t\in T_{k+1}} a_t X_{+}div_{i_{\pi +1}}\dots
X_{+}div_{i_{a +k+1}} C^{t,i_{\pi+1}\dots
i_{a+k+1}}_{g}(\Omega_1,\dots
,\Omega_p,[\omega_1,\omega_2],\phi_1,\dots,\phi_u) 
\\& +\Sum_{q\in Q} a_q X_{+}div_{i_{1}}\dots
X_{+}div_{i_{a +1}} C^{q,i_{1}\dots i_{a+1}}_{g}(\Omega_1,\dots
,\Omega_p, \nabla_{+}[\omega_1,\omega_2],\phi_1,\dots, \phi_u)
\\&+\Sum_{j\in J} a_j C^j_{g}
(\Omega_1,\dots,\Omega_p,
[\omega_1,\omega_2],\phi_1,\dots,\phi_u),
\end{split}
\end{equation}
with the same notational conventions as above.
\end{lemma}

\par We observe that if we can show the above two
claims, then Lemma \ref{vanderbi} will follow by iterative repetition of the second claim.

\par We will now show the four Lemmas above.
\newline

{\it Proof of Lemmas \ref{secondclaim} and \ref{secondclaimb}:}
Lemma \ref{secondclaim} is a direct consequence of Lemma
4.10 in \cite{alexakis4}.\footnote{Observe that
our hypotheses on the tensor fields in the equation in Lemma \ref{obote} not being ``bad''
ensure that we do not fall under the ``forbidden'' cases of Lemma 4.10 in \cite{alexakis4}.} Lemma
\ref{secondclaimb} can be proven in two steps:  Firstly, by
 Lemma \ref{addition2}
 we derive that there exists a  linear combination of acceptable $(a+k+1)$-tensor
 fields (indexed in $X$ below) with a u-simple character $\vec{\kappa}_{simp}$ so that:

\begin{equation}
\label{episkepsh4''}
\begin{split}
&\Sum_{t\in T_k} a_t C^{t,i_{1}\dots i_{a +k}}_{g}
(\Omega_1,\dots,\Omega_p,[\omega_1,\omega_2],\phi_1,\dots,\phi_u)
\nabla_{i_1}\upsilon\dots\nabla_{i_{a+k}}\upsilon-
\\&\Sum_{t\in T_{k+1}} a_t X_{*}div_{i_{a +k+1}} C^{t,i_{\pi+1}\dots
i_{a+k+1}}_{g}(\Omega_1,\dots
,\Omega_p,[\omega_1,\omega_2],\phi_1,\dots,\phi_u)
\nabla_{i_1}\upsilon\dots\nabla_{i_{a+k}}\upsilon 
\\&+\Sum_{j\in J} a_j C^j_{g}
(\Omega_1,\dots,\Omega_p,
[\omega_1,\omega_2],\phi_1,\dots,\phi_u,\upsilon^{a+k}),
\end{split}
\end{equation}
where the complete contractions indexed in $J$ have length $\sigma+a+k+1$ and are
 simply subsequent to $\vec{\kappa}_{simp}$.

\par Then, making the factors $\nabla\upsilon$ in the above
 into $X_{+}div$s, we derive Lemma \ref{secondclaimb}. $\Box$
\newline

{\it Proof of Lemma \ref{firstclaim}:}

\par We have denoted by $\vec{\kappa}_{simp}$ the simple
 character of our tensor fields.
 We distinguish two cases: In case A there is a factor
 $\nabla^{(m)}R_{ijkl}$ in $\vec{\kappa}_{simp}$, and in case $B$
 there is no such factor.

\par  We denote $\alpha+1=\gamma$, for
brevity.

Now we break the set $H^{b}_2$ into subsets: In case A we say that
$h\in H_2^{b,+}$ if and only if $\nabla Y$ is contracting against
an internal index of a factor $\nabla^{(m)}R_{ijkl}$. In  case B
we say that $h\in H_2^{b,+}$ if and only if $\nabla Y$ is
contracting against one of the indices ${}_k,{}_l$ in a factor
$S_{*}\nabla^{(\nu)}R_{ijkl}$.

We define $H_2^{b,-}=H^{b}_2\setminus H^{b,+}_2$.

\par In each of the above cases and subcases we treat
the function $\nabla Y$ as a function $\nabla\phi_{u+1}$ in our
Lemma hypothesis. Then, by applying the first claim
in Lemma 4.10 in \cite{alexakis4}\footnote{By weight considerations, since we started out with no ``bad
 terms'' in Lemma \ref{obote}, we will not encounter no ``forbidden tensor fields'' 
for Lemma 4.10 in \cite{alexakis4}.}
 to our
Lemma hypothesis and then making the $\nabla\upsilon$s into $X_{*}divs$,
 we derive that we can write:

\begin{equation}
\label{singer}
\begin{split}
& X_{*}div_{i_{\pi+1}}\dots X_{*}div_{i_\gamma} \Sum_{h\in
H^{b,+}_2} a_h C^{h,i_{\pi+1}\dots i_\alpha,i_\gamma}_{g}
(\Omega_1,\dots ,\Omega_p,Y,\phi_1,\dots ,\phi_u) =
\\&  X_{*}div_{i_{\pi+1}}\dots X_{*}div_{i_\gamma}
\Sum_{h\in H^{b,*,-}_2} a_h C^{h,i_{\pi+1}\dots
i_\alpha,i_\gamma}_{g} (\Omega_1,\dots ,\Omega_p,Y,\phi_1,\dots
,\phi_u) \\& + \Sum_{j\in J} a_j
C^j_{g}(\Omega_1,\dots,\Omega_p,Y, \phi_1,\dots ,\phi_u),
\end{split}
\end{equation}
where $\Sum_{h\in H^{b,*,-}_2} a_h C^{h,i_{\pi+1}\dots
i_\alpha,i_\gamma}_{g} (\Omega_1,\dots ,\Omega_p,Y,\phi_1,\dots
,\phi_u)$ stands for a {\it generic} linear combination as defined
above (i.e. it is in the general form $\sum_{h\in H^{b}_2}\dots$
but the factor $\nabla Y$ is {\it not} contracting against a
special index in any factor $\nabla^{(m)}R_{ijkl}$ or
$S_{*}\nabla^{(\nu)}R_{ijkl}$.\footnote{Recall that a special index in a
 factor $\nabla^{(m)}R_{ijkl}$ is an internal index, while a special
 index in a factor $S_{*}\nabla^{(\nu)}R_{ijkl}$ is an index ${}_k,{}_l$.}
 On the other hand, each
$C^j_{g}(\Omega_1,\dots,\Omega_p,Y,\phi_1,\dots ,\phi_u)$
 is a complete contraction with a simple character that is
subsequent to $\vec{\kappa}_{simp}$.

\par Thus, by virtue of (\ref{singer}),  we
reduce ourselves to the case where
$H^{b,+}_2=\emptyset$. We will then show Lemma \ref{firstclaim}
separately in cases A and B, under the assumption that
$H^{b,+}_2=\emptyset$.
\newline

{\it Proof of Lemma \ref{firstclaim} in case A:} We will define
the C-crucial factor, for the purposes of this 
proof only: We denote by $Set$ the set of numbers $u$
for which $\nabla\phi_u$ is contracting against one of the factors
$\nabla^{(m)}R_{ijkl}$. If $Set\ne\emptyset$, we define $u_{+}$ to
be the minimum element of $Set$, and we pick out the factor
$\nabla^{(m)}R_{ijkl}$ in each $C^h$ against which
$\nabla\phi_{u_{+}}$ contracts. We call that factor $\nabla^{(m)}
R_{ijkl}$ {\it C-crucial}. If $Set=\emptyset$, we will say {\it
the C-crucial factors} and will mean {\it any} of the factors
$\nabla^{(m)}R_{ijkl}$.

 Now, we pick out the subset
$H^{b,*}_2\subset H^b_2$, which is defined by the rule: $h\in
H^{b,*}_2$ if and only if $\nabla Y$ is contracting against the
(one of the) C-crucial factor(s).

\par Now, for each $h\in H^a_2$ we denote by
$$Hitdiv_{i_\gamma} C^{h,i_{\pi+1}\dots i_{\alpha+1}}_{g}
(\Omega_1,\dots ,\Omega_p,Y,\phi_1,\dots,\phi_u)$$
the sublinear combination in $X_{*}div_{i_\gamma} C^{h,i_{\pi+1}\dots
i_{\alpha+1}}_{g} (\Omega_1,\dots,\Omega_p,Y,\phi_1,
\dots,\phi_u)$ that arises when $\nabla_{i_\gamma}$ hits the (one of the) 
C-crucial factor.\footnote{Recall that ${}_{i_\gamma}$ is the
 free index that belongs to $\nabla Y$.} It then follows that:

\begin{equation}
\label{hat}
\begin{split}
&\Sum_{h\in H^a_2} a_h X_{*}div_{i_{\pi+1}}\dots X_{*}div_{i_\alpha}
Hitdiv_{i_\gamma} C^{h,i_{\pi+1}\dots i_{\alpha+1}}_{g}
(\Omega_1,\dots ,\Omega_p,Y,\phi_1,\dots,\phi_u)
\\& + \Sum_{h\in H^{b,*}_2} a_h X_{*}div_{i_{\pi+1}}\dots
X_{*}div_{i_\gamma} C^{h,i_{\pi+1}\dots
i_{\alpha+1}}_{g}(\Omega_1,\dots ,
\Omega_p,Y,\phi_1,\dots,\phi_u)+
\\&\Sum_{j\in J} a_j C^j_{g}(\Omega_1,
\dots ,\Omega_p,Y,\phi_1,\dots ,\phi_u),
\end{split}
\end{equation}
where each $C^j_{g}$ has the factor $\nabla Y$ contracting against
the C-crucial factor $\nabla^{(m)}R_{ijkl}$ and is simply
subsequent to $\vec{\kappa}_{simp}$.

  We now denote
the $(u+1)$-simple character (the one defined by
$\nabla\phi_1$,$\dots$ ,$\nabla\phi_{u+1}=\nabla Y$) of the tensor
fields $Hitdiv_{i_\gamma} C^{h,i_{\pi+1}\dots
i_\alpha,i_\gamma}_{g} (\Omega_1,\dots ,\Omega_p,Y,\phi_1,\dots
,\phi_u)$ by $\vec{\kappa}'_{simp}$. (Observe that they all have
the same $(u+1)$-simple  character).

We observe that just applying Corollary 1 in \cite{alexakis4}
 to (\ref{hat}) (all tensor fields are acceptable and have the same
 simple character $\vec{\kappa}'_{simp}$),\footnote{Notice that by 
weight considerations, since we started out with no ``bad'' terms in the hypothesis 
of Lemma \ref{obote}, there is no danger of falling under
a ``forbidden case'' of that Corollary.} we obtain an equation:

\begin{equation}
\label{ane3ignw}
\begin{split}
&\Sum_{h\in H^a_2} a_h  Hitdiv_{i_\gamma} C^{h,i_{\pi+1}\dots
i_\alpha}_{g} (\Omega_1,\dots ,\Omega_p,Y,\phi_1,\dots ,\phi_u)
\nabla_{i_{\pi+1}}\upsilon\dots \nabla_{i_\alpha}\upsilon+
\\&\Sum_{u\in U} a_u Xdiv_{i_{\alpha+1}}
C^{u,i_{\pi+1}\dots i_\alpha,i_{\alpha+1}}_{g} (\Omega_1,\dots
,\Omega_p,Y,\phi_1,\dots ,\phi_u) \nabla_{i_{\pi+1}}\upsilon\dots
\nabla_{i_\alpha}\upsilon=
\\& \Sum_{j \in J} a_j C^{j,i_{\pi+1}\dots i_\alpha}_{g}
(\Omega_1,\dots ,\Omega_p,Y,\phi_1,\dots ,\phi_u)
\nabla_{i_{\pi+1}}\upsilon\dots\nabla_{i_\alpha}\upsilon =0,
\end{split}
\end{equation}
where the tensor fields indexed in $U$ are acceptable (we are treating
 $\nabla Y$ as a factor $\nabla\phi_{u+1}$), have a simple character
$\vec{\kappa}'_{simp}$ and
 each $C^j$ is simply subsequent to $\vec{\kappa}'_{simp}$.

\par But then, our first claim follows almost immediately.
We recall the operation $Erase_{\nabla Y}[\dots]$ from the Appendix in \cite{alexakis1} which acts on
the complete contractions in the above by erasing the factor
$\nabla Y$ and the (derivative) index that it contracts
 against. Then, since (\ref{ane3ignw})
holds formally, we have that the tensor field required for Lemma
\ref{firstclaim} is:

$$\Sum_{u\in U} a_u
Erase_{\nabla Y}[C^{u,i_{\pi+1}\dots i_\alpha,i_{\alpha+1}}_{g}
(\Omega_1,\dots ,\Omega_p,Y,\phi_1,\dots ,\phi_u)]\cdot
\nabla_{i_\gamma}Y.$$
\newline

{\it Proof of Lemma \ref{firstclaim} in case B:}  We again
distinguish two subcases: In subcase $(i)$ there is some
non-simple factor $S_{*}\nabla^{(\nu)}R_{ijkl}$ in
$\vec{\kappa}_{simp}$ or a non-simple factor
$\nabla^{(B)}\Omega_x$ contracting against two factors
$\nabla\phi'_h$ in $\vec{\kappa}_{simp}$. In subcase $(ii)$ there
are no such factors.

In the subcase $(i)$, we arbitrarily pick out one factor
$S_{*}\nabla^{(\nu)}R_{ijkl}$ or $\nabla^{(B)}\Omega_x$ with the
properties described above and call it the $D$-crucial factor. In
this first subcase we will show our claim for the whole sublinear
combination $\sum_{h\in H^a_2}\dots$ in one piece.

 In the subcase $(ii)$, we will introduce some notation: We will
 examine each factor $T=S_{*}\nabla^{(\nu)}R_{ijkl}$,
 $T=\nabla^{(B)}\Omega_x$ in each tensor field $C^{h,i_{\pi+1}\dots
i_\alpha,i_{\alpha+1}}_{g}$ and define its ``measure'' as follows:
If $T=S_{*}\nabla^{(\nu)}R_{ijkl}$ then its ``measure'' will stand
for its total number of free indices {\it plus $\frac{1}{2}$}. If
$T=\nabla^{(B)}\Omega_x$ then its ``measure'' will stand for its
total number of free indices {\it plus} the number of factors
$\nabla\phi_h$ against which it is contracting.

 We divide the
index set $H_2^a$ into subsets according to the measure of any
given factor. We denote by $M$ the maximum measure among all
factors among the tensor fields $C^{h,i_{\pi+1}\dots
i_\alpha,i_{\alpha+1}}_{g}$, $h\in H_2^a$. We denote by
$H^{2,*}_a\subset H^a_2$ the index set of the tensor fields which
contain a factor of maximum measure. We will show the claim of
Lemma \ref{firstclaim} for the sublinear combination $\sum_{h\in
H^{2,*}_a}\dots$. Clearly, if we can do this, then Lemma
\ref{firstclaim} will follow by induction.

\par We will prove Lemma \ref{firstclaim} in the second subcase (which is the hardest).
The proof in the first subcase follows by the same argument, only by disregarding any
reference to $M$ free indices belonging to a given factor etc.

{\it Proof of Lemma \ref{firstclaim} in case B for the sublinear
combination $\sum_{h\in H^{2,*}_a}\dots$:}

We will further divide $H^{2,*}_a$ into subsets, $H^{2,*,k}_a,
k=1,\dots ,\sigma$, according to the factor of maximum measure:
Firstly, we order the factors $S_{*}\nabla^{(\nu)}R_{ijkl},\dots$
\\$\nabla^{(p)}\Omega_h$ in $\vec{\kappa}_{simp}$, and label them
$T_1,\dots ,T_\sigma$ (observe each factor is well-defined in
$\vec{\kappa}_{simp}$, because we are in case B). We then say that
$h\in H^{a,*,1}_2$ if in $C^{u,i_{\pi+1}\dots i_\alpha}_{g}$ the
factor $T_1$ has measure $M$. We say say that $h\in H^{a,*,2}_2$
if in $C^{u,i_{\pi+1}\dots i_\alpha}_{g}$ the factor $T_2$ has
measure $M$ and $T_1$ has measure less than $M$, etc. We will then
prove our claim for each of the index sets $h\in
H^{a,*,k}_2$:\footnote{Again we observe that if we can prove this
then Lemma \ref{firstclaim} in case B will follow by induction.}
We arbitrarily pick a $k\le K$ and show our claim for $\sum_{h\in
H^{2,*,k}_a}\dots$.

\par For the purposes of this proof, we call the factor $T_k$ the
$D$-crucial factor.

 Now, we pick out the subset
$H^{b,k}_2\subset H^b_2$, which is defined by the rule: $h\in
H^{b,k}_2$ if and only if $\nabla Y$ is contracting against the
D-crucial factor $T_k$.

\par Now, for each $h\in H^a_2$ we denote by
$$Hitdiv_{i_\gamma} C^{h,i_{\pi+1}\dots i_{\alpha+1}}_{g}
(\Omega_1,\dots ,\Omega_p,Y,\phi_1,\dots,\phi_u)$$ the sublinear
combination in $Xdiv_{i_\gamma} C^{h,i_{\pi+1}\dots
i_{\alpha+1}}_{g} (\Omega_1,\dots,\Omega_p,Y,\phi_1,
\dots,\phi_u)$ that arises when $\nabla_{i_\gamma}$ hits the
D-crucial factor.\footnote{Recall that ${}_{i_\gamma}={}_{i_{\alpha+1}}$
 belongs to $\nabla Y$ by hypothesis.} It then follows that:

\begin{equation}
\label{hat'}
\begin{split}
&\Sum_{h\in H^{a}_2} a_h Xdiv_{i_{\pi+1}}\dots Xdiv_{i_\alpha}
Hitdiv_{i_\gamma} C^{h,i_{\pi+1}\dots i_{\alpha+1}}_{g}
(\Omega_1,\dots ,\Omega_p,Y,\phi_1,\dots,\phi_u)
\\& + \Sum_{h\in H^{b,k}_2} a_h Xdiv_{i_{\pi+1}}\dots
Xdiv_{i_\gamma} C^{h,i_{\pi+1}\dots
i_{\alpha+1}}_{g}(\Omega_1,\dots ,
\Omega_p,Y,\phi_1,\dots,\phi_u)+
\\&\Sum_{j\in J} a_j C^j_{g}(\Omega_1,
\dots ,\Omega_p,Y,\phi_1,\dots ,\phi_u),
\end{split}
\end{equation}
where each $C^j_{g}$ has the factor $\nabla Y$ contracting against
the D-crucial factor and is simply subsequent to
$\vec{\kappa}_{simp}$.

  We now denote
the $(u+1)$-simple character (the one defined by
$\nabla\phi_1,\dots $, $\nabla\phi_{u+1}=\nabla Y$) of the tensor
fields $Hitdiv_{i_\gamma} C^{h,i_{\pi+1}\dots
i_\alpha,i_\gamma}_{g} (\Omega_1,\dots ,\Omega_p,Y,\phi_1,\dots
,\phi_u)$ by $\vec{\kappa}'_{simp}$. (Observe that they all have
the same $(u+1)$-simple  character).

We  apply Corollary 1 in \cite{alexakis4}
 to (\ref{hat'}) (all tensor fields are acceptable and have the same
 simple character $\vec{\kappa}'_{simp}$) and then pick out the
 sublinear combination where there are $M$ factors
 $\nabla\upsilon$ or $\nabla\phi_h$ or $\nabla\phi'_h$
 contracting against $T_k$, we obtain an equation:

\begin{equation}
\label{ane3ignw'}
\begin{split}
&\Sum_{h\in H^{a,*,k}_2} a_h  Hitdiv_{i_\gamma}
C^{h,i_{\pi+1}\dots i_\alpha}_{g} (\Omega_1,\dots
,\Omega_p,Y,\phi_1,\dots ,\phi_u) \nabla_{i_{\pi+1}}\upsilon\dots
\nabla_{i_\alpha}\upsilon+
\\&\Sum_{u\in U} a_u Xdiv_{i_{\alpha+1}}
C^{h,i_{\pi+1}\dots i_\alpha,i_{\alpha+1}}_{g} (\Omega_1,\dots
,\Omega_p,Y,\phi_1,\dots ,\phi_u) \nabla_{i_{\pi+1}}\upsilon\dots
\nabla_{i_\alpha}\upsilon=
\\& \Sum_{j \in J} a_j C^{j,i_{\pi+1}\dots i_\alpha}_{g}
(\Omega_1,\dots ,\Omega_p,Y,\phi_1,\dots ,\phi_u)
\nabla_{i_{\pi+1}}\upsilon\dots\nabla_{i_\alpha}\upsilon =0,
\end{split}
\end{equation}
where the tensor fields indexed in $U$ are acceptable and  have a simple character
$\vec{\kappa}'_{simp}$ and
 each $C^j$ is simply subsequent to $\vec{\kappa}'_{simp}$.

\par Now, observe that if $M\ge \frac{3}{2}$,
 we can apply the eraser to $\nabla Y$ (see the 
Appendix of \cite{alexakis1}) and the index it is
contracting against in the $D$-crucial factor and derive our
conclusion as in case A.

\par On the other hand, in the remaining
 cases\footnote{Observe that the remaining cases
 are when $M=0$, $M=\frac{1}{2}$, $M=1$.}
 the above argument cannot be directly applied. In those cases, we
 derive our claim as follows:

In the case $M=1$ the $D$-crucial factor is of the form
$\nabla^{(p)}\Omega_h$, then we cannot directly derive our claim
by the above argument, because if for some tensor fields in $U$
above we have $\nabla Y$ contracting according to the pattern
$\nabla_iY\nabla^{ij}\Omega_h\nabla_j\psi$ (where $\psi=\upsilon$
or $\psi=\phi_h$), then we will not obtain acceptable tensor
fields after we apply the eraser. Therefore, if $M=1$ and the
$D$-crucial factor is of the form $\nabla^{(p)}\Omega_h$, we apply  Lemma 
4.6 in \cite{alexakis4} to (\ref{ane3ignw'}) (treating the
factors $\nabla\upsilon$ as factors
$\nabla\phi$)\footnote{Furthermore, we can observe 
that we do not fall under a  ``forbidden case'' of Lemma 4.1 in \cite{alexakis4}, 
by weight considerations, and since the tensor fields in our Lemma assumption are not ``bad''.}
 to obtain a new
equation in the form (\ref{ane3ignw'}), where each tensor field in
$U$ has the factor $\nabla Y$ is contracting against a factor
$\nabla^{(l)}\Omega_h$, $l\ge 3$.\footnote{Note that the weight becomes less negative,
hence Lemma 4.10 in \cite{alexakis4} applies.}
Then, applying the eraser as
explained, we derive our Lemma \ref{firstclaim} in this case.

\par When $M=\frac{1}{2}$ or $M=0$, then we first
 apply the inductive assumptions of
Corollaries 3,2 in \cite{alexakis4} (respectively) to
(\ref{ane3ignw'}),\footnote{By our assumptions there will be a removable
 index in these cases. Hence our extra requirements of those Lemmas are fulfilled.}
 in order to assume with no loss of generality
that for each tensor field indexed in $U$ there, the factor
$\nabla Y$ is either contracting against a factor $\nabla^{(B)}\Omega_h$,
$B\ge 3$ or a factor $S_{*}\nabla^{(\nu)}R_{ijkl}$, $\nu\ge 1$.
Then the eraser can be applied and produces acceptable tensor
fields. Hence, applying $Erase_{\nabla Y}$ to (\ref{ane3ignw'}) we
derive our claim. $\Box$
\newline

{\it Proof of Lemma \ref{firstclaimb}:}
\newline

\par We re-write the hypothesis of
Lemma \ref{vanderbi} (which is also the hypothesis of Lemma
\ref{firstclaimb}) in the following form:

\begin{equation}
\label{mayrh9b} \begin{split} &\Sum_{h\in H_2} a_h
X_{*}div_{i_{\pi+1}}\dots X_{*}div_{i_{\alpha+1}}
\{C^{h,i_{1}\dots i_{\alpha+1}}_{g}(\Omega_1,\dots ,\Omega_p,
\omega_1,\omega_2,\phi_1,\dots,\phi_u)-
\\&Switch[C]^{h,i_{\pi+1}\dots i_{\alpha+1}}_{g}
(\Omega_1,\dots ,\Omega_p,\omega_1,\omega_2,\phi_1,\dots,
\phi_u)\}
\\&=\Sum_{j\in J} a_j C^j_{g}(\Omega_1,\dots ,
\Omega_p,[\omega_1,\omega_2],\phi_1,\dots,\phi_u).
\end{split}
\end{equation}
Here the operation $Switch$ interchanges the indices ${}_a$
and ${}_b$ in the two factors $\nabla_a\omega_1$,
$\nabla_b\omega_2$.

{\it Notational conventions:} We have again denoted by
$H^a_2\subset H_2$ the index set of those vector fields for which
one of the free indices (say ${}_{i_{\alpha+1}}$) belongs to a
factor $\nabla\omega_1$ or $\nabla\omega_2$. With no loss of
generality we assume that for each $h\in H_2^a$
${}_{i_{\alpha+1}}$ belongs to the factor $\nabla\omega_1$.
 We can clearly do this, due to the antisymmetry
of the factors $\nabla\omega_1,\nabla\omega_2$.

\par We have defined $H^b_2=H_2\setminus H_2^a$. For each $h\in H^b_2$ we denote by
 $T_{\omega_1},T_{\omega_2}$ the factors against which
 $\nabla\omega_1,\nabla\omega_2$ are contracting. Also, for each
 $h\in H_2^a$ we will denote by $T_{\omega_2}$ the factor against
 which $\nabla\omega_2$ is
 contracting.\footnote{Note that the definition of $T_{\omega_1},T_{\omega_2}$
 depends on $h$; however, to simplify notation we
  suppress the index $h$ that should appear in $T_{\omega_1},T_{\omega_2}$.}

\par For each $h\in H_2$, we will call
the factors $T_{\omega_1},T_{\omega_2}$ against which
$\nabla\omega_1$ or $\nabla\omega_2$ are contracting
``problematic'' in the following cases: If $T_{\omega_1}$ or
$T_{\omega_2}$ is of the form $\nabla^{(m)}R_{ijkl}$ and
$\nabla\omega_1$ or $\nabla\omega_2$ is contracting against
 an internal index. Alternatively, if $T_{\omega_1}$ or $T_{\omega_2}$ is of
 the form
$S_{*}\nabla^{(\nu)} R_{ijkl}$ and the factor $\nabla\omega_1$ or
$\nabla\omega_2$ is contracting one of the indices ${}_k$ or
${}_l$.

We then define a few subsets of $H^a_2$, $H^b_2$:

\begin{definition}
\label{subsets} We define $H^b_{2,**}$ to stand for the index set
of the tensor fields $C^{h,i_{\pi+1}\dots i_{\alpha+1}}_{g}$'s for
which $\nabla\omega_1,\nabla\omega_2$ are contracting against
different factors and both $T_{\omega_1},T_{\omega_2}$ are
problematic.

We define $H^a_{2,*}\subset H^a_2$ to be the index set of the
tensor fields $C^{h,i_{\pi+1}\dots i_{\alpha+1}}_{g}$'s for which
$T_{\omega_2}$ is problematic.

We define $H^b_{2,*}$ to stand for the index set of the tensor
fields $C^{h,i_{\pi+1}\dots i_{\alpha+1}}_{g}$'s for which either
$T_{\omega_1}=T_{\omega_2}$ or $T_{\omega_1}\ne T_{\omega_2}$ and one
of the factors $T_{\omega_1},T_{\omega_2}$ is problematic.
\end{definition}

\par Abusing notation, we will be using the symbols
$\Sum_{h\in H^b_{2,*}}\dots$ etc to denote {\it generic}
 linear combinations as explained above, when these symbols appear in the right hand sides of
 the equations below.
\newline

\par We then state three preparatory claims:

Firstly, we claim that we can write:

\begin{equation}
\label{zekkos1}
\begin{split} &\Sum_{h\in H^b_{2,**}} a_h
X_{+}div_{i_{\pi+1}}\dots X_{+}div_{i_{\alpha+1}}
\{C^{h,i_{\pi+1}\dots i_{\alpha+1}}_{g}(\Omega_1,\dots ,\Omega_p,
\omega_1,\omega_2,\phi_1,\dots,\phi_u)
\\&-Switch[C]^{h,i_{\pi+1}\dots i_{\alpha+1}}_{g}
(\Omega_1,\dots ,\Omega_p,\omega_1,\omega_2,\phi_1,\dots,
\phi_u)\}=
\\&\Sum_{h\in H^b_{2,*}} a_h
X_{+}div_{i_{\pi+1}}\dots X_{+}div_{i_{\alpha+1}}
\{C^{h,i_{\pi+1}\dots i_{\alpha+1}}_{g}(\Omega_1,\dots ,\Omega_p,
\omega_1,\omega_2,\phi_1,\dots,\phi_u)
\\&-Switch[C]^{h,i_{\pi+1}\dots i_{\alpha+1}}_{g}
(\Omega_1,\dots ,\Omega_p,\omega_1,\omega_2,\phi_1,\dots,
\phi_u)\} \\&+\Sum_{j\in J} a_j C^j_{g}(\Omega_1,\dots ,
\Omega_p,[\omega_1,\omega_2],\phi_1,\dots,\phi_u),
\end{split}
\end{equation}
where the linear combination $\sum_{h\in H^b_{2,*}}\dots$ on the
RHS stands for a generic linear combination in the form described
above. Observe that if we can show (\ref{zekkos1}) then we may
assume with no loss of generality that $H^b_{2,**}=\emptyset$ in
our Lemma hypothesis.

\par Then, assuming that $H^b_{2,**}=\emptyset$ in
our Lemma hypothesis we will show that there exists a linear
combination of $(\alpha-\pi+1)$-tensor fields (indexed in $X$
below) which are in the form (\ref{form2''}) with a simple
character $\vec{\kappa}_{simp}$ so that:

\begin{equation}
\label{zekkos2}
\begin{split} &\Sum_{h\in H^a_{2,*}} a_h
\{C^{h,i_{\pi+1}\dots i_{\alpha+1}}_{g}(\Omega_1,\dots ,\Omega_p,
\omega_1,\omega_2,\phi_1,\dots,\phi_u)-
\\&Switch[C]^{h,i_{\pi+1}\dots i_{\alpha+1}}_{g}
(\Omega_1,\dots ,\Omega_p,\omega_1,\omega_2,\phi_1,\dots,
\phi_u)\}\nabla_{i_{\pi+1}}\upsilon\dots
\nabla_{i_{\alpha+1}}\upsilon-
\\&X_{*}div_{i_{\alpha+2}}\Sum_{x\in X} a_x
\{C^{x,i_{1}\dots i_{\alpha+1}i_{\alpha+2}}_{g}(\Omega_1,\dots
,\Omega_p, \omega_1,\omega_2,\phi_1,\dots,\phi_u)-
\\&Switch[C]^{h,i_{\pi+1}\dots i_{\alpha+1}}_{g}
(\Omega_1,\dots ,\Omega_p,\omega_1,\omega_2,\phi_1,\dots,
\phi_u)\}\nabla_{i_{\pi+1}}\upsilon\dots
\nabla_{i_{\alpha+1}}\upsilon
\\&+\Sum_{h\in H^b_{2,*}} a_h
\{C^{h,i_{1}\dots i_{\alpha+1}}_{g}(\Omega_1,\dots ,\Omega_p,
\omega_1,\omega_2,\phi_1,\dots,\phi_u)-
\\&Switch[C]^{h,i_{\pi+1}\dots i_{\alpha+1}}_{g}
(\Omega_1,\dots ,\Omega_p,\omega_1,\omega_2,\phi_1,\dots,
\phi_u)\}\nabla_{i_{\pi+1}}\upsilon\dots
\nabla_{i_{\alpha+1}}\upsilon \\& =\Sum_{j\in J} a_j
C^j_{g}(\Omega_1,\dots ,
\Omega_p,[\omega_1,\omega_2],\phi_1,\dots,\phi_u,\upsilon^{\alpha-\pi}).
\end{split}
\end{equation}

\par We observe that if we can show the above, we may then assume
that $H^a_{2,*}=\emptyset$ (and $H^b_{2,**}=\emptyset$) in the hypothesis of Lemma
\ref{firstclaimb}.

\par Finally, under the assumption that
$H^b_{2,**}=H^a_{2,*}=\emptyset$ in our Lemma hypothesis, we will
show that we can write:

\begin{equation}
\label{zekkos3}
\begin{split} &\Sum_{h\in H^b_{2,*}} a_h
X_{+}div_{i_{\pi+1}}\dots X_{+}div_{i_{\alpha+1}}
\{C^{h,i_{1}\dots i_{\alpha+1}}_{g}(\Omega_1,\dots ,\Omega_p,
\omega_1,\omega_2,\phi_1,\dots,\phi_u)-
\\&Switch[C]^{h,i_{\pi+1}\dots i_{\alpha+1}}_{g}
(\Omega_1,\dots ,\Omega_p,\omega_1,\omega_2,\phi_1,\dots,
\phi_u)\}=
\\&\Sum_{h\in H^b_{2,OK}} a_h
X_{+}div_{i_{\pi+1}}\dots X_{+}div_{i_{\alpha+1}}
\{C^{h,i_{1}\dots i_{\alpha+1}}_{g}(\Omega_1,\dots ,\Omega_p,
\omega_1,\omega_2,\phi_1,\dots,\phi_u)-
\\&Switch[C]^{h,i_{\pi+1}\dots i_{\alpha+1}}_{g}
(\Omega_1,\dots ,\Omega_p,\omega_1,\omega_2,\phi_1,\dots,
\phi_u)\} \\&+\Sum_{j\in J} a_j C^j_{g}(\Omega_1,\dots ,
\Omega_p,[\omega_1,\omega_2],\phi_1,\dots,\phi_u),
\end{split}
\end{equation}
where the sublinear combination $\Sum_{h\in H^b_{2,OK}}\dots$ on
the right hand side stands for  a generic linear combination of
acceptable  tensor fields in the form (\ref{form2''}) with simple
character $\vec{\kappa}_{simp}$, with no free indices in the
factors $\nabla\omega_1,\nabla\omega_2$ and where the factors
$T_{\omega_1},T_{\omega_2}$ are not problematic. Therefore, if we
can show the above equations, we are reduced to showing Lemma
\ref{firstclaimb} under the additional assumptions that
$H^2_{a,*}=H^2_{b,**}=H^2_{b,*}=\emptyset$.
\newline

{\it (Sketch of the) Proof of (\ref{zekkos1}), (\ref{zekkos2}), (\ref{zekkos3}):}
 (\ref{zekkos1}) follows by re-iterating the proof
 of the first claim of Lemma 4.10 in \cite{alexakis4}.\footnote{By the additional restrictions
imposed on the assumption of Lemma \ref{vanderbi} there is no
danger of falling under a ``forbidden case'' of Corollary
1 in \cite{alexakis4}.} (\ref{zekkos2}) follows by re-iterating the
 proof of the first claim of Lemma 4.10 in \cite{alexakis4}, but
rather than applying Corollary 1 \cite{alexakis4} in
that proof, we now apply Lemma \ref{firstclaim} (which
we have shown).\footnote{Observe that the assumption that Lemma \ref{vanderbi} does not
include ``forbidden cases'' ensures that we will not need to apply Lemma \ref{firstclaim}
in a ``forbidden case''.} Finally, the claim of (\ref{zekkos3})
{\it for the sublinear combination in $H^b_{2,*}$ where $T_{\omega_1}\ne T_{\omega_2}$}
follows by applying  Lemma \ref{addition2}.\footnote{In
this case there will be a factor $\nabla\omega_1$
or $\nabla\omega_2$ contracting against a non-special index;
therefore there is no danger of falling under a
``forbidden'' case of Lemma \ref{firstclaim}.} We can then show that
the remaining sublinear combination in
$\sum_{h\in H^b_{2,*}}\dots$ must vanish separately
(modulo a linear combination $\sum_{j\in J} \dots$) by just picking
out the sublinear combination in the hypothesis of Lemma \ref{secondclaimb}
where both factors $\nabla\omega_1,\nabla\omega_2$
are contracting against the same factor. $\Box$
\newline

 Now, under these additional assumptions that
 $H^2_{a,*}=H^2_{b,**}=H^2_{b,*}=\emptyset$, we will show our
 claim by distinguishing two cases: In case A there is a factor
 $\nabla^{(m)}R_{ijkl}$ in $\vec{\kappa}_{simp}$; in case $B$
 there is no such factor. An important note: We may now {\it use}
  Lemma \ref{firstclaim},
  which we have proven earlier in this section.
\newline

{\it Proof of Lemma \ref{firstclaimb} in case $A$.}
 \newline

\par We define the (set of) C-crucial factor(s) (which will necessarily be of
the form $\nabla^{(m)}R_{ijkl}$) as in the setting of Lemma
\ref{firstclaim}. Firstly  a mini-claim which only applies to the
case where the $C$-crucial factor is unique:

{\it Mini-claim, when the $C$-crucial factor is unique:} We then
consider the tensor fields $C^{h,i_{\pi+1}\dots i_{\alpha+1}}_g$,
$h\in H_2^a$ for which $\nabla\omega_2$ is contracting against the
 C-crucial factor. Notice that by our hypothesis that
$H^2_{a,*}=\emptyset$, we will have that $\nabla\omega_2$ is
contracting against a derivative index in the
C-crucial factor. Denote by $H^{a,+}_2\subset H^a_2$ the index set
of these tensor fields.

\par We observe that for each $h\in H^{a,+}_2$ we can now construct
a tensor field by erasing the index in the factor
$\nabla^{(m)}R_{ijkl}$ that contracts against the factor
$\nabla\omega_2$ and making the index in $\nabla\omega_2$ into a
free index
 ${}_{i_\beta}$. We denote this tensor field by $C^{h,i_{\pi+1}\dots
i_{\alpha+1}i_\beta}_{g}(\Omega_1,\dots ,\Omega_p,
\omega_1,\omega_2,\phi_1,\dots ,\phi_u)$. By the analogous
operation we obtain a tensor field $Switch[C^{h,i_{\pi+1}\dots
i_{\alpha+1}i_\beta}_{g}(\Omega_1,\dots ,\Omega_p,
\omega_1,\omega_2,\phi_1,\dots ,\phi_u)]$.

It follows that in the case where the C-crucial factor is
unique, for each $h\in H_2^{a,+}$:

\begin{equation}
\label{irineosgo}
\begin{split}
&X_{*}div_{i_{\pi+1}}\dots X_{*}div_{i_{\alpha+1}}\{
C^{h,i_{\pi+1}\dots i_{\alpha+1}}_{g}(\Omega_1,\dots ,\Omega_p,
\omega_1,\omega_2,\phi_1,\dots ,\phi_u)
\\& -Switch[C]^{h,i_{\pi+1}\dots i_{\alpha+1}}_{g}(\Omega_1,\dots
,\Omega_p, \omega_1,\omega_2,\phi_1,\dots ,\phi_u)\}=
 \\& X_{*}div_{i_{\pi+1}}\dots X_{*}div_{i_{\alpha+1}}X_{*}div_{i_\beta}\{
C^{h,i_{\pi+1}\dots i_{\alpha+1}i_\beta}_{g}(\Omega_1,\dots
,\Omega_p, \omega_1,\omega_2,\phi_1,\dots ,\phi_u)
\\& -Switch[C]^{h,i_{\pi+1}\dots
i_{\alpha+1}i_\beta}_{g}(\Omega_1,\dots ,\Omega_p,
\omega_1,\omega_2,\phi_1,\dots ,\phi_u)\}+
\\& \Sum_{r\in R} a_r X_{*}div_{i_{\pi+1}}\dots X_{*}div_{i_{\alpha+1}}\{ C^{r,i_{\pi+1}\dots
i_{\alpha+1}}_{g}(\Omega_1,\dots ,\Omega_p,
\omega_1,\omega_2,\phi_1,\dots ,\phi_u)
\\& -Switch[C]^{r,i_{\pi+1}\dots i_{\alpha+1}}(\Omega_1,\dots
,\Omega_p, \omega_1,\omega_2,\phi_1,\dots ,\phi_u)\}+
\\& \Sum_{j\in J} a_j C^j_{g}(\Omega_1,\dots ,\Omega_p,
\omega_1,\omega_2,\phi_1,\dots ,\phi_u),
\end{split}
\end{equation}
where each tensor field $C^{r,i_{\pi+1}\dots
i_{\alpha+1}}_{g}(\Omega_1,\dots ,\Omega_p,
\omega_1,\omega_2,\phi_1,\dots ,\phi_u)$ has the factor
$\nabla\omega_2$ contracting against some factor other than the
C-crucial factor.

\par But we observe that:
\begin{equation}
\label{seligaki}
\begin{split}
&X_{*}div_{i_{\pi+1}}\dots
X_{*}div_{i_{\alpha+1}}X_{*}div_{i_\beta}\{ C^{h,i_{\pi+1}\dots
i_{\alpha+1}i_\beta}_{g}(\Omega_1,\dots ,\Omega_p,
\omega_1,\omega_2,\phi_1,\dots ,\phi_u)
\\& -Switch[C]^{h,i_{\pi+1}\dots
i_{\alpha+1}i_\beta}_{g}(\Omega_1,\dots ,\Omega_p,
\omega_1,\omega_2,\phi_1,\dots ,\phi_u)\}=0.
\end{split}
\end{equation}
 Therefore, in the case $Set\ne\emptyset$ or $Set=\emptyset$
and $\sigma_1=1$, we have now reduced Lemma \ref{firstclaimb} to
 the case where $H^{a,+}_2=\emptyset$.
\newline

\par Now, (under the assumption that $H^{a,+}_2=\emptyset$ when the C-crucial factor is unique)
we consider the sublinear combination $Special$ in the equation
hypothesis of Lemma \ref{firstclaimb} that consists of complete
 contractions with $\nabla\omega_1$ contracting against the
 $C$-crucial factor while the factor $\nabla\omega_2$ is
 contracting against some other factor. (If $Set=\emptyset$
 and $\sigma_1>1$ $Special$ stands for the sublinear combination where
  $\nabla\omega_1$ is contracting against a generic
   $C$-crucial factor and $\nabla\omega_2$ is contracting against some other factor).
 In particular, for each $h\in
H^a_2$,
 since $H^{a,+}_2=\emptyset$ we see that the sublinear
 combination in
\begin{equation}
\label{laikhsofia} \begin{split} &\Sum_{h\in H^a_2} a_h
X_{*}div_{i_{\pi+1}}\dots X_{*}div_{i_{\alpha+1}}
\{C^{h,i_{\pi+1}\dots i_{\alpha+1}}_{g} (\Omega_1,\dots ,\Omega_p,
\omega_1,\omega_2,\phi_1,\dots ,\phi_u)
\\&-Switch[C]^{h,i_{\pi+1}\dots i_{\alpha+1}}_{g} (\Omega_1,\dots
,\Omega_p, \omega_1,\omega_2,\phi_1,\dots ,\phi_u)\}
\end{split}
\end{equation}
 that belongs to $Special$ is precisely:

$$\Sum_{h\in H^a_2} a_h X_{*}div_{i_{\pi+1}}\dots X_{*}div_{i_{\alpha}}
Hitdiv_{i_{\alpha+1}}C^{h,i_{\pi+1}\dots i_{\alpha+1}}_{g}
(\Omega_1,\dots ,\Omega_p, \omega_1,\omega_2,\phi_1,\dots
,\phi_u);$$ (in the case $Set=\emptyset$ and $\sigma_1>1$
 $Hitdiv_{i_{\alpha+1}}$ just means that $\nabla_{i_\gamma}$ can
 hit any factor $\nabla^{(m)}R_{ijkl}$ that is not contracting
 against $\nabla\omega_2$; recall that in the other cases it means that
 it must hit the unique C-crucial factor).

\par We also consider the tensor fields $C^{h,i_{\pi+1}\dots i_{\alpha+1}}$,
$Switch[C]^{h,i_{\pi+1}\dots i_{\alpha+1}}$, $h\in H^b_2$  for
which $\nabla\omega_1$ is contracting against the C-crucial factor
and $\nabla\omega_2$ is not (or, if there are multiple C-crucial
factors, where $\nabla\omega_1,\nabla\omega_2$ are contracting
against different C-crucial factors). For this proof, we index all
those tensor fields in $H^{b,\Psi}_2$ and we will denote them by
$C^{h,i_{\pi+1}\dots i_{\alpha+1}}_g$.

Thus we derive an equation:

\begin{equation}
\label{split}
\begin{split}
&\Sum_{h\in H^a_2} a_h X_{*}div_{i_{\pi+1}}\dots
X_{*}div_{i_{\alpha}} Hitdiv_{i_{\alpha+1}}C^{h,i_{\pi+1}\dots
i_{\alpha+1}}_{g} (\Omega_1,\dots ,\Omega_p,
\omega_1,\omega_2,\phi_1,\dots ,\phi_u)
\\&+\Sum_{h\in H^{b,\Psi}_2} a_h X_{*}div_{i_{\pi+1}}\dots X_{*}div_{i_{\alpha+1}}
C^{h,i_{\pi+1}\dots i_{\alpha+1}}_{g} (\Omega_1,\dots ,\Omega_p,
\omega_1,\omega_2,\phi_1,\dots ,\phi_u)=
\\&\Sum_{j\in J} a_j
C^j_{g}(\Omega_1,\dots ,\Omega_p, \omega_1,\omega_2,\phi_1,\dots
,\phi_u).
\end{split}
\end{equation}

We group up the vector fields on the left hand side according to
 their {\it weak $(u+2)$-characters}\footnote{See \cite{alexakis4} 
for a definition of this notion.} (defined by
 $\nabla\phi_1,\dots,\nabla\phi_u,\nabla\omega_1,\nabla\omega_2$).
(Recall that we started off with
 complete contractions with the same $u$-simple
 characters-so the only new information that we are taking into account
is what {\it type} of factor is $\nabla\omega_2$ contracting
 against).  We consider the set of
weak simple characters that we have obtained. We denote the set by
$\{\vec{\kappa}_1,\dots \vec{\kappa}_B\}$, and we respectively
have the index sets $H^{a,\vec{\kappa}_f}_2$ and
$H^{b,\vec{\kappa}_f}_2$.

 We will show our Lemma \ref{firstclaimb} by replacing the index set
$H^a_2$ by any $H^{a,\vec{\kappa}_f}_2$, $f\le B$.

It follows that for each $f\le B$:

\begin{equation}
\label{papar}
\begin{split}
&\Sum_{h\in H^{a,\vec{\kappa}_f}_2} a_h X_{*}div_{i_{\pi+1}}\dots
X_{*}div_{i_{\alpha}} Hitdiv_{i_{\alpha+1}}C^{h,i_{\pi+1}\dots
i_{\alpha+1}}_{g} \\&(\Omega_1,\dots ,\Omega_p,
\omega_1,\omega_2,\phi_1,\dots ,\phi_u)+
\\&\Sum_{h\in H^{b,\vec{\kappa}_f}_2} a_h X_{*}div_{i_{\pi+1}}\dots X_{*}div_{i_{\alpha+1}}
C^{h,i_{\pi+1}\dots i_{\alpha+1}}_{g} (\Omega_1,\dots ,\Omega_p,
\omega_1,\omega_2,\phi_1,\dots ,\phi_u)=
\\&\Sum_{j\in J} a_j
C^j_{g}(\Omega_1,\dots ,\Omega_p, \omega_1,\omega_2,\phi_1,\dots
,\phi_u),
\end{split}
\end{equation}
where the complete contractions $C^j_{g}$ have a $u$-simple
character that is subsequent to $\vec{\kappa}_{simp}$. We will
show our claim for each of
 the index sets $H^{b,\vec{\kappa}_f}_2$ separately.

\par Now, we treat the factors
$\nabla\omega_1,\nabla\omega_2$ in the above as factors
$\nabla\phi_{u+1},\nabla\phi_{u+2}$. We see that since
$H^b_{2,**}=H^2_{b,*}=H^2_{a,*}=\emptyset$, all the tensor fields in
the above have the same $(u+2)$-simple character.

\par Our claim (Lemma \ref{firstclaimb}) for 
the index set  $H^{a,\vec{\kappa}_f}_2$ then follows:
Firstly, apply the operator
$Erase_{\nabla\omega_1}[\dots]$ to (\ref{papar}).\footnote{See 
the relevant Lemma in the Appendix of \cite{alexakis1}.} 
We are then left with tensor fields (denote them by
$$C^{h,i_{\pi+1}\dots i_\alpha}_{g} (\Omega_1,\dots ,\Omega_p,
\omega_2,\phi_1,\dots ,\phi_u),h\in H^{a,\vec{\kappa}_f}_2,$$
 $$C^{h,i_{\pi+1}\dots i_{\alpha+1}}_{g} (\Omega_1,\dots ,\Omega_p,
\omega_2,\phi_1,\dots ,\phi_u), h\in H^{b,\vec{\kappa}_f}_2,$$
respectively) with the same $(u+1)$-simple character say
$\vec{\kappa}_{simp,f}$. We can then apply Corollary
1 from \cite{alexakis4} (since we have weight $-n+2k$, $k>0$ by virtue of
the eraser--notice that by weight considerations, since 
we started out with no ``bad'' tensor fields, there is no 
danger of falling under a ``fobidden case''). 
 to derive that there is a
linear combination of acceptable $\alpha$-tensor fields indexed in
$V$ below,  with $(u+1)$-simple character $\vec{\kappa}_{simp,f}$
so that:

\begin{equation}
\label{shost}
\begin{split}
&\Sum_{h\in H^{a,\vec{\kappa}_f}_2} a_h C^{h,i_{\pi+1}\dots
i_{\alpha}}_{g} (\Omega_1,\dots ,\Omega_p, \omega_2,\phi_1,\dots
,\phi_u)\nabla_{i_{\pi+1}}\upsilon\dots
\nabla_{i_{\alpha}}\upsilon-
\\&\Sum_{v\in V} a_v X_{*}div_{i_{\alpha+1}}C^{v,i_{\pi+1}\dots
i_{\alpha+1}}_{g} (\Omega_1,\dots ,\Omega_p, \omega_2,\phi_1,\dots
,\phi_u)\nabla_{i_{\pi+1}}\upsilon\dots
\nabla_{i_{\alpha}}\upsilon=
\\& \Sum_{j\in J} a_j C^{j,i_{\pi+1}\dots
i_{\alpha}}_{g} (\Omega_1,\dots ,\Omega_p, \omega_2,\phi_1,\dots
,\phi_u)\nabla_{i_{\pi+1}}\upsilon\dots
\nabla_{i_{\alpha}}\upsilon,
\end{split}
\end{equation}
where each complete contraction indexed in $J$ is $(u+1)$-subsequent to
$\vec{\kappa}_{simp,f}$. In this setting $X_{*}div_i$ just means that
 in addition to the restrictions imposed on $Xdiv_i$ we
  are not allowed to hit the factor $\nabla\omega_2$.

\par Then, if we multiply the above equation by an expression $\nabla_i\omega_1\nabla^i\upsilon$
and then anti-symmetrize the indices $a,b$ in the factors $\nabla_a\omega_1,\nabla_b\omega_2$ and
finally make all $\nabla\upsilon$s into $X_{+}div$s, we derive our
claim. $\Box$
\newline

{\it Proof of Lemma \ref{firstclaimb} in case B (when
$\sigma_1=0$).}
\newline

\par Our proof follows the same pattern as the proof of Lemma
\ref{firstclaim} in case B.

 We again define the ``measure'' of each factor in each
  tensor field $C^{h,i_{\pi+1}\dots i_{\alpha+1}}_g$ as in the
  proof of case B in Lemma \ref{firstclaim}.
  Again, let $M$ stand for the
maximum measure among all factors in all tensor fields
 $C^{h,i_{\pi+1}\dots i_{\alpha+1}}_g$, $h\in H^a_2$. We
denote by $H^{a,M}_2\subset H^a_2$ the index set of the tensor
fields for which some factor has  measure $M$.

We will further divide $H^{2,M}_a$ into subsets, $H^{2,M,k}_a,
k=1,\dots ,\sigma$, according to the factor which has measure $M$:
Firstly, we order the factors $S_{*}\nabla^{(\nu)}R_{ijkl},\dots
\nabla^{(p)}\Omega_h$ in $\vec{\kappa}_{simp}$, and label them
$T_1,\dots ,T_\sigma$ (observe each factor is well-defined in
$\vec{\kappa}_{simp}$, because we are in case B). We then say that
$h\in H^{a,M,1}_2$ if in $C^{h,i_{\pi+1}\dots i_\alpha}_{g}$ $T_1$
has measure $M$. We say say that $h\in H^{a,M,2}_2$ if in
$C^{h,i_{\pi+1}\dots i_\alpha}_{g}$ $T_2$ has measure $M$ and
$T_1$ has measure less than $M$, etc. We will then prove our claim
for each of the index sets $h\in H^{a,M,k}_2$.\footnote{Again we
observe that if we can prove this then Lemma \ref{firstclaimb} in
case B will follow by induction.} We arbitrarily pick a $k\le \sigma$
and show our claim for $\sum_{h\in H^{2,M,k}_a}\dots$.

\par For the purposes of this proof, we call the factor $T_k$ the
$D$-crucial factor (in this setting the $D$-crucial factor is
unique).

 Now, we pick out the subset
$H^{b,k,+}_2\subset H^b_2$, which is defined by the rule: $h\in
H^{b,k}_2$ if and only if $\nabla\omega_1$  is contracting against
the D-crucial factor $T_k$. We also pick out the subset
$H^{b,k,-}_2\subset H^b_2$, which is defined by the rule: $h\in
H^{b,k}_2$ if and only if $\nabla\omega_2$  is contracting against
the D-crucial factor $T_k$. Finally, we define
$H^{a,\tilde{}}_2\subset H^a_2$, $H^{a,-}_2\subset H^a_2$ to stand
for the index set of tensor fields for which $\nabla\omega_2$
 contracts against the  D-crucial factor.

\par Now, for each $h\in H^a_2$ we denote by
$$Hitdiv_{i_\gamma} C^{h,i_{\pi+1}\dots i_{\alpha+1}}_{g}
(\Omega_1,\dots ,\Omega_p,\omega_1,\omega_2,\phi_1,\dots,\phi_u)$$
the sublinear combination in $Xdiv_{i_\gamma} C^{h,i_{\pi+1}\dots
i_{\alpha+1}}_{g} (\Omega_1,\dots,\Omega_p,\omega_1,\omega_2,\phi_1,
\dots,\phi_u)$ that arises when $\nabla_{i_\gamma}$ hits the
D-crucial factor. It then follows that:

\begin{equation}
\label{hat''}
\begin{split}
&\Sum_{h\in H^{a}_2} a_h X_{*}div_{i_{\pi+1}}\dots
X_{*}div_{i_\alpha} Hitdiv_{i_\gamma} C^{h,i_{\pi+1}\dots
i_{\alpha+1}}_{g} (\Omega_1,\dots
,\Omega_p,\omega_1,\omega_2,\phi_1,\dots,\phi_u)
\\&-\Sum_{h\in H^{a,\tilde{}}_2} a_h X_{*}div_{i_{\pi+1}}\dots
X_{*}div_{i_{\alpha+1}}  Switch[C]^{h,i_{\pi+1}\dots
i_{\alpha+1}}_{g} (\Omega_1,\dots
,\Omega_p,\omega_1,\omega_2,\phi_1,\dots,\phi_u)
\\& + \Sum_{h\in H^{b,k,+}_2} a_h Xdiv_{i_{\pi+1}}\dots
Xdiv_{i_\gamma} C^{h,i_{\pi+1}\dots
i_{\alpha+1}}_{g}(\Omega_1,\dots ,
\Omega_p,\omega_1,\omega_2,\phi_1,\dots,\phi_u)
\\&-\Sum_{h\in H^{b,k,-}_2} a_h Xdiv_{i_{\pi+1}}\dots
Xdiv_{i_\gamma} Switch[C]^{h,i_{\pi+1}\dots
i_{\alpha+1}}_{g}(\Omega_1,\dots ,
\Omega_p,\omega_1,\omega_2,\phi_1,\dots,\phi_u)
\\&+\Sum_{j\in J} a_j C^j_{g}(\Omega_1,
\dots ,\Omega_p,\omega_1,\omega_2,\phi_1,\dots ,\phi_u),
\end{split}
\end{equation}
where each $C^j_{g}$ has the factor $\nabla \omega_1$ contracting
against the D-crucial factor and is simply subsequent to
$\vec{\kappa}_{simp}$.

  We now denote
the $(u+1)$-simple character (the one defined by
$\nabla\phi_1,\dots ,\nabla\omega_{1}$) of the tensor fields
$Hitdiv_{i_\gamma} C^{h,i_{\pi+1}\dots i_\alpha,i_\gamma}_{g}
(\Omega_1,\dots ,\Omega_p,\omega_1,\omega_2,\phi_1,\dots ,\phi_u)$
by $\vec{\kappa}'_{simp}$. (Observe that they all have the same
$(u+1)$-simple  character).

We observe that just applying Lemma \ref{obote}
 to (\ref{hat''}) (all tensor fields are acceptable and have the same
 simple character $\vec{\kappa}'_{simp}$--we treat $\nabla\omega_1$ as
  a factor $\nabla\phi_{u+1}$ and the factor $\nabla\omega_2$
 as a factor $\nabla Y$) and we then pick out the
 sublinear combination where there are $M$ factors
 $\nabla\upsilon$ contracting against $T_k$, we obtain an equation:

\begin{equation}
\label{ane3ignw''}
\begin{split}
&\Sum_{h\in H^{a,*,k}_2} a_h  Hitdiv_{i_\gamma}
C^{h,i_{\pi+1}\dots i_\alpha}_{g} (\Omega_1,\dots
,\Omega_p,\omega_1,\omega_2,\phi_1,\dots ,\phi_u) \nabla_{i_{\pi+1}}\upsilon\dots
\nabla_{i_\alpha}\upsilon+
\\&\Sum_{x\in X} a_x Xdiv_{i_{\alpha+1}}
C^{x,i_{\pi+1}\dots i_\alpha,i_{\alpha+1}}_{g} (\Omega_1,\dots
,\Omega_p,\omega_1,\omega_2,\phi_1,\dots ,\phi_u) \nabla_{i_{\pi+1}}\upsilon\dots
\nabla_{i_\alpha}\upsilon+
\\& \Sum_{j \in J} a_j C^{j,i_{\pi+1}\dots i_\alpha}_{g}
(\Omega_1,\dots ,\Omega_p,\omega_1,\omega_2,\phi_1,\dots ,\phi_u)
\nabla_{i_{\pi+1}}\upsilon\dots\nabla_{i_\alpha}\upsilon =0,
\end{split}
\end{equation}
where the tensor fields indexed in $X$ are acceptable and have a $(u+1)$-simple
character $\vec{\kappa}'_{simp}$ and each $C^j$ is simply
subsequent to $\vec{\kappa}'_{simp}$.

\par Now, observe that  {\it if} $M\ge \frac{3}{2}$ then we can apply
the Eraser (from the Appendix in \cite{alexakis1}) to $\nabla\omega_1$ and the index it is contracting
against in the $D$-crucial factor and derive our conclusion as in
case A.

\par The remaining cases are when $M=1,M=\frac{1}{2}$ and $M=0$.
The first one is easier, so we proceed to show our claim in that
case. The two subcases $M=\frac{1}{2},M=0$ will be discussed in
the next subsection.

In the case $M=1$, i.e. the $D$-crucial factor is of the form
$\nabla^{(p)}\Omega_h$, then we cannot derive our claim, because
if for some tensor fields in $X$ above we have $\nabla \omega_1$
contracting according to the pattern:
$\nabla_i\omega_1\nabla^{ij}\Omega_h\nabla_j\psi$, where
$\psi=\upsilon$ or $\psi=\phi_h$.  Therefore, in this setting, we
first apply the eraser twice to remove the expression
$\nabla^{(2)}_{ij}\Omega_h\nabla^i\psi\nabla^j\omega_1$ and then
 apply Corollary 2 from \cite{alexakis4}\footnote{Recall 
that we showed in \cite{alexakis4} that this is a Corollary 
of Lemma 4.6 in \cite{alexakis4}, which we have now shown.}
 to (\ref{ane3ignw''}) (observe that
  (\ref{ane3ignw''}) now falls under the inductive assumption of
Lemma 4.6 in \cite{alexakis4} since we have lowered the weight\footnote{There is no
danger of falling under a ``forbidden case'' 
of Lemma \ref{obote} by weight considerations 
 since we are assuming that none of the 
tensor fields of minimum rank in the assumption of Lemma
\ref{vanderbi} are ``bad''.}
 to obtain a new equation in the
form (\ref{ane3ignw''}), where each tensor field in $X$ has the
factor $\nabla \omega_1$ contracting against a factor
$\nabla^{(l)}\Omega_h$, $l\ge 3$. Then, applying the eraser as
explained, we derive our Lemma \ref{firstclaimb} in this case.

{\it The cases $M=\frac{1}{2}$, $M=0$:} Notice that in this case
we must have $\alpha=\pi$, by virtue of the the definition of maximal ``measure''
above. We will then prove our claim by proving a more
general claim by induction, in the next subsection. $\Box$

\subsection{The remaining cases of Lemma \ref{firstclaimb}.}

\par We prove our claim in this case by an induction.
 In order to give a detailed proof,
 we will re-state our Lemma hypothesis in this case
(with a slight change of notation):
\newline

{\it The hypothesis of the remaining cases of Lemma
\ref{firstclaimb}:} We are assuming an equation:

\begin{equation}
\label{klm}
\begin{split}
&\sum_{x\in X_a} a_x X_{*}div_{i_1}C^{x,i_1}_g(\Omega_1,\dots,
\Omega_p,\phi_1,\dots,\phi_u,[\omega_1,\omega_2])+
\\&\sum_{x\in X_b} a_x X_{*}div_{i_1}C^{x,i_1}_g(\Omega_1,\dots,
\Omega_p,\phi_1,\dots,\phi_u,[\omega_1,\omega_2])+
\\&\sum_{j\in J} a_j C^{j}_g(\Omega_1,\dots,\Omega_p,\phi_1,\dots,\phi_u)=0,
\end{split}
\end{equation}
which holds modulo complete contractions of length
$\ge\sigma+u+3$ ($\sigma\ge 3$--here $\sigma$ stands for $u+p$--see the next equation).
 We denote the weight of the complete contractions
in the above by $-K$.
 The tensor fields in the above equation are each
 in the form:

\begin{equation}
\label{presevo}
\begin{split}
&pcontr(S_{*}\nabla^{(\nu_1)}R_{x_1jkl}\otimes\dots\otimes
S_{*}\nabla^{(\nu_u)}R_{x_zj'k'l'}\otimes
\\&\nabla^{(a_1)}\Omega_1\otimes\dots\nabla^{(a_p)}\Omega_p\otimes [\nabla\omega_1\otimes\nabla\omega_2]
\\&\otimes\nabla^{x_1}\tilde{\phi}_1\otimes\dots\otimes\nabla^{x_u}\tilde{\phi}_z).
\end{split}
\end{equation}
We recall that the $u$-simple character of the above has been
denoted by $\vec{\kappa}_{simp}$. Recall that we are now assuming
that all the
 factors  $\nabla^{(a_i)}\Omega_x$
 in $\vec{\kappa}_{simp}$ are acceptable.\footnote{meaning that each $a_i\ge 2$.}
The complete contractions indexed in $J$ in (\ref{klm}) are simply
subsequent to $\vec{\kappa}_{simp}$. We also recall that
$X_{*}div_i$ stands
 for the sublinear combination in $Xdiv_i$
 where $\nabla_i$ is not allowed to hit either
of the factors $\nabla\omega_1,\nabla\omega_2$.

\par We recall that the tensor fields indexed in $X_a$
 have the free index ${}_{i_1}$ belonging to the factor
$\nabla\omega_1$. The tensor fields indexed in $X_b$ have the free
index ${}_{i_1}$ {\it not} belonging to any of the factors
$\nabla\omega_1,\nabla\omega_2$.

\par We recall the key assumption that for each of the tensor
 fields indexed in $X_a$, there is at least one removable index in each
tensor field \\$C^{x,i_1}_g(\Omega_1,\dots,
\Omega_p,\phi_1,\dots,\phi_u,[\omega_1,\omega_2])$,
 $x\in X_a$.\footnote{Recall the definition of
 a ``removable'' index from Definition \ref{proextremovable}.}
\newline

\par In order to complete our proof of Lemma \ref{firstclaimb}, we
will show that we can write:

\begin{equation}
\label{proarendt}
\begin{split}
&\sum_{x\in X_a} a_x C^{x,i_1}_g(\Omega_1,\dots,
\Omega_p,\phi_1,\dots,\phi_u,[\omega_1,\omega_2])\nabla_{i_1}\upsilon=
\\&\sum_{x\in X'} a_x X_{*}div_{i_2}\dots X_{*}div_{i_a}C^{x,i_1\dots i_a}_g(\Omega_1,\dots,
\Omega_p,\phi_1,\dots,\phi_u,[\omega_1,\omega_2])\nabla_{i_1}\upsilon+
\\&\sum_{j\in J} a_j C^{j}_g(\Omega_1,\dots,\Omega_p,\phi_1,\dots,\phi_u),
\end{split}
\end{equation}
where the tensor fields indexed in $X'$ are acceptable in the form
(\ref{presevo}), each with  rank $a\ge 2$. Note that this will imply
the remaining cases of  Lemma \ref{firstclaimb}, completing the proof of Lemma \ref{vanderbi}.

\par We recall that we are proving this claim when the assumption
(\ref{klm}) formally falls under our inductive assumption of
Proposition \ref{giade} (if we formally treat
$\nabla\omega_1,\nabla\omega_2$ as factors
 $\nabla\phi_{z+1},\nabla\phi_{z+2}$).

\par We will prove (\ref{proarendt}) by inductively proving
 a more general statement.
The more general statement is as follows:
\newline

{\it The general statement:}
\newline

{\it Assumptions:}  We consider vector fields $C^{\zeta,i_1}_g
(\Omega_1,\dots,\Omega_b,\phi_1,\dots,\phi_v,
Y,\psi_1,\dots,\psi_\tau)$,
\\ $C^{\zeta,i_1}_g(\Omega_1,\dots,\Omega_b,\phi_1,\dots,\phi_v,
[\chi_1,\chi_2],\psi_1,\dots,\psi_\tau)$ in
 the following forms, respectively:

\begin{equation}
\label{proherodium}
\begin{split}
&pcontr(S_{*}\nabla^{(\nu_1)}R_{x_1jkl}\otimes\dots\otimes
S_{*}\nabla^{(\nu_v)}R_{x_vj'k'l'}\otimes
\\&\nabla^{(a_1)}\Omega_1\otimes\dots\nabla^{(a_b)}
\Omega_b\otimes \nabla Y\otimes\nabla
\psi_1\otimes\dots\otimes\nabla\psi_\tau
\\&\otimes\nabla^{x_1}\tilde{\phi}_1\otimes\dots\otimes\nabla^{x_v}\tilde{\phi}_v),
\end{split}
\end{equation}

\begin{equation}
\label{herodium}
\begin{split}
&pcontr(S_{*}\nabla^{(\nu_1)}R_{x_1jkl}\otimes\dots\otimes
S_{*}\nabla^{(\nu_v)}R_{x_vj'k'l'}\otimes
\\&\nabla^{(a_1)}\Omega_1\otimes\dots\nabla^{(a_b)}
\Omega_b\otimes [\nabla\chi_1\otimes\nabla\chi_2]\otimes\nabla
\psi_1\otimes\dots\otimes\nabla\psi_\tau
\\&\otimes\nabla^{x_1}\tilde{\phi}_1\otimes\dots\otimes\nabla^{x_v}\tilde{\phi}_v),
\end{split}
\end{equation}
for which the weight is $-W+1, W\le K$. We also assume $v+b\ge 2$.
{\it Note:} the bracket $[\dots]$ stands for the
 anti-symmetrization of the indices ${}_a,{}_b$ in the expression $\nabla_a\omega_1\nabla_b\omega_2$.

 We assume (respectively) the equations:

\begin{equation}
\label{procareless}
\begin{split}
&\sum_{\zeta\in Z_a} a_\zeta X_{*}div_{i_1}
C^{\zeta,i_1}_g(\Omega_1,\dots,\Omega_b,\phi_1,\dots,\phi_v,
Y,\psi_1,\dots,\psi_\tau)+
\\&\sum_{\zeta\in \overline{Z}_a} a_\zeta X_{*}div_{i_1}\dots X_{*}div_{i_\gamma}
C^{\zeta,i_1\dots
i_\gamma}_g(\Omega_1,\dots,\Omega_b,\phi_1,\dots,\phi_v,
Y,\psi_1,\dots,\psi_\tau)+
\\&\sum_{\zeta\in Z_b} a_\zeta X_{*}div_{i_1}
C^{\zeta,i_1}_g(\Omega_1,\dots,\Omega_b,\phi_1,\dots,\phi_v,
Y,\psi_1,\dots,\psi_\tau)+
\\&\sum_{j\in J} a_j C^{j}_g(\Omega_1,\dots,\Omega_b,\phi_1,\dots,\phi_v,
Y,\psi_1,\dots,\psi_\tau)=0,
\end{split}
\end{equation}

\begin{equation}
\label{careless}
\begin{split}
&\sum_{\zeta\in Z_a} a_\zeta X_{*}div_{i_1}
C^{\zeta,i_1}_g(\Omega_1,\dots,\Omega_b,\phi_1,\dots,\phi_v,
[\chi_1,\chi_2],\psi_1,\dots,\psi_\tau)+
\\&\sum_{\zeta\in \overline{Z}_a} a_\zeta X_{*}div_{i_1}\dots X_{*}div_{i_\gamma}
C^{\zeta,i_1\dots
i_\gamma}_g(\Omega_1,\dots,\Omega_b,\phi_1,\dots,\phi_v,
[\chi_1,\chi_2],\psi_1,\dots,\psi_\tau)+
\\&\sum_{\zeta\in Z_b} a_\zeta X_{*}div_{i_1}
C^{\zeta,i_1}_g(\Omega_1,\dots,\Omega_b,\phi_1,\dots,\phi_v,
[\chi_1,\chi_2],\psi_1,\dots,\psi_\tau)+
\\&\sum_{j\in J} a_j C^{j}_g(\Omega_1,\dots,\Omega_b,\phi_1,\dots,\phi_v,
[\chi_1,\chi_2],\psi_1,\dots,\psi_\tau)=0,
\end{split}
\end{equation}
which holds modulo complete contractions of length $\ge v+b+\tau+3$.

\par The tensor fields indexed in $Z_a$ are assumed to have
a free index in one of the factors $\nabla
Y,\nabla\psi_1,\dots,\nabla\psi_\tau$, or one of the factors
 $\nabla\chi_1,\nabla\chi_2,\nabla\psi_1,\dots,\nabla\psi_\tau$, respectively.
 The tensor fields indexed in $\overline{Z}_a$ have rank $\gamma\ge 2$ and
 all their free indices belong to the factors $\nabla
Y,\nabla\psi_1,\dots,\nabla\psi_\tau$, or the factors
 $\nabla\chi_1,\nabla\chi_2,\nabla\psi_1,\dots,\nabla\psi_\tau$,
 respectively.
The tensor fields indexed in $Z_b$ have the property that
${}_{i_1}$ does not belong to any of the factors $\nabla
Y,\nabla\psi_1,\dots,\nabla\psi_\tau$, $\nabla\chi_1,
\nabla\chi_2,\nabla\psi_1,\dots,\nabla\psi_\tau$, respectively. We
furthermore assume that for the tensor fields indexed in
$Z_a\bigcup Z_b\bigcup \overline{Z}_a$, none of the factors
 $\nabla\psi_1,\dots,\nabla\psi_\tau$ are contracting
 against a special index in any factor
 $S_{*}\nabla^{(\nu)}R_{ijkl}$ and none of them are contracting
 against the rightmost index in each $\nabla^{(a_h)}\Omega_h$
  (we will refer to this property as the $@$-property).
We assume that $v+b\ge 2$, and furthermore if $v+b=2$ then for
each $\zeta\in Z_a\bigcup Z_b$, the factor(s) $\nabla Y$ (or
$\nabla\chi_1,\nabla\chi_2$) are also not contracting against a
special index in any $S_{*}\nabla^{(\nu)}R_{ijkl}$ and are not
contracting against the rightmost index in any
$\nabla^{(a_h)}\Omega_h$. Finally (and importantly) we assume that
for the tensor fields
 indexed in $Z_a$, there is at least one removable index in
 each $C^{\zeta,i_1}$. (In this setting, for a tensor field indexed in $Z_a$,
  a ``removable'' index is either a non-special index in
   a factor $S_{*}\nabla^{(\nu)}R_{ijkl}$, with $\nu>0$ or an
  index in a factor $\nabla^{(B)}\Omega_h$, $B\ge 3$).

{\it Convention:} In this subsection only, for tensor fields in the
forms (\ref{procareless}), (\ref{careless}) we say then an index
is {\it special} if it is one of the indices ${}_k,{}_l$ in a
factor $S_{*}\nabla^{(\nu)}R_{ijkl}$ (this is the usual
convention), {\it or} if it is an index in a factor
$\nabla^{(B)}_{r_1\dots r_B}\Omega_h$ {\it for which all the other
indices are contracting against factors
$\nabla\psi_1,\dots,\nabla\psi_\tau$}.

\par All tensor fields in (\ref{procareless}), (\ref{careless}) have a given
$v$-simple character $\overline{\kappa}_{simp}$. The complete
contractions indexed in $J$
 are assumed to have a weak $v$-character
 $Weak(\overline{\kappa}_{simp})$ and to be simply
 subsequent to $\overline{\kappa}_{simp}$. Here $X_{*}div_i$
 stands for the sublinear combination in $Xdiv_i$ where $\nabla_i$
 is not allowed to hit any of the factors $\nabla Y,
\nabla\psi_1,\dots,\nabla\psi_\tau$ or
$\nabla\chi_1,\nabla\chi_2,\nabla\psi_1,\dots,\nabla\psi_\tau$,
 respectively.
\newline

{\it The Claims of the general statement:} We claim that under the
assumption (\ref{careless}), there exists a linear combination of
acceptable 2-tensor fields in the form (\ref{proherodium}),
(\ref{herodium}) respectively (indexed in $W$ below), for which
 the $@$-property is satisfied,
 so that (respectively):

\begin{equation}
\label{procareless2}
\begin{split}
&\sum_{\zeta\in Z_a} a_\zeta
C^{\zeta,i_1}_g(\Omega_1,\dots,\Omega_b,\phi_1,\dots,\phi_v,
Y,\psi_1,\dots,\psi_\tau)\nabla_{i_1}\upsilon-
\\&\sum_{w\in W} a_w X_{*}div_{i_2}
C^{w,i_1i_2}_g(\Omega_1,\dots,\Omega_b,\phi_1,\dots,\phi_v,
Y,\psi_1,\dots,\psi_\tau)\nabla_{i_1}\upsilon+
\\&\sum_{j\in J} a_j C^{j,i_1}_g(\Omega_1,\dots,\Omega_b,\phi_1,\dots,\phi_v,
Y,\psi_1,\dots,\psi_\tau)\nabla_{i_1}\upsilon=0,
\end{split}
\end{equation}

\begin{equation}
\label{careless2}
\begin{split}
&\sum_{\zeta\in Z_a} a_\zeta
C^{\zeta,i_1}_g(\Omega_1,\dots,\Omega_b,\phi_1,\dots,\phi_v,
[\chi_1,\chi_2],\psi_1,\dots,\psi_\tau)\nabla_{i_1}\upsilon+
\\&\sum_{w\in W} a_w X_{*}div_{i_2}
C^{w,i_1i_2}_g(\Omega_1,\dots,\Omega_b,\phi_1,\dots,\phi_v,
[\chi_1,\chi_2],\psi_1,\dots,\psi_\tau)\nabla_{i_1}\upsilon+
\\&\sum_{j\in J} a_j C^{j,i_1}_g(\Omega_1,\dots,\Omega_b,\phi_1,\dots,\phi_v,
[\chi_1,\chi_2],\psi_1,\dots,\psi_\tau)\nabla_{i_1}\upsilon=0.
\end{split}
\end{equation}

\par We observe that when $\tau=0$ and $v+b\ge 3$, (\ref{careless2}) coincides
with (\ref{proarendt}).\footnote{Also, the assumption of existence 
of a non-removable index coincides with
 the corresponding assumption of Lemma \ref{vanderbi}.}
Therefore, if we can prove this general
statement, we will have shown Lemma \ref{firstclaimb} in full
generality, thus also completing the proof of Lemma \ref{vanderbi}.
\newline

\par We also have a further claim, when we assume
(\ref{procareless}), (\ref{careless}) with $v+b=2$.

In that case, we also claim that we can write:

\begin{equation}
\label{balladur}
\begin{split}
&X_{+}div_{i_1}\sum_{\zeta\in Z_a\bigcup Z_b\bigcup \overline{Z}_a} a_\zeta
C^{\zeta,i_1}_g(\Omega_1,\dots,\Omega_b,\phi_1,\dots,\phi_v,
Y,\psi_1,\dots,\psi_\tau)=
\\&\sum_{q\in Q} a_q X_{+}div_{i_1}
C^{q,i_1}_g(\Omega_1,\dots,\Omega_b,\phi_1,\dots,\phi_v,
Y,\psi_1,\dots,\psi_\tau)+
\\&\sum_{j\in J} a_j C^{j}_g(\Omega_1,\dots,\Omega_b,\phi_1,\dots,\phi_v,
Y,\psi_1,\dots,\psi_\tau),
\end{split}
\end{equation}

\begin{equation}
\label{balladur'}
\begin{split}
&X_{+}div_{i_1}\sum_{\zeta\in Z_a\bigcup \overline{Z}_a\bigcup Z_b} a_\zeta
C^{\zeta,i_1}_g(\Omega_1,\dots,\Omega_b,\phi_1,\dots,\phi_v,
[\chi_1,\chi_2],\psi_1,\dots,\psi_\tau)=
\\&\sum_{q\in Q} a_q X_{+}div_{i_1}
C^{q,i_1}_g(\Omega_1,\dots,\Omega_b,\phi_1,\dots,\phi_v,
[\chi_1,\chi_2],\psi_1,\dots,\psi_\tau)+
\\&\sum_{j\in J} a_j  C^{j}_g(\Omega_1,\dots,\Omega_b,\phi_1,\dots,\phi_v,
[\chi_1,\chi_2],\psi_1,\dots,\psi_\tau),
\end{split}
\end{equation}
where the tensor fields indexed in $Q$ are in the same form as
(\ref{proherodium}) or (\ref{herodium}) respectively, but have a
factor (expression)
 $\nabla^{(2)}Y$ or $\nabla^{(2)}_{a[i}\omega_1\nabla_{j]}\omega_2$, respectively, and
satisfy all the other properties of the tensor fields in $Z_a$.
 \newline

{\it Consequence of (\ref{procareless2}), (\ref{careless2}) when $v+b\ge 3$:} We
here codify a conclusion one can derive from (\ref{procareless2}),
(\ref{careless2}). This implication will be useful further down in
this subsection. We see that by making the factors
$\nabla\upsilon$ into $X_{*} div$'s in (\ref{procareless}),
(\ref{careless}) and replacing into (\ref{procareless2}),
(\ref{careless2}), we obtain new equations:

\begin{equation}
\label{batis1}
\begin{split}
&\sum_{\zeta\in Z'_a} a_\zeta X_{*}div_{i_1}
C^{\zeta,i_1}_g(\Omega_1,\dots,\Omega_b,\phi_1,\dots,\phi_v,
Y,\psi_1,\dots,\psi_\tau)+
\\&\sum_{\zeta\in Z_b} a_\zeta X_{*}div_{i_1}
C^{\zeta,i_1}_g(\Omega_1,\dots,\Omega_b,\phi_1,\dots,\phi_v,
Y,\psi_1,\dots,\psi_\tau)+
\\&\sum_{j\in J} a_j C^{j}_g(\Omega_1,\dots,\Omega_b,\phi_1,\dots,\phi_v,
Y,\psi_1,\dots,\psi_\tau)=0,
\end{split}
\end{equation}

\begin{equation}
\label{batis2}
\begin{split}
&\sum_{\zeta\in Z'_a} a_\zeta X_{*}div_{i_1}
C^{\zeta,i_1}_g(\Omega_1,\dots,\Omega_b,\phi_1,\dots,\phi_v,
[\chi_1,\chi_2],\psi_1,\dots,\psi_\tau)+
\\&\sum_{\zeta\in Z_b} a_\zeta X_{*}div_{i_1}
C^{\zeta,i_1}_g(\Omega_1,\dots,\Omega_b,\phi_1,\dots,\phi_v,
[\chi_1,\chi_2],\psi_1,\dots,\psi_\tau)+
\\&\sum_{j\in J} a_j C^{j}_g(\Omega_1,\dots,\Omega_b,\phi_1,\dots,\phi_v,
[\chi_1,\chi_2],\psi_1,\dots,\psi_\tau)=0,
\end{split}
\end{equation}
where here the  tensor fields indexed in $Z'_a$ are like the
tensor fields indexed in $Z_a$ in (\ref{procareless}),
(\ref{careless}) but have the additional feature that no free
index belongs to the factor $\nabla\psi_1$ (and all the other
assumptions of equations (\ref{procareless}), (\ref{careless})
continue to hold).

\par We then claim that we can derive new equations:

\begin{equation}
\label{corbatis1}
\begin{split}
&\sum_{\zeta\in Z'_a} a_\zeta X_{+}div_{i_1}
C^{\zeta,i_1}_g(\Omega_1,\dots,\Omega_b,\phi_1,\dots,\phi_v,
Y,\psi_1,\dots,\psi_\tau)+
\\&\sum_{\zeta\in Z_b} a_\zeta X_{+}div_{i_1}
C^{\zeta,i_1}_g(\Omega_1,\dots,\Omega_b,\phi_1,\dots,\phi_v,
Y,\psi_1,\dots,\psi_\tau)=
\\&\sum_{q\in Q} a_q X_{+}div_{i_1}
C^{q,i_1}_g(\Omega_1,\dots,\Omega_b,\phi_1,\dots,\phi_v,
Y,\psi_1,\dots,\psi_\tau)+
\\&\sum_{j\in J} a_j C^{j}_g(\Omega_1,\dots,\Omega_b,\phi_1,\dots,\phi_v,
Y,\psi_1,\dots,\psi_\tau),
\end{split}
\end{equation}

\begin{equation}
\label{corbatis2}
\begin{split}
&\sum_{\zeta\in Z'_a} a_\zeta X_{+}div_{i_1}
C^{\zeta,i_1}_g(\Omega_1,\dots,\Omega_b,\phi_1,\dots,\phi_v,
[\chi_1,\chi_2],\psi_1,\dots,\psi_\tau)+
\\&\sum_{\zeta\in Z_b} a_\zeta X_{+}div_{i_1}
C^{\zeta,i_1}_g(\Omega_1,\dots,\Omega_b,\phi_1,\dots,\phi_v,
[\chi_1,\chi_2],\psi_1,\dots,\psi_\tau)=
\\&\sum_{q\in Q} a_q X_{+}div_{i_1}
C^{q,i_1}_g(\Omega_1,\dots,\Omega_b,\phi_1,\dots,\phi_v,
[\chi_1,\chi_2],\psi_1,\dots,\psi_\tau)+
\\&\sum_{j\in J} a_j C^{j}_g(\Omega_1,\dots,\Omega_b,\phi_1,\dots,\phi_v,
[\chi_1,\chi_2],\psi_1,\dots,\psi_\tau),
\end{split}
\end{equation}
where here $X_{+}div_i$ stands for the  sublinear combination in
$Xdiv_i$ where $\nabla_i$ {\it is} allowed to hit the factor
$\nabla Y$ or $\nabla\chi_1$ (respectively), but not the
factors $\nabla\psi_1,\dots,\nabla\phi_\tau,(\nabla\chi_2)$. Furthermore, the
linear combinations indexed  in $Q$ stand for generic linear
combinations of vector fields in the form (\ref{proherodium}) or
(\ref{herodium}), only with the expressions $\nabla Y$ or
$\nabla_{[a}\omega_1\nabla_{b]}\omega_2$ replaced by expressions
$\nabla^{(2)}Y$, $\nabla^{(2)}_{c[a}\omega_1\nabla_{b]}\omega_2$.
\newline

{\it Proof that (\ref{corbatis1}), (\ref{corbatis2}) follow from
(\ref{procareless2}), (\ref{careless2}):} We prove the above by an
induction. We will firstly subdivide $Z'_a,Z_b$ into subsets as
follows: $\zeta\in Z'_{a,@}$ or $\zeta\in Z_{b,@}$ if the factor
$\nabla Y$ (or one of the factors $\nabla\chi_1,\nabla\chi_2$) is
contracting against a special index in the same factor against
which $\nabla\psi_1$ is contracting.

\par Now, if $Z'_{a,@}\bigcup Z_{b,@}\ne \emptyset$ our inductive
statement will be the following:

We inductively assume that we can write:

\begin{equation}
\label{olmert1}
\begin{split}
&\sum_{\zeta\in Z'_{a,@}} a_\zeta X_{+}div_{i_1}\dots
X_{+}div_{i_\gamma} C^{\zeta,i_1\dots
i_\gamma}_g(\Omega_1,\dots,\Omega_b,\phi_1,\dots,\phi_v,
Y,\psi_1,\dots,\psi_\tau)=
\\&\sum_{\zeta\in Z_{b,@}} a_\zeta X_{+}div_{i_1}
C^{\zeta,i_1}_g(\Omega_1,\dots,\Omega_b,\phi_1,\dots,\phi_v,
Y,\psi_1,\dots,\psi_\tau)+
\\&\sum_{t\in T^k} a_t X_{+}div_{i_1}\dots X_{+}div_{i_k}
C^{t,i_1\dots i_k}_g(\Omega_1,\dots,\Omega_b,\phi_1,\dots,\phi_v,
Y,\psi_1,\dots,\psi_\tau)+
\\&\sum_{\zeta\in Z'_{a,No@}} a_\zeta X_{+}div_{i_1}\dots
X_{+}div_{i_\gamma} C^{\zeta,i_1\dots
i_\gamma}_g(\Omega_1,\dots,\Omega_b,\phi_1,\dots,\phi_v,
Y,\psi_1,\dots,\psi_\tau)+
\\&\sum_{q\in Q} a_q X_{+}div_{i_1}
C^{q,i_1}_g(\Omega_1,\dots,\Omega_b,\phi_1,\dots,\phi_v,
Y,\psi_1,\dots,\psi_\tau)+
\\&\sum_{j\in J} a_j C^{j}_g(\Omega_1,\dots,\Omega_b,\phi_1,\dots,\phi_v,
Y,\psi_1,\dots,\psi_\tau),
\end{split}
\end{equation}
and

\begin{equation}
\label{olmert2}
\begin{split}
&\sum_{\zeta\in Z'_{a,@}} a_\zeta X_{+}div_{i_1}
C^{\zeta,i_1}_g(\Omega_1,\dots,\Omega_b,\phi_1,\dots,\phi_v,
[\chi_1,\chi_2],\psi_1,\dots,\psi_\tau)=
\\&\sum_{\zeta\in Z_{b,@}} a_\zeta X_{+}div_{i_1}
C^{\zeta,i_1}_g(\Omega_1,\dots,\Omega_b,\phi_1,\dots,\phi_v,
[\chi_1,\chi_2],\psi_1,\dots,\psi_\tau)+
\\&\sum_{t\in T^k} a_t X_{+}div_{i_1}\dots X_{+}div_{i_k}
C^{t,i_1\dots i_k}_g(\Omega_1,\dots,\Omega_b,\phi_1,\dots,\phi_v,
[\chi_1,\chi_2],\psi_1,\dots,\psi_\tau)+
\\&\sum_{\zeta\in Z'_{a,No@}} a_\zeta X_{+}div_{i_1}\dots
X_{+}div_{i_\gamma} C^{\zeta,i_1\dots
i_\gamma}_g(\Omega_1,\dots,\Omega_b,\phi_1,\dots,\phi_v,
[\chi_1,\chi_2],\psi_1,\dots,\psi_\tau)+
\\&\sum_{q\in Q} a_q X_{+}div_{i_1}
C^{q,i_1}_g(\Omega_1,\dots,\Omega_b,\phi_1,\dots,\phi_v,
[\chi_1,\chi_2],\psi_1,\dots,\psi_\tau)+
\\&\sum_{j\in J} a_j C^{j}_g(\Omega_1,\dots,\Omega_b,\phi_1,\dots,\phi_v,
[\chi_1,\chi_2],\psi_1,\dots,\psi_\tau),
\end{split}
\end{equation}
where the tensor fields indexed in $T^k$ have all the properties
of the tensor fields indexed in $Z'_{a,@}$ (in particular the index in
$\nabla\psi_1$ is not free) and in addition have rank $k$. The
tensor fields indexed in $Z'_{a,No@}$ in the RHS have all the
regular features of the terms indexed in $Z'_a$ (in particular rank $\gamma\ge
1$ and the factor $\nabla\psi_1$ does not contain a free index)
and in addition none of the factors $\nabla Y$ (or
$\nabla\chi_1,\nabla\chi_2$) are contracting against a special
index.

\par Our inductive claim is that we can write:

\begin{equation}
\label{proolmert1'}
\begin{split}
&\sum_{\zeta\in Z'_a} a_\zeta X_{+}div_{i_1}
C^{\zeta,i_1}_g(\Omega_1,\dots,\Omega_b,\phi_1,\dots,\phi_v,
Y,\psi_1,\dots,\psi_\tau)=
\\&\sum_{\zeta\in Z_b} a_\zeta X_{+}div_{i_1}
C^{\zeta,i_1}_g(\Omega_1,\dots,\Omega_b,\phi_1,\dots,\phi_v,
Y,\psi_1,\dots,\psi_\tau)+
\\&\sum_{t\in T^{k+1}} a_t X_{+}div_{i_1}\dots X_{+}div_{i_{k+1}}
C^{t,i_1\dots i_{k+1}}_g(\Omega_1,\dots,\Omega_b,\phi_1,\dots,\phi_v,
Y,\psi_1,\dots,\psi_\tau)+
\\&\sum_{\zeta\in Z'_{a,No@}} a_\zeta X_{+}div_{i_1}\dots
X_{+}div_{i_\gamma} C^{\zeta,i_1\dots
i_\gamma}_g(\Omega_1,\dots,\Omega_b,\phi_1,\dots,\phi_v,
Y,\psi_1,\dots,\psi_\tau)+
\\&\sum_{q\in Q} a_q X_{+}div_{i_1}
C^{q,i_1}_g(\Omega_1,\dots,\Omega_b,\phi_1,\dots,\phi_v,
Y,\psi_1,\dots,\psi_\tau)+
\\&\sum_{j\in J} a_j C^{j}_g(\Omega_1,\dots,\Omega_b,\phi_1,\dots,\phi_v,
Y,\psi_1,\dots,\psi_\tau),
\end{split}
\end{equation}

\begin{equation}
\label{proolmert2'}
\begin{split}
&\sum_{\zeta\in Z'_a} a_\zeta X_{+}div_{i_1}
C^{\zeta,i_1}_g(\Omega_1,\dots,\Omega_b,\phi_1,\dots,\phi_v,
[\chi_1,\chi_2],\psi_1,\dots,\psi_\tau)=
\\&\sum_{\zeta\in Z_b} a_\zeta X_{+}div_{i_1}
C^{\zeta,i_1}_g(\Omega_1,\dots,\Omega_b,\phi_1,\dots,\phi_v,
[\chi_1,\chi_2],\psi_1,\dots,\psi_\tau)+
\\&\sum_{t\in T^{k+1}} a_t X_{+}div_{i_1}\dots X_{+}div_{i_k}
C^{t,i_1\dots
i_{k+1}}_g(\Omega_1,\dots,\Omega_b,\phi_1,\dots,\phi_v,
[\chi_1,\chi_2],\psi_1,\dots,\psi_\tau)+
\\&\sum_{\zeta\in Z'_{a,No@}} a_\zeta X_{+}div_{i_1}\dots
X_{+}div_{i_\gamma} C^{\zeta,i_1\dots
i_\gamma}_g(\Omega_1,\dots,\Omega_b,\phi_1,\dots,\phi_v,
[\chi_1,\chi_2],\psi_1,\dots,\psi_\tau)
\\&+\sum_{q\in Q} a_q X_{+}div_{i_1}
C^{q,i_1}_g(\Omega_1,\dots,\Omega_b,\phi_1,\dots,\phi_v,
[\chi_1,\chi_2],\psi_1,\dots,\psi_\tau)+
\\&\sum_{j\in J} a_j C^{j}_g(\Omega_1,\dots,\Omega_b,\phi_1,\dots,\phi_v,
[\chi_1,\chi_2],\psi_1,\dots,\psi_\tau)=0.
\end{split}
\end{equation}

\par We will derive (\ref{proolmert1'}), (\ref{proolmert2'}) momentarily.
For now, we observe that by iterative repetition of the above inductive step  we
are reduced to showing (\ref{corbatis1}), (\ref{corbatis2}) under
the additional assumption that $Z'_{a,@}=\emptyset$.

\par Under that assumption, we denote by $Z_{b,@}\subset Z_b$ the
index set of vector fields for which the factor $\nabla Y$ (or one
of the factors $\nabla\chi_1,\nabla\chi_2$) is contracting against
a special index. We will then assume that we can write:

\begin{equation}
\label{proolmert1b'}
\begin{split}
&\sum_{\zeta\in Z_{b,@}} a_\zeta X_{+}div_{i_1}
C^{\zeta,i_1}_g(\Omega_1,\dots,\Omega_b,\phi_1,\dots,\phi_v,
Y,\psi_1,\dots,\psi_\tau)=
\\&\sum_{t\in V^{k}} a_t X_{+}div_{i_1}\dots X_{+}div_{i_{k}}
C^{t,i_1\dots i_k}_g(\Omega_1,\dots,\Omega_b,\phi_1,\dots,\phi_v,
Y,\psi_1,\dots,\psi_\tau)+
\\&\sum_{\zeta\in Z_{b,No@}} a_\zeta X_{+}div_{i_1} C^{\zeta,i_1}_g
(\Omega_1,\dots,\Omega_b,\phi_1,\dots,\phi_v,
Y,\psi_1,\dots,\psi_\tau)+
\\&\sum_{q\in Q} a_q X_{+}div_{i_1}
C^{q,i_1}_g(\Omega_1,\dots,\Omega_b,\phi_1,\dots,\phi_v,
Y,\psi_1,\dots,\psi_\tau)+
\\&\sum_{j\in J} a_j C^{j}_g(\Omega_1,\dots,\Omega_b,\phi_1,\dots,\phi_v,
Y,\psi_1,\dots,\psi_\tau),
\end{split}
\end{equation}

\begin{equation}
\label{proolmertb2'}
\begin{split}
&\sum_{\zeta\in Z_{b,@}} a_\zeta X_{+}div_{i_1}
C^{\zeta,i_1}_g(\Omega_1,\dots,\Omega_b,\phi_1,\dots,\phi_v,
[\chi_1,\chi_2],\psi_1,\dots,\psi_\tau)=
\\&\sum_{t\in V^{k}} a_t X_{+}div_{i_1}\dots X_{+}div_{i_{k}}
C^{t,i_1\dots i_k}_g(\Omega_1,\dots,\Omega_b,\phi_1,\dots,\phi_v,
[\chi_1,\chi_2],\psi_1,\dots,\psi_\tau)+
\\&\sum_{\zeta\in Z_{b,No@}} a_\zeta X_{+}div_{i_1} C^{\zeta,i_1}_g
(\Omega_1,\dots,\Omega_b,\phi_1,\dots,\phi_v,
[\chi_1,\chi_2],\psi_1,\dots,\psi_\tau)+
\\&\sum_{q\in Q} a_q X_{+}div_{i_1}
C^{q,i_1}_g(\Omega_1,\dots,\Omega_b,\phi_1,\dots,\phi_v,
[\chi_1,\chi_2],\psi_1,\dots,\psi_\tau)+
\\&\sum_{j\in J} a_j C^{j}_g(\Omega_1,\dots,\Omega_b,\phi_1,\dots,\phi_v,
[\chi_1,\chi_2],\psi_1,\dots,\psi_\tau),\end{split}
\end{equation}
where the tensor fields indexed in $V^k$ have all the features of
the tensor fields indexed in $Z_{b,@}$ but in addition have all
the $k$ free indices {\it not} belonging to factors
$\nabla\psi_1,\dots,\nabla\psi_\tau$. The tensor fields indexed in
$Z_{b,No@}$ have all the regular features of the tensor fields in
$Z_b$ and in addition have the factor $\nabla Y$ (or the factors
$\nabla\chi_1,\nabla\chi_2$) {\it not} contracting against special
indices. The terms indexed in $Q$ are as required in the RHS of
(\ref{corbatis1}), (\ref{corbatis2}) (which are the equations that
we are proving).

\par We will then show that we can write:

\begin{equation}
\label{proolmert1c'}
\begin{split}
&\sum_{\zeta\in Z_{b,@}} a_\zeta X_{+}div_{i_1}
C^{\zeta,i_1}_g(\Omega_1,\dots,\Omega_b,\phi_1,\dots,\phi_v,
Y,\psi_1,\dots,\psi_\tau)=
\\&\sum_{t\in V^{k+1}} a_t X_{+}div_{i_1}\dots X_{+}div_{i_{k+1}}
C^{t,i_1\dots
i_{k+1}}_g(\Omega_1,\dots,\Omega_b,\phi_1,\dots,\phi_v,
Y,\psi_1,\dots,\psi_\tau)+
\\&\sum_{\zeta\in Z_{b,No@}} a_\zeta X_{+}div_{i_1} C^{\zeta,i_1}_g
(\Omega_1,\dots,\Omega_b,\phi_1,\dots,\phi_v,
Y,\psi_1,\dots,\psi_\tau)+
\\&\sum_{q\in Q} a_q X_{+}div_{i_1}
C^{q,i_1}_g(\Omega_1,\dots,\Omega_b,\phi_1,\dots,\phi_v,
Y,\psi_1,\dots,\psi_\tau)+
\\&\sum_{j\in J} a_j C^{j}_g(\Omega_1,\dots,\Omega_b,\phi_1,\dots,\phi_v,
Y,\psi_1,\dots,\psi_\tau),
\end{split}
\end{equation}

\begin{equation}
\label{proolmertc2'}
\begin{split}
&\sum_{\zeta\in Z_{b,@}} a_\zeta X_{+}div_{i_1}
C^{\zeta,i_1}_g(\Omega_1,\dots,\Omega_b,\phi_1,\dots,\phi_v,
[\chi_1,\chi_2],\psi_1,\dots,\psi_\tau)=
\\&\sum_{t\in V^{k+1}} a_t X_{+}div_{i_1}\dots X_{+}div_{i_{k+1}}
C^{t,i_1\dots
i_{k+1}}_g(\Omega_1,\dots,\Omega_b,\phi_1,\dots,\phi_v,
[\chi_1,\chi_2],\psi_1,\dots,\psi_\tau)
\\&+\sum_{\zeta\in Z_{b,No@}} a_\zeta X_{+}div_{i_1} C^{\zeta,i_1}_g
(\Omega_1,\dots,\Omega_b,\phi_1,\dots,\phi_v,
[\chi_1,\chi_2],\psi_1,\dots,\psi_\tau)+
\\&\sum_{q\in Q} a_q X_{+}div_{i_1}
C^{q,i_1}_g(\Omega_1,\dots,\Omega_b,\phi_1,\dots,\phi_v,
[\chi_1,\chi_2],\psi_1,\dots,\psi_\tau)+
\\&\sum_{j\in J} a_j C^{j}_g(\Omega_1,\dots,\Omega_b,\phi_1,\dots,\phi_v,
[\chi_1,\chi_2],\psi_1,\dots,\psi_\tau).\end{split}
\end{equation}
(Here the tensor fields indexed in $V^{k+1}$ have all 
the features described above and moreover have rank $k+1$). 

\par Thus, by iterative repetition of this step we are reduced to
showing our claim under the additional assumption that
$Z'_{a,@}=Z_{b,@}=\emptyset$.
\newline

We prove (\ref{proolmert1c'}), (\ref{proolmertc2'}) below. Now, we 
 present the rest of our claims under the
assumption that $Z'_{a,@}=Z_{b,@}=\emptyset$. For the rest of this
proof we will be assuming that all tensor fields have the factor
$\nabla Y$ (or the factors $\nabla\chi_1,\nabla\chi_2$) not
contracting against special indices.

\par We then perform a new induction: We assume that we can write:

\begin{equation}
\label{olmert1}
\begin{split}
&\sum_{\zeta\in Z'_a} a_\zeta X_{+}div_{i_1}
C^{\zeta,i_1}_g(\Omega_1,\dots,\Omega_b,\phi_1,\dots,\phi_v,
Y,\psi_1,\dots,\psi_\tau)=
\\&\sum_{\zeta\in Z_b} a_\zeta X_{+}div_{i_1}
C^{\zeta,i_1}_g(\Omega_1,\dots,\Omega_b,\phi_1,\dots,\phi_v,
Y,\psi_1,\dots,\psi_\tau)+
\\&\sum_{t\in T^k} a_t X_{+}div_{i_1}\dots X_{+}div_{i_k}
C^{t,i_1\dots i_k}_g(\Omega_1,\dots,\Omega_b,\phi_1,\dots,\phi_v,
Y,\psi_1,\dots,\psi_\tau)+
\\&\sum_{q\in Q} a_q X_{+}div_{i_1}
C^{q,i_1}_g(\Omega_1,\dots,\Omega_b,\phi_1,\dots,\phi_v,
Y,\psi_1,\dots,\psi_\tau)+
\\&\sum_{j\in J} a_j C^{j}_g(\Omega_1,\dots,\Omega_b,\phi_1,\dots,\phi_v,
Y,\psi_1,\dots,\psi_\tau),
\end{split}
\end{equation}

\begin{equation}
\label{olmert2}
\begin{split}
&\sum_{\zeta\in Z'_a} a_\zeta X_{+}div_{i_1}
C^{\zeta,i_1}_g(\Omega_1,\dots,\Omega_b,\phi_1,\dots,\phi_v,
[\chi_1,\chi_2],\psi_1,\dots,\psi_\tau)=
\\&\sum_{\zeta\in Z_b} a_\zeta X_{+}div_{i_1}
C^{\zeta,i_1}_g(\Omega_1,\dots,\Omega_b,\phi_1,\dots,\phi_v,
[\chi_1,\chi_2],\psi_1,\dots,\psi_\tau)+
\\&\sum_{t\in T^k} a_t X_{+}div_{i_1}\dots X_{+}div_{i_k}
C^{t,i_1\dots i_k}_g(\Omega_1,\dots,\Omega_b,\phi_1,\dots,\phi_v,
[\chi_1,\chi_2],\psi_1,\dots,\psi_\tau)+
\\&\sum_{q\in Q} a_q X_{+}div_{i_1}
C^{q,i_1}_g(\Omega_1,\dots,\Omega_b,\phi_1,\dots,\phi_v,
[\chi_1,\chi_2],\psi_1,\dots,\psi_\tau)+
\\&\sum_{j\in J} a_j C^{j}_g(\Omega_1,\dots,\Omega_b,\phi_1,\dots,\phi_v,
[\chi_1,\chi_2],\psi_1,\dots,\psi_\tau),
\end{split}
\end{equation}
where the tensor fields indexed in $T^k$ have all the properties of the
tensor fields indexed in $Z'_a$ (in particular the index in
$\nabla\psi_1$ is not free) and in addition have rank $k$. We then
show that we can write:

\begin{equation}
\label{olmert1'}
\begin{split}
&\sum_{\zeta\in Z'_a} a_\zeta X_{+}div_{i_1}
C^{\zeta,i_1}_g(\Omega_1,\dots,\Omega_b,\phi_1,\dots,\phi_v,
Y,\psi_1,\dots,\psi_\tau)=
\\&\sum_{\zeta\in Z_b} a_\zeta X_{+}div_{i_1}
C^{\zeta,i_1}_g(\Omega_1,\dots,\Omega_b,\phi_1,\dots,\phi_v,
Y,\psi_1,\dots,\psi_\tau)+
\\&\sum_{t\in T^{k+1}} a_t X_{+}div_{i_1}\dots X_{+}div_{i_{k+1}}
C^{t,i_1\dots i_k}_g(\Omega_1,\dots,\Omega_b,\phi_1,\dots,\phi_v,
Y,\psi_1,\dots,\psi_\tau)+
\\&\sum_{q\in Q} a_q X_{+}div_{i_1}
C^{q,i_1}_g(\Omega_1,\dots,\Omega_b,\phi_1,\dots,\phi_v,
Y,\psi_1,\dots,\psi_\tau)+
\\&\sum_{j\in J} a_j C^{j}_g(\Omega_1,\dots,\Omega_b,\phi_1,\dots,\phi_v,
Y,\psi_1,\dots,\psi_\tau),
\end{split}
\end{equation}

\begin{equation}
\label{olmert2'}
\begin{split}
&\sum_{\zeta\in Z'_a} a_\zeta X_{+}div_{i_1}
C^{\zeta,i_1}_g(\Omega_1,\dots,\Omega_b,\phi_1,\dots,\phi_v,
[\chi_1,\chi_2],\psi_1,\dots,\psi_\tau)=
\\&\sum_{\zeta\in Z_b} a_\zeta X_{+}div_{i_1}
C^{\zeta,i_1}_g(\Omega_1,\dots,\Omega_b,\phi_1,\dots,\phi_v,
[\chi_1,\chi_2],\psi_1,\dots,\psi_\tau)+
\\&\sum_{t\in T^{k+1}} a_t X_{+}div_{i_1}\dots X_{+}div_{i_{k+1}}
C^{t,i_1\dots
i_{k+1}}_g(\Omega_1,\dots,\Omega_b,\phi_1,\dots,\phi_v,
[\chi_1,\chi_2],\psi_1,\dots,\psi_\tau)
\\&+\sum_{q\in Q} a_q X_{+}div_{i_1}
C^{q,i_1}_g(\Omega_1,\dots,\Omega_b,\phi_1,\dots,\phi_v,
[\chi_1,\chi_2],\psi_1,\dots,\psi_\tau)+
\\&\sum_{j\in J} a_j C^{j}_g(\Omega_1,\dots,\Omega_b,\phi_1,\dots,\phi_v,
[\chi_1,\chi_2],\psi_1,\dots,\psi_\tau).
\end{split}
\end{equation}

\par We will derive (\ref{olmert1'}), (\ref{olmert2'}) momentarily.
For now, we observe that by iterative repetition of the above we
are reduced to showing (\ref{corbatis1}), (\ref{corbatis2}) under
the additional assumption that $Z'_a=\emptyset$. In that setting,
we can just repeatedly apply the eraser (see the Appendix in \cite{alexakis4} 
for a definition of this notion) to as many factors $\nabla\psi_\tau$
as needed in order to reduce ourselves to a new true equation where
each of the real factors is contracting against at most one of the factors
 $\nabla\psi_1,\dots,\nabla\psi_\tau,\nabla Y$
(or $\nabla\chi_1,\nabla\chi_2$).\footnote{All
remaining factors $\nabla\psi_1,\dots,\nabla\psi_\tau$ and also
 the factor(s) $\nabla Y$ (or $\nabla\chi_1,\nabla\chi_2$) are treated as factors $\nabla\phi_h$}
 Then, by invoking Corollary 1 from \cite{alexakis4}\footnote{Notice that there will necessarily be
at least one non-simple factor $S_{*}\nabla^{(\nu)}R_{ijkl}$ or
$\nabla^{(B)}\Omega_h$, by virtue of the factor(s) $\nabla Y$ (or
$\nabla\omega_1,\nabla\omega_2$), therefore that Corollary
 can be applied.} and then re-introducing the
factors we erased, we derive our claim.
\newline

{\it Proof of (\ref{olmert1'}), (\ref{olmert2'}):} Picking out the
sublinear combination in (\ref{olmert1}), (\ref{olmert2}) with one
derivative on $\nabla Y$ or $\nabla\chi_1$ and substituting into
(\ref{batis1}), (\ref{batis2}) we derive a new equation:

\begin{equation}
\label{expected1}
\begin{split}
&\sum_{t\in T^k} a_t X_{*}div_{i_1}\dots X_{*}div_{i_k}
C^{t,i_1\dots i_k}_g(\Omega_1,\dots,\Omega_b,\phi_1,\dots,\phi_v,
Y,\psi_1,\dots,\psi_\tau)+
\\&\sum_{\zeta\in Z_b} a_\zeta X_{*}div_{i_1}
C^{\zeta,i_1}_g(\Omega_1,\dots,\Omega_b,\phi_1,\dots,\phi_v,
Y,\psi_1,\dots,\psi_\tau)=
\\&\sum_{j\in J} a_j C^{j}_g(\Omega_1,\dots,\Omega_b,\phi_1,\dots,\phi_v,
Y,\psi_1,\dots,\psi_\tau),
\end{split}
\end{equation}

\begin{equation}
\label{expected2}
\begin{split}
&\sum_{t\in T^k} a_t X_{*}div_{i_1}\dots X_{*}div_{i_k}
C^{t,i_1\dots i_k}_g(\Omega_1,\dots,\Omega_b,\phi_1,\dots,\phi_v,
[\chi_1,\chi_2],\psi_1,\dots,\psi_\tau)+
\\&\sum_{\zeta\in Z_b} a_\zeta X_{*}div_{i_1}
C^{\zeta,i_1}_g(\Omega_1,\dots,\Omega_b,\phi_1,\dots,\phi_v,
[\chi_1,\chi_2],\psi_1,\dots,\psi_\tau)=
\\&\sum_{j\in J} a_j C^{j}_g(\Omega_1,\dots,\Omega_b,\phi_1,\dots,\phi_v,
[\chi_1,\chi_2],\psi_1,\dots,\psi_\tau);
\end{split}
\end{equation}
(the sublinear combination $\sum_{\zeta\in Z_b}\dots$ above is {\it generic}).

\par We now divide the index set $T^k$ according to which of the
factors $\nabla\psi_2$,$\dots$,$\nabla\psi_\tau$,$\nabla Y$ (or
$\nabla\psi_1,\dots,\nabla\psi_\tau,\nabla \chi_1$) contain the
$k$ free indices. Thus we write: $T^k=\bigcup_{\alpha\in A}
T^{k,\alpha}$ (each $\alpha\in A$ corresponds to a $k$-subset of
the set of factors $\nabla\psi_1$,$\dots,\nabla\psi_\tau$,$\nabla Y$
or $\nabla\psi_1$,$\dots$,$\nabla\psi_\tau$,$\nabla \chi_1$). We will
then show that for each $\alpha\in A$ there exists a tensor field
$\sum_{b\in B^\alpha} a_b C^{b,i_1\dots i_{k+1}}_g$ in the form
(\ref{proherodium}) or (\ref{herodium}) with the first $k$ free
indices belonging to the factors in the set $\alpha$, and the free
index ${}_{i_{k+1}}$ {\it not} belonging to $\nabla\psi_1$, so
that:

\begin{equation}
\label{goldfarm1}
\begin{split}
&\sum_{t\in T^{k,\alpha}} a_t C^{t,i_1\dots
i_k}_g(\Omega_1,\dots,\Omega_b,\phi_1,\dots,\phi_v,
Y,\psi_1,\dots,\psi_\tau)\nabla_{i_1}\upsilon\dots\nabla_{i_k}\upsilon-
\\&X_{*}div_{i_{k+1}}\sum_{b\in B^\alpha} a_b
C^{b,i_1\dots
i_{k+1}}_g(\Omega_1,\dots,\Omega_b,\phi_1,\dots,\phi_v,
Y,\psi_1,\dots,\psi_\tau)\nabla_{i_1}\upsilon\dots\nabla_{i_k}\upsilon
\\&=\sum_{j\in J} a_jC^{j,i_1\dots
i_k}_g(\Omega_1,\dots,\Omega_b,\phi_1,\dots,\phi_v,
Y,\psi_1,\dots,\psi_\tau)\nabla_{i_1}\upsilon\dots\nabla_{i_k}\upsilon,
\end{split}
\end{equation}

\begin{equation}
\label{goldfarm2}
\begin{split}
&\sum_{t\in T^{k,\alpha}} a_t C^{t,i_1\dots
i_k}_g(\Omega_1,\dots,\Omega_b,\phi_1,\dots,\phi_v,
[\chi_1,\chi_2],\psi_1,\dots,\psi_\tau)\nabla_{i_1}\upsilon\dots\nabla_{i_k}\upsilon-
\\&X_{*}div_{i_{k+1}}\sum_{b\in B^\alpha} a_b
C^{b,i_1\dots
i_{k+1}}_g(\Omega_1,\dots,\Omega_b,\phi_1,\dots,\phi_v,
[\chi_1,\chi_2],\psi_1,\dots,\psi_\tau)\nabla_{i_1}\upsilon\dots\nabla_{i_k}\upsilon
\\&=\sum_{j\in J} a_jC^{j,i_1\dots
i_k}_g(\Omega_1,\dots,\Omega_b,\phi_1,\dots,\phi_v,
[\chi_1,\chi_2],\psi_1,\dots,\psi_\tau)\nabla_{i_1}\upsilon\dots\nabla_{i_k}\upsilon.
\end{split}
\end{equation}
If we can show the above for every $\alpha\in A$, then replacing
the factor $\nabla\upsilon$ by $X_{+}div$'s we can derive our
claim (\ref{olmert1'}), (\ref{olmert2'}).

{\it Proof of (\ref{goldfarm1}), (\ref{goldfarm2}):} Refer to
(\ref{expected1}) and (\ref{expected2}). Denote $Y$ or $\chi_1$ by
$\psi_{\tau+1}$ for uniformity. We pick out any $\alpha\in A$;
assume that $\alpha=\{\nabla\psi_{x_1},\dots ,\nabla\psi_{x_k}\}$.

Pick out the sublinear combination where the factors
$\nabla\psi_{x_1},\dots ,\nabla\psi_{x_k}$ which belong to
$\alpha$ are contracting against the same factor as
$\nabla\psi_1$. This sublinear combination $Z_g$ vanishes
separately (i.e. $Z_g=0$). We then apply the eraser to the factors
$\nabla\psi_2,\dots,\nabla Y\in A$ (notice this is {\it
well-defined}, since all the above factors {\it and} the factor
$\nabla\psi_1$ are contracting against non-special indices). We
obtain a new true equation, which we denote by $Erase[Z_g]=0$. It
then follows that $Erase[Z_g]\cdot(\nabla_{i_1}\psi_{x_1}
\nabla^{i_1}\upsilon\dots\nabla_{i_k}\psi_{x_k}\nabla^{i_k}\upsilon)=0$
is our desired conclusion (\ref{goldfarm1}), (\ref{goldfarm2}).
$\Box$
\newline

{\it (Sketch of) Proof of (\ref{proolmert1'}), (\ref{proolmert2'})
(\ref{proolmert1c'}), (\ref{proolmertc2'}):} These equations can
be proven by only a slight modification of the idea above. We
again subdivide the index sets $T^k,V^k$ according to the set of
factors $\nabla\psi_2,\dots,\nabla\psi_\tau$ or
$\nabla\psi_2,\dots,\nabla\psi_\tau,\nabla\omega_1$ which contain
the $k$ free indices (so we write $T^{k}=\bigcup_{\alpha\in A}
T^{k,\alpha}$ and $V^k=\bigcup_{\alpha\in A} V^{k,\alpha}$)
 and we prove the claims above separately for those sublinear
combinations.

\par To prove this, we pick out the sublinear combination in our hypotheses with
the factors $\nabla\psi_h$, $h\in \alpha$ contracting against the
same factor against which $\nabla\psi_1$ and $\nabla Y$ (or
$\nabla\psi_1$ and $\nabla\omega_1$) are contracting. Say
$\alpha=\{h_1,\dots,h_k\}$; we then formally replace the
expressions $S_{*}\nabla^{(\nu)}_{r_1\dots r_\mu l_1\dots
l_k}R_{ijkl}\nabla_{l_1}\psi_{h_1}\dots
\nabla^{l_k}\psi_{h_k}\nabla^i\tilde{\phi}_1\nabla^j\psi_1\nabla^kY$
or $\nabla^{(A)}_{r_1\dots r_\mu l_1\dots
l_kst}\Omega_1\nabla^{l_1}\psi_{h_1}\dots
\nabla^{l_k}\psi_{h_k}\nabla^s\psi_1\nabla^tY$ etc,  by
expressions \\$S_{*}\nabla^{(\nu-k)}_{r_1\dots
r_\mu}R_{ijkl}\nabla^i\tilde{\phi}_1\nabla^j\psi_1\nabla^kY$,
$\nabla^{(A-k)}_{r_1\dots r_\mu
st}\Omega_1\nabla^s\psi_1\nabla^tY$ and derive our
claims (\ref{proolmert1'}), (\ref{proolmert2'})
(\ref{proolmert1c'}), (\ref{proolmertc2'}) as above. $\Box$
\newline

{\it Proof of the claims of our general statement (i.e.
(\ref{procareless2}), (\ref{careless2}) by induction):} We will
prove these claims by an induction. Our inductive assumptions are
that (\ref{procareless2}), (\ref{careless2}) follow from
(\ref{procareless}), (\ref{careless}) for any weight $-W'$, $W'<
K$ and when $W'=K$ they hold for any
 length $v+b\ge \gamma\ge 2$. We will then show the claim when the
 weight is $-K$, and $v+b=\gamma+1$. In the end,
 we will check our claims for the base case $v+b=2$.

{\it Proof of the inductive step:} Refer back to
(\ref{procareless}), (\ref{careless}). We will prove this claim in
four steps.
\newline

{\it Step 1:} Firstly, we will denote by
$Z_a^{spec},\overline{Z}^{spec}_a, Z_b^{spec}$ the index sets of
the tensor fields for which $\nabla Y$ or one of the factors
$\nabla\chi_1$, $\nabla\chi_2$
 (respectively) is contracting against
 a special index.
 Then {\it using the inductive assumptions
 of our general claim},
 we will show that there exists a linear
 combination of 2-tensor fields (indexed in $W$ below)
 which satisfies all the requirements of
 (\ref{procareless}), (\ref{procareless2}) so that:

\begin{equation}
\label{yirgachef}
\begin{split}
&\sum_{\zeta\in Z_a^{spec}} a_\zeta
C^{\zeta,i_1}_g\nabla_{i_1}\upsilon-X_{*}div_{i_2}\sum_{w\in W}
a_w C^{w,i_1i_2}_g\nabla_{i_1}\upsilon=
\\&\sum_{\zeta\in Z_a^{OK}} a_\zeta
C^{\zeta,i_1}_g\nabla_{i_1}\upsilon+ \sum_{j\in J} a_j
C^{j,i_1}_g\nabla_{i_1}\upsilon,
\end{split}
\end{equation}
where the tensor fields n $Z_a^{OK}$ are {\it generic} linear
combinations of tensor fields of the same general type as the ones
indexed in $Z_a$ in (\ref{procareless}), (\ref{procareless2}) and
where in addition none of the factors $\nabla Y$ or
$\nabla\chi_1,\nabla\chi_2$ are contracting against a special
index.

\par Thus, if we can show the above, by replacing $\nabla\upsilon$
by an $X_{*}div_i$, and substituting back into
(\ref{procareless}), (\ref{procareless2}), we are reduced to
showing (\ref{careless}), (\ref{careless2}) under the additional
assumption that $Z_a^{spec}=\emptyset$.
\newline

{\it Step 2:} Then, under the assumption that
$Z_a^{spec}=\emptyset$, we will show that we can write:

\begin{equation}
\label{yirgachef2}
\begin{split}
&\sum_{\zeta\in Z_b^{spec}} a_\zeta X_{*}div_{i_1}C^{\zeta,i_1}_g+
\sum_{\zeta\in \overline{Z}_a^{spec}} a_\zeta
X_{*}div_{i_1}\dots X_{*}div_{i_c}C^{\zeta,i_1\dots i_c}_g=
\\&X_{*}div_{i_1}\dots X_{*}div_{i_b}\sum_{c\in C} a_c
C^{c,i_1\dots i_b}_g + \sum_{j\in J} a_j C^{j,i_1}_g,
\end{split}
\end{equation}
where the tensor fields on the RHS are of the general form as the
 ones indexed in $Z_b,\overline{Z}_a$ in our hypothesis,
 and moreover the factor $\nabla Y$ (or the
  factors $\nabla\chi_1,\nabla\chi_2$) is (are) not
contracting against special indices.

\par Notice that if we can show (\ref{yirgachef}), (\ref{yirgachef2})
then we are reduced to showing our claim under the additional
assumption that for each $\zeta\in Z_a\bigcup \overline{Z}_a\bigcup Z_b$ the factor(s)
$\nabla Y$ (or $\nabla\chi_1,\nabla\chi_2$) are not contracting
against special indices. We will show (\ref{yirgachef}),
(\ref{yirgachef2}) below.
\newline

{\it Proof of (\ref{careless}), (\ref{careless2}) under the
additional assumption that for each
 $\zeta\in Z_a\bigcup \overline{Z}_a\bigcup Z_b$ the factor
$\nabla Y$ or ($\nabla\chi_1,\nabla\chi_2$) is not contracting
against special indices:}
\newline

{\it Step 3: Proof of (\ref{princecharles}) below:}
\newline

 We note that for all
 the tensor fields in the rest of this proof
will {\it not} have the factor $\nabla Y$ (or any of the factors
$\nabla\chi_1,\nabla\chi_2$) contracting against a special index
in any factor $S_{*}\nabla^{(\nu)}R_{ijkl}$ or
$\nabla^{(B)}\Omega_h$.
 Now,  we arbitrarily pick out one factor
$T=S_{*}\nabla^{(\nu)}R_{ijkl}$ or $T=\nabla^{(B)}\Omega_x$
 in $\overline{\kappa}_{simp}$ and call it the
``chosen factor'' for the rest of this subsection.

\par  We  will say that the factor $\nabla Y$ (or $\nabla\omega_2$)
is contracting against a good index in $T$, if it is contracting
against
 a non-special index in  $T$
when $T$ is of the form $S_{*}\nabla^{(\nu)}R_{ijkl}$ with
$\nu>0$;  when $T$ is of the form $\nabla^{(B)}\Omega_x$,
 then it is contracting against a good index provided $B\ge 3$.

\par We will say that the factor $\nabla Y$ (or $\nabla\omega_2$)
 is contracting against a bad index if it is
contracting against the index ${}_j$ in a factor $T=S_{*}R_{ijkl}$
or an index in a factor $T=\nabla^{(2)}\Omega_x$.  We denote by
$Z_a^{BAD}\subset Z_a$ the index set of tensor fields
 for which $\nabla Y$ (or $\nabla\omega_2$) is
contracting against a bad index. We also denote by
$Z_b^{BAD}\subset Z_b$ the index set of the vector fields for
which $\nabla Y$ is contracting against a bad index in $T$ and $T$
also contains a free index. We will show that we can write:

\begin{equation}
\label{princecharles}
\begin{split}
&\sum_{\zeta\in Z_a^{BAD}\bigcup Z_b^{BAD}} a_\zeta
 C^{\zeta,i_1}_g\nabla_{i_1}\upsilon-
X_{*}div_{i_2}\sum_{h\in H} a_h C^{i_1i_2}_g\nabla_{i_1}\upsilon=
\\&\sum_{\zeta\in Z_a'^{GOOD}\bigcup Z_b'^{GOOD}} a_\zeta
 C^{\zeta,i_1}_g\nabla_{i_1}\upsilon+
\sum_{j\in J} a_j C^j_g,
\end{split}
\end{equation}
where all the tensor fields indexed in
 $Z_a'^{GOOD}\bigcup Z_b'^{GOOD}$ are
 {\it generic} vector fields of the forms indexed in
$Z_a,Z_b$, only with the factors $\nabla Y$ or $\nabla \omega_2$
contracting against a {\it good} index in the factor $T$. The
tensor fields indexed in $H$ are as required in the claim of our
general statement (they correspond to the index set $W$ in our
general statement).
\newline

{\it Step 4: Proof that (\ref{princecharles}) implies our claims
 (\ref{procareless2}), (\ref{careless2}).}
\newline

\par We start by proving (\ref{princecharles}) (i.e.~we prove Step 3).
Then, we will show how we can derive our claim from
(\ref{princecharles}) (i.e.~we then prove Step 4).
\newline

{\it Proof of Step 3: Proof of (\ref{princecharles}):}
 We can prove this equation by virtue of our inductive assumption on our general claim.
First, we define $\overline{Z}_a^{BAD}\subset \overline{Z}_a$ to
stand for the index set of tensor fields where the factor $\nabla
Y$ (or $\nabla\omega_2$) is contracting against a bad index in the
chosen factor.   We pick out the sublinear combination in our
Lemma assumption where $\nabla Y$ (or $\nabla \omega_2$) are
contracting against the chosen factor $T=S_{*}R_{ijkl}$ or
$T=\nabla^{(2)}\Omega_x$). This sublinear combination must vanish
separately, and we thus derive an equation:

\begin{equation}
\label{princecharles'}
\begin{split}
&\sum_{\zeta\in Z_a^{BAD}\bigcup Z_b^{BAD}} a_\zeta
 X_{**}div_{i_1} C^{\zeta,i_1}_g+\sum_{\zeta\in \overline{Z}_a^{BAD}}
  a_\zeta X_{**}div_{i_1}\dots X_{**}div_{i_c} C^{\zeta,i_1\dots i_c}_g+
  \\& \sum_{\zeta\in Z^{nvBAD}_b} a_f C^{f,i_1}_g=
\sum_{j\in J} a_j C^j_g,
\end{split}
\end{equation}
where $X_{**}div_{i_1}$  stands for the sublinear combination
for which $\nabla_{i_1}$ is not allowed to hit the chosen factor $T$.
$Z^{nvBAD}_b\subset Z_b$ stands for the index set of tensor fields
 indexed in $Z_b$ with the free index ${}_{i_1}$
{\it not} belonging to the chosen factor and  also with the factor
$\nabla  Y$ (or $\nabla\omega_2$) contracting against a bad index.

\par Now, define an operation $Op[\dots]$ which acts on the complete
 contractions above by formally replacing any expression
$\nabla^{(2)}_{ij}\Omega_x\nabla^iY$ (or
$\nabla^{(2)}_{ij}\Omega_x\nabla^i\chi_2$) by $\nabla_jD$ ($D$ is
a scalar function), or any expression
$S_{*}R_{ijkl}\nabla^i\tilde{\phi}_1\nabla^j Y$ (or
\\$S_{*}R_{ijkl}\nabla^i\tilde{\phi}_1\nabla^j \chi_2$)
 by $\nabla_{[k}\theta_1\nabla_{l]}\theta_2$. (Denote by
 $\tilde{\kappa}_{simp}$ the simple character of these resulting vector
 fields).
Acting on (\ref{princecharles'}) by $Op[\dots]$ produces a true
equation, which we may write out as:

\begin{equation}
\label{princecharles''}
\begin{split}
&\sum_{\zeta\in Z_a^{BAD}\bigcup Z_b^{BAD}} a_\zeta
 X_{**}div_{i_1} Op[C]^{\zeta,i_1}_g+X_{**}div_{i_1}
 \sum_{f\in F} a_f C^{f,i_1}_g
\\&+\sum_{\zeta\in \overline{Z}_a^{BAD}}
  a_\zeta X_{**}div_{i_1}\dots X_{**}div_{i_c} C^{\zeta,i_1\dots i_c}_g=
\sum_{j\in J} a_j C^j_g.
\end{split}
\end{equation}
Here $X_{**}div_i$ stands for the sublinear combination in $div_i$
where $\nabla_i$ is not allowed to hit the factor
 to which $\nabla_i$ belongs, nor any of the factors $\nabla\phi_1,\dots,\nabla\phi_u$,
$\nabla\psi_1,\dots,\nabla\psi_\tau$, nor any factors $\nabla
D,\nabla\theta_1,\nabla\theta_2$. The vector fields indexed in $F$
are generic vector fields with a simple character
$\tilde{\kappa}_{simp}$, for which the free index ${}_{i_1}$ {\it
does not} belong to any of the factors
$\nabla\psi_1,\dots,\nabla\psi_\tau$ or any of the factors $\nabla
D,(\nabla\chi_1), \nabla\theta_1,\nabla\theta_2$.

\par Now, {\it observe that the above equation
falls under our inductive assumption of the general statement we
are proving}: We now either have factors
$\nabla\psi_1$,$\dots$,$\nabla\psi_\tau$,$\nabla D$, or
$\nabla\psi_1$,$\dots$,$\nabla\psi_\tau$,$\nabla\chi_1$,$\nabla D$ or
$\nabla\psi_1,$\dots,$\nabla\psi_\tau$,$[\nabla
\theta_1,\nabla\theta_2]$ or
$\nabla\psi_1$,$\dots$,$\nabla\psi_\tau$,$\nabla\chi_1$,$[\nabla
\theta_1,\nabla\theta_2]$. Notice that the tensor fields indexed
in $H_a^{BAD},H_b^{BAD}$ are precisely the ones that contain a
free index in one of these factors. Therefore,
by our inductive assumption of the ``general claim'' we derive that
there exists a linear
 combination of 2-tensor fields, $\sum_{v\in V}\dots$, (with factors
$\nabla\psi_1$,$\dots$,$\nabla\psi_\tau$,$\nabla D$ etc, and which
satisfy the $@$-property for the
 factors $\nabla\psi_1$,$\dots$,$\nabla\psi_\tau$) so that:

\begin{equation}
\label{princecharles'''}
\begin{split}
&\sum_{\zeta\in Z_a^{BAD}\bigcup Z_b^{BAD}} a_\zeta
  Op[C]^{\zeta,i_1}_g\nabla_{i_1}\upsilon-
 X_{**}div_{i_2}\sum_{v\in V} a_v C^{v,i_1i_2}_g\nabla_{i_1}\upsilon=
\\&\sum_{j\in J} a_j C^{j,i_1}_g\nabla_{i_1}\upsilon.
\end{split}
\end{equation}
Now, we define an operation $Op^{-1}[\dots]$, which acts on the
complete
 contractions in the above equation by replacing the factor $\nabla_j D$ by an expression
$\nabla_{ij}\Omega_x\nabla^jY$ (or
$\nabla_{ij}\Omega_x\nabla^j\omega_2$) or the expression
$\nabla_{[a}\theta_1\nabla_{b]}\theta_2$ by
$S_{*}R_{ijab}\nabla^i\tilde{\phi}_1\nabla^jY$ (or
$S_{*}R_{ijab}\nabla^i\tilde{\phi}_1\nabla^j\omega_2$). The
operation $Op^{-1}$ clearly produces a true equation, which
 is our desired conclusion, (\ref{princecharles}). $\Box$
\newline

{\it Proof of Step 4:} We derive our conclusions
 (\ref{procareless2}), (\ref{careless2}) in pieces.
Firstly, we show these equations with the sublinear combinations
$Z_a$ replaced by the index set $Z_{a,spec}$, which index the
terms with the free index ${}_{i_1}$ belonging to the factor
$\nabla Y$ or $\nabla\omega_1$ (this will be sub-step A). After
proving this claim, we will  show
(\ref{procareless2}), (\ref{careless2}) under the additional
assumption that $Z_{a,spec}=\emptyset$ (this will be sub-step B).

{\it Proof of sub-step A:} We make the $\nabla\upsilon$'s into
$X_{*}div$'s in (\ref{princecharles}) and replace the resulting equations into our
Lemma hypothesis. We thus derive a new equation:

\begin{equation}
\label{tolook}
\begin{split}
&\sum_{\zeta\in Z_a} a_\zeta X_{*}div_{i_1} C^{\zeta,i_1}_g +
\sum_{\zeta\in Z_b^1} a_\zeta  X_{*}div_{i_1} C^{\zeta,i_1}_g+
\\& \sum_{\zeta\in Z_b^2} a_\zeta
 X_{*}div_{i_1}\dots X_{*}div_{i_a} C^{\zeta,i_1\dots i_a}_g
+\sum_{j\in J} a_j C^j_g=0,
\end{split}
\end{equation}
where we now have the tensor fields indexed in $Z_a$ have a free
index
 among the factors $\nabla\psi_1,\dots,\nabla\psi_\tau,\nabla Y$
(or $\nabla\psi_1,\dots,\nabla\psi_\tau,\nabla
\chi_1,\nabla\chi_2$), and furthermore the factor $\nabla Y$ (or
the factors
 $\nabla\omega_1,\nabla\omega_2$) are not contracting against
 a bad index in the chosen factor $T$.
The tensor fields indexed in $Z_b^1$ have a free index that {\it
does not}
 belong to one of the factors $\nabla\psi_1,\dots,\nabla\psi_\tau,\nabla Y$
(or $\nabla\psi_1,\dots,\nabla\psi_\tau,\nabla
\chi_1,\nabla\chi_2$), and furthermore {\it if} the factor $\nabla
Y$ (or one of
 the factors $\nabla\omega_1,\nabla\omega_2$) is contracting
 against a bad index in the chosen factor $T$, then $T$ does
 not contain the free index ${}_{i_1}$. Finally the tensor fields
 indexed in $Z_b^2$ each have rank $a\ge 2$ and all free indices
  belong to the factors $\nabla\psi_1,\dots,\nabla\psi_\tau,\nabla Y$,
 $(\nabla\omega_1,\nabla\omega_2)$.
We may then re-write our equation (\ref{tolook}) in the form:

\begin{equation}
\label{tolook'}
\begin{split}
&\sum_{\zeta\in Z_a} a_\zeta X_{*}div_{i_1} C^{\zeta,i_1}_g +
\sum_{\zeta\in Z_b^1} a_\zeta  X_{*}div_{i_1} C^{\zeta,i_1}_g+
\\& \sum_{\zeta\in {Z_b^2}'} a_\zeta
 X_{*}div_{i_1}\dots X_{*}div_{i_a} C^{\zeta,i_1\dots i_a}_g
+\sum_{j\in J} a_j C^j_g=0,
\end{split}
\end{equation}
where now for the tensor fields indexed in ${Z_b^2}'$, each $a\ge
1$ and the factor $\nabla\psi_1$ does not contain a free index for
any of the tensor fields  for which $\nabla Y$ (or one of
$\nabla\omega_1,\nabla\omega_2$) is contracting against a bad
index in the chosen factor.

\par We will denote by $Z_{b,\sharp}^1\subset Z_b^1$ and
${Z_{b,\sharp}^2}'\subset {Z_b^2}'$ the index sets of tensor fields
 where $\nabla Y$ (or one of
$\nabla\omega_1,\nabla\omega_2$) is contracting
 against a bad index in the chosen factor $T$.

From (\ref{tolook'}) we derive an equation:

\begin{equation}
\label{mizera}
\begin{split}
&\sum_{\zeta\in Z_{b,\sharp}^1} a_\zeta  X_{**}div_{i_1}
C^{\zeta,i_1}_g+
\\& \sum_{\zeta\in {Z_{b,\sharp}^2}'} a_\zeta
 X_{**}div_{i_1}\dots X_{**}div_{i_a} C^{\zeta,i_1\dots i_a}_g
+\sum_{j\in J} a_j C^j_g=0,
\end{split}
\end{equation}
where $X_{**}div_i$ stands for the sublinear combination in
$X_{*}div_i$ for which $\nabla_i$ is in addition no allowed to hit
the chosen factor $T$.

\par Then, applying operation $Op$ as in Step 3 and the 
the inductive assumption of the general claim we are 
proving,\footnote{The resulting equation falls under the inductive assumption, as in Step 3.} 
and then using the operation $Op^{-1}[\dots]$ as in the proof of Step 3,
we derive a new equation:

\begin{equation}
\label{mizeracor}
\begin{split}
&\sum_{\zeta\in Z_{b,\sharp}^1} a_\zeta  X_{*}div_{i_1}
C^{\zeta,i_1}_g+  \sum_{\zeta\in {Z_{b,\sharp}^2}'} a_\zeta
 X_{*}div_{i_1}\dots X_{*}div_{i_a} C^{\zeta,i_1\dots i_a}_g=
 \\&\sum_{\zeta\in {Z_{OK}}} a_\zeta
 X_{*}div_{i_1}\dots X_{*}div_{i_a} C^{\zeta,i_1\dots i_a}_g
+\sum_{j\in J} a_j C^j_g=0,
\end{split}
\end{equation}
where the tensor fields indexed in $Z_{OK}$ have rank $a\ge 1$
(all free indices {\it not} belonging to factors
$\nabla\psi_1,\dots ,\nabla Y$ or
$\nabla\psi_1,\dots,\nabla\chi_2$) and furthermore have the
property that the one index in  $\nabla Y$ or $\nabla\omega_1$ is not
contracting against a bad index in the chosen factor (and it is
also not free). Thus, replacing the above back into (\ref{tolook'}),
 we derive:

\begin{equation}
\label{tolook''}
\begin{split}
&\sum_{\zeta\in Z_a} a_\zeta X_{*}div_{i_1} C^{\zeta,i_1}_g +
\sum_{\zeta\in {Z_b^1}'} a_\zeta  X_{*}div_{i_1} C^{\zeta,i_1}_g+
\\& \sum_{\zeta\in {Z_b^2}''} a_\zeta
 X_{*}div_{i_1}\dots X_{*}div_{i_a} C^{\zeta,i_1\dots i_a}_g
+\sum_{j\in J} a_j C^j_g=0,
\end{split}
\end{equation}
where the tensor fields indexed in ${Z_b^1}',{Z_b^2}''$ have the
additional restriction that if the factor $\nabla Y$ (or
$\nabla\omega_1,\nabla\omega_2$) is contracting against the chosen
factor $T$ then it is not contracting against a bad index in $T$.

\par We are now in a position to derive sub-step A from the above:
 To see this claim, we just apply
$Erase_{\nabla Y}$ or $Erase_{\nabla\omega_1}$ to (\ref{tolook''})
and multiply the resulting equation by
$\nabla_{i_1}Y\nabla^{i_1}\upsilon$.
\newline

{\it Sub-step B:} Now, we are reduced to showing our claim when
$Z_{a,spec}=\emptyset$. In that setting, we denote by
$Z_{a,s}\subset Z_a$ the index set of vector fields in $Z_a$ for
which the free index ${}_{i_1}$ belongs to the factor
$\nabla\psi_s$; we prove our claim separately for each of the
sublinear combinations $\sum_{\zeta\in Z_{a,s}}\dots$. This claim
is proven by picking out the sublinear combinations in
(\ref{procareless}), (\ref{careless}) where the factors
$\nabla\psi_s$ and $\nabla Y$ (or $\nabla\chi_1$) are contracting
against the same factor; we then apply the eraser to
$\nabla\psi_s$ (this is well-defined and produces a true
equation), and multiply by
$\nabla_{i_1}\psi_s\nabla^{i_1}\upsilon$. The resulting equation
is precisely our claim for the sublinear combination
$\sum_{\zeta\in Z_{a,s}}\dots$.
\newline

{\it (Sketch of the) Proof of Steps 1 and 2 (i.e. of (\ref{yirgachef}) and
 (\ref{yirgachef2})):} We will sketch the proof of these claims for the
sublinear combinations in \\$Z_a^{spec}\bigcup Z_b^{spec}\bigcup
 \overline{Z}^a_{spec}$ where one
of the special indices in $C^{\zeta,i_1}$ is an index ${}_k$ or
${}_l$ that belongs to a factor $S_{*}\nabla^{(\nu)}R_{ijkl}$. The
remaining case (where the special indices belong to factors
$\nabla^{(a)}\Omega_h$) can be seen by  a similar (simpler)
argument.\footnote{The only extra feature in this setting is that
one must prove the claim by a separate induction on the {\it
number} of factors $\nabla\psi_z$ that are contracting against
$\nabla^{(a)}\Omega_h$.}

 For each
$\zeta\in Z_a^{spec}\bigcup Z_b^{spec}\bigcup
 \overline{Z}^a_{spec}$ We denote by
$\overline{C}^{\zeta,{i_1}}_g$, $\overline{C}^{\zeta,i_1\dots i_\gamma}_g$
 the tensor fields that arise from $C^{\zeta,i_1}$
 $C^{\zeta,i_1\dots i_\gamma}_g$in (\ref{procareless}),
(\ref{procareless2})  by replacing the expressions
$S_{*}\nabla^{(\nu)}_{r_1\dots r_\nu}R_{ijkl}
\nabla^i\tilde{\phi}_1\nabla^k Y$, $S_{*}\nabla^{(\nu)}_{r_1\dots
r_\nu}R_{ijkl} \nabla^i\tilde{\phi}_1\nabla^k \chi_2$
 by a factor $\nabla^{(\nu+2)}_{r_1\dots r_\nu jl}\Omega_{b+1}$.
 We denote by $\tilde{\kappa}_{simp}$ the
 resulting simple character. We derive an equation:

\begin{equation}
\label{drivethere}
\begin{split}
&\sum_{\zeta\in Z_a^{spec}\bigcup Z_b^{spec}} a_\zeta
X_{*}div_{i_1} \overline{C}^{\zeta,i_1}_g +
\\&\sum_{\zeta\in \overline{Z}^{spec}_a} a_\zeta
X_{*}div_{i_1} \dots X_{*}div_{i_\gamma}\overline{C}^{\zeta,i_1}_g
+\sum_{j\in J} a_j
\overline{C}^j_g=0.
\end{split}
\end{equation}

\par Now, again applying the inductive assumption of our general statement,
 we derive that there is a linear combination of tensor fields (indexed in $W$ below)
with a free index ${}_{i_1}$ belonging to one of the factors
$\nabla\psi_1,\dots,\nabla\psi_\tau$ or
 $\nabla\psi_1,\dots,\nabla\psi_\tau, \nabla \chi_1$
  so that:

\begin{equation}
\label{fifteenth}
\begin{split}
& \sum_{\zeta\in Z_a^{spec}} a_\zeta
 \overline{C}^{\zeta,i_1}_g\nabla_{i_1}\upsilon
-X_{*}div_{i_2} \sum_{w\in W}
a_w C^{w,i_1i_2}_g\nabla_{i_1}\upsilon=\sum_{j\in J} a_j \overline{C}^j_g.
\end{split}
\end{equation}

\par Now, applying an operation $Op^{*}$ to the above which formally
 replaces the factor $\nabla^{(A)}_{r_1\dots r_A}\Omega_x$ by a factor
 $S_{*}\nabla^{(A-2)}_{r_1\dots r_{A-2}}R_{ir_{A-1}kr_A}\nabla^i\tilde{\phi}_1\nabla^kY$ or
 \\$S_{*}\nabla^{(A-2)}_{r_1\dots r_{A-2}}R_{ir_{A-1}kr_A}\nabla^i\tilde{\phi}_1\nabla^k\chi_2$,
   we derive (\ref{yirgachef}) (since we
  can repeat the permutations by which (\ref{fifteenth}) is made to
  hold formally, modulo introducing correction terms that allowed in the RHS of (\ref{yirgachef})).
\newline

\par We will now prove (\ref{yirgachef2}) by
repeating the induction performed in the ``Consequence'' we derived
above (where we showed that inductively assuming (\ref{olmert1}),
(\ref{olmert2}) we can derive (\ref{proolmert1'}),
(\ref{proolmert2'})):

\par We will show the claim of Step 2 in pieces: First consider
the tensor fields indexed in $\overline{Z}_{a,@}$ of minimum rank
$2$ (denote the corresponding index set is
$\overline{Z}_{a,@}^2$); we then show that we can write:

\begin{equation}
\label{synerghsw}
\begin{split}
&\sum_{\zeta\in \overline{Z}_{a,@}^2} a_\zeta Xdiv_{i_1}
Xdiv_{i_2} C^{\zeta,i_1 i_2}_g=
\\& \sum_{\zeta\in \overline{Z}_{a,@}^{3}} a_\zeta Xdiv_{i_1}\dots Xdiv_{i_{3}}
C^{\zeta,i_1\dots i_{3}}_g+\sum_{\zeta\in Z_{b,@}} a_\zeta
Xdiv_{i_1} C^{\zeta,i_1}_g+
\\&\sum_{\zeta\in Z_{OK}} a_\zeta Xdiv_{i_1}\dots
 Xdiv_{i_a}C^{\zeta,i_1\dots i_a}_g+
\sum_{j\in } a_j C^j.
\end{split}
\end{equation}
The tensor fields indexed in $\overline{Z}_{a,@}^{3}, Z_{b,@}$ in
the RHS are generic linear combinations in those forms (the first
with rank $3$). The tensor fields indexed in $Z_{OK}$ are generic
linear combinations as allowed in the RHS of our Step 2. Assuming we can prove (\ref{synerghsw}), 
we are
then reduced to showing our claim when the minimum rank among the
tensor fields indexed in $Z_{a,@}$ is 3. We may then ``forget''
about any $Xdiv_{i_h}$ where ${}_{i_h}$ belongs to the factor
$\nabla\psi_1$. Therefore, we are reduced to showing our claim
when the minimum rank is 2 and the factor $\nabla\psi_1$ does not
contain a free index. We then show our claim by an induction 
(for the rest of this derivation, 
 all tensor fields will {\it not} have a free index in the
factor $\nabla\psi_1$): Assume that the minimum rank of the tensor
fields indexed in $\overline{Z}_{a,@}$ is $k$, and they are
indexed in $\overline{Z}^k_{a,@}$. We then show that we can write:

\begin{equation}
\label{synerghsw'}
\begin{split}
&\sum_{\zeta\in \overline{Z}_{a,@}^k} a_\zeta Xdiv_{i_1}\dots
Xdiv_{i_k} C^{\zeta,i_1\dots i_k}_g=
\\& \sum_{\zeta\in \overline{Z}_{a,@}^{k+1}} a_\zeta Xdiv_{i_1}\dots Xdiv_{i_{k+1}}
C^{\zeta,i_1\dots i_{k+1}}_g+\sum_{\zeta\in Z_{b,@}} a_\zeta
Xdiv_{i_1} C^{\zeta,i_1}_g+
\\&\sum_{\zeta\in Z_{OK}} a_\zeta Xdiv_{i_1}\dots
 Xdiv_{i_a}C^{\zeta,i_1\dots i_a}_g+
\sum_{j\in } a_j C^j.
\end{split}
\end{equation}
The tensor fields indexed in $\overline{Z}_{a,@}^{3}, Z_{b,@}$ in
the the RHS are generic linear combinations in those forms (the
first with rank $k+1$). The tensor fields indexed in $Z_{OK}$ are
generic linear combinations as allowed in the RHS of our Step 2.

\par Iteratively repeating this step we are reduced
 to showing our Step 2 when $Z_{a,@}=\emptyset$.

\par In that case we then assume that the tensor fields indexed
 in $Z_{b,@}$ have minimum rank $k$ (and the corresponding index set is $Z^k_{b,@}$)
 and we show that we can write:

\begin{equation}
\label{synerghsw2}
\begin{split}
&\sum_{\zeta\in Z_{b,@}^k} a_\zeta Xdiv_{i_1}\dots Xdiv_{i_k}
C^{\zeta,i_1\dots i_k}_g=
\\& \sum_{\zeta\in Z_{b,@}^{k+1}} a_\zeta
Xdiv_{i_1}\dots Xdiv_{i_{k+1}}
C^{\zeta,i_1\dots i_{k+1}}_g+
\\&\sum_{\zeta\in Z_{OK}} a_\zeta Xdiv_{i_1}\dots
 Xdiv_{i_a}C^{\zeta,i_1\dots i_a}_g+\sum_{j\in } a_j C^j,
\end{split}
\end{equation}
 (with the same conventions as in the above equation).

\par If we can prove (\ref{synerghsw}) and (\ref{synerghsw2})
we will have shown our step 2.

{\it Proof of (\ref{synerghsw}), (\ref{synerghsw'}),
(\ref{synerghsw2}):} We start with a small remark: If the chosen factor is of the form
$S_{*}\nabla^{(\nu)}R_{ijkl}$, we replace our assumption by a more
convenient equation: Consider the tensor fields $C^{\zeta,i_1\dots
i_a}_g$, $\zeta\in \overline{Z}_{a,@}\bigcup Z_{b,@}$; we denote
by $\tilde{C}^{\zeta,i_1\dots i_a}_g$ the tensor fields that arise
 from $C^{\zeta,i_1\dots i_a}_g$ by replacing the expression
 $\nabla^{(\nu)}_{r_1\dots
 r_\nu}R_{ijkl}\nabla^i\tilde{\phi}_1\nabla^kY$ (or
$\nabla^{(\nu)}_{r_1\dots
r_\nu}R_{ijkl}\nabla^i\tilde{\phi}_1\nabla^k\chi_2$) by a factor
$\nabla^{(\nu+2)}_{r_1\dots r_\nu jl}\Omega_{p+1}$. We then derive
an equation:

\begin{equation}
\label{sfragis}
\begin{split}
&\sum_{\zeta\in \overline{Z}_a\bigcup Z_b} a_\zeta
X_{*}div_{i_1}\dots X_{*}div_{i_a} \tilde{C}^{\zeta,i_1\dots
i_a}_g(\Omega_1,\dots,\Omega_{p+1},\phi_2,\dots
,\phi_u,(\chi_1),\psi_1,\dots ,\psi_\tau)
\\&+\sum_{j\in J} a_j C^j_g(\Omega_1,\dots,\Omega_{p+1},\phi_2,\dots
,\phi_u,(\chi_1),\psi_1,\dots ,\psi_\tau).
\end{split}
\end{equation}
Now we can derive our claims:
\newline

{\it Proof of (\ref{synerghsw'}):} We divide the index set
$Z_{\overline{Z}_{a,@}^2}$ according to the two factors which
contain the two free indices  and we show our claim for each of
those tensor fields separately. The proof goes as follows: We pick
out the sublinear combination in our hypothesis (or in
(\ref{sfragis})) where the factors $\nabla\psi_h,\nabla\psi_{h'}$
(or $\nabla\psi_h,\nabla\chi_2$) are contracting against the same
factor. Clearly, this sublinear combination, $X_g$, vanishes
separately. We then formally erase the factor $\nabla\psi_h$.
Then, we may apply the inductive assumption of our general claim
to the resulting equation (the minimum rank of the tensor fields
will be 1), and (in case our assumption is (\ref{sfragis}) we also
apply an operation $Op^{-1}$ which replaces the factor
$\nabla^{(y)}_{r_1\dots r_y}\Omega_{p+1}$ by
$S_{*}\nabla^{(y-2)}_{r_1\dots
r_{y-2}}R_{ir_{y-1}kr_y}\nabla^i\tilde{\phi}_1\nabla^kY(\nabla^k\chi_1)$).
This is our desired conclusion.
\newline

{\it Proof of (\ref{synerghsw}), (\ref{synerghsw2}):} Now, we show
(\ref{synerghsw}) for the subset $Z_{a,@}^{k,\alpha}$ (which
indexes the $k$-tensor fields for which the free indices
${}_{i_1},\dots, {}_{i_k}$ belong to a chosen subset of the
factors $\nabla\psi_1,\dots,\nabla\psi_\tau, (\nabla\chi_1)$
(hence the label $\alpha$ designates the chosen subset). To
prove this equation, we pick out the sublinear combination in the
equation (\ref{sfragis}) where the factors
$\nabla\psi_2,\dots,\nabla\psi_\tau,(\nabla\chi_1)$ (indexed in
$\alpha$) are contracting against the same factor as
$\nabla\psi_1$. Then we apply the eraser to these factors and the
indices they contract against.  This is our desired conclusion. To
show (\ref{synerghsw2}), we only have to treat the factors
$\nabla\psi_h$ as factors $\nabla\phi_h$. The claim then follows
by applying Corollary 1 in \cite{alexakis4} and making the
 factors $\nabla\upsilon$ into $Xdiv$'s.\footnote{Observe that by virtue
of the factor $\nabla\psi_1$, we must have at least one non-simple
factor $S_{*}\nabla^{(\nu)}R_{ijkl}$ or $\nabla^{(B)}\Omega_h$
in (\ref{sfragis})--hence (\ref{sfragis}) does not fall under
any of the ``forbidden cases'' of Corollary 1 in \cite{alexakis4}, 
by inspection.} $\Box$
\newline

{\it Proof of the base case ($v+b=2$) of the general
 claim:} We firstly prove our claim when our hypothesis is (\ref{careless})
  (as opposed to (\ref{procareless})).

{\it Proof of the base case under the hypothesis
(\ref{careless}):} We observe that the weight $-K$ in our
assumption must satisfy $K\ge 2\tau +8$ if $v>0$ and $K\ge
2\tau+6$ if $v=0$.

\par First consider the case where we have the strict inequalities
$K> 2\tau +8$ if $v>0$ and $K> 2\tau+6$ if $v=0$. In that case our
first claim of the base case can be proven straightforwardly, by
picking out a removable index in each $C^{\zeta,i_a}_g, \zeta\in
Z_a$ and treating it as an $X_{*}div$ (which can be done when we
only have two real factors). Thus, in this setting we only have to
show our second claims (\ref{balladur}), (\ref{balladur'}).

In this setting, by using the ``manual'' constructions as in
\cite{alexakis3}, we can construct explicit tensor fields which
satisfy all the assumptions of our claim in the base case (each
with rank $\ge 2$), so that:

\begin{equation}
\label{colgate}
\begin{split}
&X_{+}div_{i_1} \sum_{\zeta\in Z'_a\bigcup Z_b} a_\zeta
C^{\zeta,i_1}_g(\Omega_1,\dots,\Omega_b,\phi_1,\dots,\phi_v,
[\chi_1,\chi_2],\psi_1,\dots,\psi_\tau)=
\\&\sum_{q\in Q} a_q X_{+}div_{i_1}C^{q,i_q}_g(\Omega_1,\dots,\Omega_b,\phi_1,\dots,\phi_v,
[\chi_1,\chi_2],\psi_1,\dots,\psi_\tau)+
\\&\sum_{p\in P} a_p X_{+}div_{i_1}\dots X_{+}div_{i_{c+1}} C^{p,i_1\dots
i_{c+1}}_g(\Omega_1,\dots,\Omega_b,\phi_1,\dots,\phi_v,
[\chi_1,\chi_2],\\&\psi_1,\dots,\psi_\tau)
+\sum_{j\in J} a_j C^j_g
\end{split}
\end{equation}
Here the tensor field $C^{p,i_1\dots i_{c+1}}_g$ will be in one of
 three forms:

If $v=2$ then each $C^{p,i_1\dots i_{c+1}}_g$ will be in the form:
\begin{equation}
\label{sadamhus}
\begin{split}
&pcontr(S_{*}{\nabla^{(\nu_1)}}^{f_{b_1}\dots f_{b_h}}_{i_1\dots
i_{c-1}} R_{x_1ji_cl}\otimes S_{*}{\nabla^{(\nu_2)}}^{f_{d_1}\dots
f_{d_y}}{{{R_{x_v}}^{j'}}_{i_{c+1}}}^l\otimes
\\&\otimes
[\nabla^j\chi_1\otimes\nabla_{j'}\chi_2]\otimes\nabla_{f_1}
\psi_1\dots\otimes\nabla_{f_\tau}\psi_\tau\otimes
\nabla^{x_1}\tilde{\phi}_1\otimes\nabla^{x_2}\tilde{\phi}_2),
\end{split}
\end{equation}
where $\{b_1,\dots ,b_h,d_1,\dots ,d_y\}=\{1,\dots,\tau\}$.

If $v=1$ then $\sum_{p\in P}\dots=0$ (this can be arranged because
of the two antisymmetric indices ${}_k,{}_l$ in the one factor
$S_{*}\nabla^{(\nu)}R_{ijkl}$).

If $v=0$ then each $C^{p,i_1\dots i_{c+1}}_g$ will be in the form:
\begin{equation}
\label{sadamhus2}
\begin{split}
&pcontr({\nabla^{(A_1)}}^{f_{b_1}\dots f_{b_h}}_{i_1\dots
i_{c-1}ji_c} \Omega_1\otimes {\nabla^{(A_2)}}^{f_{d_1}\dots
f_{d_y}j'i_{c+1}}\Omega_2
\\&\otimes
[\nabla^j\chi_1\otimes\nabla_{j'}\chi_2]\otimes\nabla_{f_1}
\psi_1\dots\otimes\nabla_{f_\tau}\psi_\tau\otimes
\nabla^{x_1}\tilde{\phi}_1\otimes\nabla^{x_2}\tilde{\phi}_2),
\end{split}
\end{equation}
where $\{b_1,\dots ,b_h,d_1,\dots ,d_y\}=\{1,\dots,\tau\}$.

\par Then, picking out the sublinear combination in
(\ref{sadamhus}), (\ref{sadamhus2}) with factors
$\nabla\psi_1,\dots,\nabla\psi_\tau,\nabla\chi_1,\nabla\chi_2$ we
derive that $\sum_{p\in P}\dots=0$. This is precisely our desired
conclusion in this case.
\newline

\par Now, the case where we have the equalities in our Lemma
hypothesis, $K= 2\tau +8$ if $v>0$ and $K= 2\tau+6$ if $v=0$. In
this case we note that in our hypothesis $Z_b=\emptyset$ if $v\ne
1$, while $Z_a=\overline{Z}_a=\emptyset$ if $v=1$.

\par Then, if $v\ne 1$, by the ``manual'' constructions as in \cite{alexakis3}, it
follows that we can construct tensor fields (as required in the
claim of our ``general claim''), so that:

\begin{equation}
\label{fronhsh}
\begin{split}
&\sum_{\zeta\in Z_a} a_\zeta
C^{\zeta,i_1}_g(\Omega_1,\dots,\Omega_b,\phi_1,\dots,\phi_v,
[\chi_1,\chi_2],\psi_1,\dots,\psi_\tau)\nabla_{i_1}\upsilon
\\&\-X_{*}div_{i_2} a_\zeta
C^{\zeta,i_1i_2}_g(\Omega_1,\dots,\Omega_b,\phi_1,\dots,\phi_v,
[\chi_1,\chi_2],\psi_1,\dots,\psi_\tau)\nabla_{i_1}\upsilon=
\\& a_{*}
C^{*,i_1}_g(\Omega_1,\dots,\Omega_b,\phi_1,\dots,\phi_v,
[\chi_1,\chi_2],\psi_1,\dots,\psi_\tau)\nabla_{i_1}\upsilon+
\\&\sum_{j\in J} a_j C^{j,i_1}_g(\Omega_1,\dots,\Omega_b,\phi_1,\dots,\phi_v,
[\chi_1,\chi_2],\psi_1,\dots,\psi_\tau)\nabla_{i_1}\upsilon,
\end{split}
\end{equation}
where the tensor field $C^{*,i_1}_g$ is in the form:

\begin{equation}
\label{invading}
\begin{split}
&pcontr(S_{*}{\nabla^{(\nu_1)}}^{f_1\dots f_{\tau-1}}
{{R_{x_1}}^{f_\tau}}_{kl}\otimes {{R_{x_2}}^{j'kl}}\otimes
\otimes\\&
[\nabla_{i_1}\chi_1\otimes\nabla_{j'}\chi_2]\otimes\nabla_{f_1}
\psi_1\dots\otimes\nabla_{f_\tau}\psi_\tau\otimes
\nabla^{x_1}\tilde{\phi}_1\otimes\nabla^{x_2}\tilde{\phi}_2),
\end{split}
\end{equation}
if $v=2$, and in the form:
\begin{equation}
\label{invading2}
\begin{split}
&pcontr({{\nabla^{(\tau+1)}}^{f_1\dots f_{\tau}}}_s\Omega_1
\otimes \nabla^{j's}\Omega_2
\otimes \\&[\nabla_{i_1}\chi_1\otimes\nabla_{j'}\chi_2]\otimes\nabla_{f_1}
\psi_1\dots\otimes\nabla_{f_\tau}\psi_\tau\otimes
\nabla^{x_1}\tilde{\phi}_1\otimes\nabla^{x_2}\tilde{\phi}_2),
\end{split}
\end{equation}
if $v=0$.

\par Thus, we are reduced to the case where $Z_a$ only consists of
the vector field (\ref{invading}) or (\ref{invading2}), and all
other tensor fields in our Lemma hypothesis have rank $\ge 2$ (we
have denoted their index set by $Z'_a$). We then show that we can
write:

\begin{equation}
\label{colgate2}
\begin{split}
&X_{+}div_{i_1} \sum_{\zeta\in Z'_a} a_\zeta
C^{\zeta,i_1}_g(\Omega_1,\dots,\Omega_b,\phi_1,\dots,\phi_v,
[\chi_1,\chi_2],\psi_1,\dots,\psi_\tau)=
\\&\sum_{q\in Q} a_q X_{+}div_{i_1}C^{q,i_1}_g(\Omega_1,\dots,\Omega_b,\phi_1,\dots,\phi_v,
[\chi_1,\chi_2],\psi_1,\dots,\psi_\tau)+
\\&\sum_{p\in P} a_p X_{+}div_{i_1}\dots X_{+}div_{i_{c+1}} C^{p,i_1\dots
i_{c+1}}_g(\Omega_1,\dots,\Omega_b,\phi_1,\dots,\phi_v,
[\chi_1,\chi_2],\psi_1,\dots,\psi_\tau)
\\&+\sum_{j\in J} a_j C^j_g,
\end{split}
\end{equation}
where the tensor fields indexed in $P$ here each have rank $\ge 2$
and are all in the form:

\begin{equation}
\label{sadamhus'}
\begin{split}
&pcontr(S_{*}{\nabla^{(\nu_1)}}^{f_1\dots f_{\tau-1}}
 {{R_{x_1}}^{f_\tau}}_{i_kl}\otimes
S_{*}{{R_{x_v}}^{j'kl}}\otimes
\\&[\nabla_{i_1}\chi_1\otimes\nabla_{j'}\chi_2]\otimes\nabla_{y_1}
\psi_1\dots\otimes\nabla_{y_\tau}\psi_\tau
\otimes\nabla^{x_1}\tilde{\phi}_1\otimes\nabla^{x_2}\tilde{\phi}_2),
\end{split}
\end{equation}

\begin{equation}
\label{sadamhus''}
\begin{split}
&pcontr({{\nabla^{(\nu_1)}}^{f_1\dots f_{\tau}}}_s
 \Omega_1\otimes
\nabla^{j's}\Omega_2\otimes
\\&[\nabla_{i_1}\chi_1\otimes\nabla_{j'}\chi_2]\otimes\nabla_{y_1}
\psi_1\dots\otimes\nabla_{y_\tau}\psi_\tau),
\end{split}
\end{equation}
where each of the indices ${}^{f_h}$ is contracting against one of
the indices ${}_{y_q}$. The indices ${}_{y_q}$ that are {\it not}
contracting against an index ${}^{f_h}$ are free indices.

\par Then, replacing the above into our Lemma hypothesis
(and making all the $\nabla\upsilon$'s ito $X_{+}div$'s), we
derive that $a_p=0$ for every $p\in P$ and $a_{*}=0$. This
concludes the proof of the base case when $v+b=2$, $v\ne 1$. In
the case $v=1$ we show our claim by just observing that we can
write:

\begin{equation}
\label{colgate3}
\begin{split}
&X_{+}div_{i_1} \sum_{\zeta\in Z_b} a_\zeta
C^{\zeta,i_1}_g(\Omega_1,\dots,\Omega_b,\phi_1,
[\chi_1,\chi_2],\psi_1,\dots,\psi_\tau)=
\\&\sum_{q\in Q} a_q X_{+}div_{i_1}C^{q,i_1}_g(\Omega_1,\dots,\Omega_b,\phi_1,
[\chi_1,\chi_2],\psi_1,\dots,\psi_\tau)
\\&+ \sum_{j\in J} a_j
C^j_g(\Omega_1,\dots,\Omega_b,\phi_1,
[\chi_1,\chi_2],\psi_1,\dots,\psi_\tau);
\end{split}
\end{equation}
this concludes the proof of the base case, when the tensor fields
in our Lemma hypothesis are in the form (\ref{careless}).
\newline

\par Now, we consider the setting where our hypothesis
is (\ref{procareless}). We again observe that if $v=0$
 then the weight $-K$ in our hypothesis must satisfy $K\ge 2\tau+4$.
  If $v>0$ it must satisfy $K\ge 2\tau+6$. We then again
  first consider the case where we have the strict inequalities
   in the hypothesis of our general claim.

\par In this case (where we have the strict inequalities
$K> 2\tau+4$ if $v=0$ and $K> 2\tau+6$ if $v\ne 0$) our first
claim follows straightforwardly (as above, we just pick out one
removable index in each $C^{\zeta,i_1}_g, \zeta\in Z_a$ and treat
it as an $X_{*}div$). To show the second claim we proceed much as
before:

\par We can ``manually'' construct tensor fields in order to write:

\begin{equation}
\label{colgate}
\begin{split}
&X_{+}div_{i_1} \sum_{\zeta\in Z'_a\bigcup Z_b} a_\zeta
C^{\zeta,i_1}_g(\Omega_1,\dots,\Omega_b,\phi_1,\dots,\phi_v,
Y,\psi_1,\dots,\psi_\tau)=
\\&\sum_{q\in Q} a_q X_{+}div_{i_1}C^{q,i_q}_g(\Omega_1,\dots,\Omega_b,\phi_1,\dots,\phi_v,
[\chi_1,\chi_2],\psi_1,\dots,\psi_\tau)+
\\&\sum_{p\in P} a_p X_{+}div_{i_1}\dots X_{+}div_{i_{c+1}} C^{p,i_1\dots
i_{c+1}}_g(\Omega_1,\dots,\Omega_b,\phi_1,\dots,\phi_v,
[\chi_1,\chi_2],\psi_1,\dots,\psi_\tau)
\\&+\sum_{j\in J} a_j C^j_g.
\end{split}
\end{equation}
Here the tensor field $C^{p,i_1\dots i_{c+1}}_g$ will be in one of
three forms:

If $v=2$ then each $C^{p,i_1\dots i_{c+1}}_g$ will be:
\begin{equation}
\label{sadamhus}
\begin{split}
&pcontr(S_{*}{\nabla^{(\nu_1)}}^{f_{b_1}\dots f_{b_h}}_{i_1\dots
i_{c-1}} {{R_{x_1}}^{f_{b_{h+1}}}}_{i_cl}\otimes
 S_{*}{\nabla^{(\nu_2)}}^{f_{d_1}\dots f_{d_y}}
 {{{R_{x_v}}^{f_{d_{y+1}}}}_{i_{c+1}}}^l\otimes
\\&\nabla_{f_{\tau+1}}Y\otimes\nabla_{f_1}
\psi_1\dots\otimes\nabla_{f_\tau}\psi_\tau\otimes
\nabla^{x_1}\tilde{\phi}_1\otimes\nabla^{x_2}\tilde{\phi}_2),
\end{split}
\end{equation}
where $\{b_1,\dots ,b_{h+1},d_1,\dots ,d_{y+1}\}=\{1,\dots,\tau+1\}$.

If $v=1$ then $\sum_{p\in P}\dots=0$ (this is because of the two
antisymmetric indices ${}_k,{}_l$ in the one factor
$S_{*}\nabla^{(\nu)}R_{ijkl}$).

If $v=0$ then each $C^{p,i_1\dots i_{c+1}}_g$ will be in the form:
\begin{equation}
\label{sadamhus2}
\begin{split}
&pcontr({\nabla^{(A_1)}}^{f_{b_1}\dots f_{b_{h}}}_{i_1\dots
i_{c-1}i_c}\Omega_1 \otimes {\nabla^{(A_2)}}^{f_{d_1}\dots
f_{d_{y}}i_{c+1}}\Omega_2
\\&\otimes\nabla_{f_{\tau+1}} Y\otimes\nabla_{f_1}
\psi_1\dots\otimes\nabla_{f_\tau}\psi_\tau),
\end{split}
\end{equation}
where $\{b_1,\dots ,b_h,d_1,\dots ,d_y\}=\{1,\dots,\tau+1\}$.

\par Then, picking out the sublinear combination in
(\ref{sadamhus}), (\ref{sadamhus2}) with factors
$\nabla\psi_1,\dots,\nabla\psi_\tau,\nabla Y$ we
derive that $\sum_{p\in P}\dots=0$. This is precisely our desired
conclusion in this case.
\newline

\par Finally, we prove our claim when we have the equalities
$K= 2\tau+4$ if $v<2$ and $K= 2\tau+6$ if $v =2$) in the
hypothesis of our general claim.

\par In this case by ``manually'' constructing $X_{+}div$'s
 so that we can write:

\begin{equation}
\label{sygkino}
\begin{split}
&\sum_{\zeta\in Z'_a\bigcup Z_b\bigcup \overline{Z}_a} a_\zeta
X_{+}div_{i_1}\dots X_{+}div_{i_a} C^{\zeta,i_1\dots
i_a}_g(\Omega_1,\dots \Omega_b,Y,\psi_1,\dots,\psi_\tau)=
\\&\sum_{q\in Q} a_q X_{+}div_{i_1}\dots X_{+}div_{i_a}
C^{q,i_1\dots i_a}_g(\Omega_1,\dots
\Omega_b,Y,\psi_1,\dots,\psi_\tau)+
\\&\sum_{p\in P} a_p X_{+}div_{i_1}\dots X_{+}div_{i_a}
C^{p,i_1\dots i_a}_g(\Omega_1,\dots
\Omega_b,Y,\psi_1,\dots,\psi_\tau)
\\&+\sum_{j\in J} a_j
C^j_g(\Omega_1,\dots \Omega_b,Y,\psi_1,\dots,\psi_\tau).
\end{split}
\end{equation}

Here the tensor fields indexed in $P$ are in the following form:

If $v=0$ then they will either be in the form:

\begin{equation}
\label{antestrateyeto1}
\begin{split}
&pcontr(\nabla_{i_{*}}Y\otimes {\nabla^{(A)}}^{f_{x_1}\dots f_{x_a}
s}\Omega_1 \otimes{\nabla^{(B)}}^{f_{x_{a+1}}\dots f_{x_\tau}}_s\Omega_2
\otimes\nabla_{f_1}\psi_1\otimes\dots\otimes\nabla_{f_\tau}\phi_\tau),
\end{split}
\end{equation}
(where $\{x_1,\dots ,x_\tau\}=\{1,\dots ,\tau\}$), or in the form:

\begin{equation}
\label{antestrateyeto2}
\begin{split}
&pcontr(\nabla_{q}Y\otimes
{\nabla^{(A)}}^{f_{x_1}\dots f_{x_a}}_{i_{*}}\Omega_1
\otimes{\nabla^{(B)}}^{f_{x_{a+1}}\dots f_{x_\tau} q}\Omega_2
\otimes\nabla_{f_1}\psi_1\otimes\dots\otimes\nabla_{f_\tau}\phi_\tau),
\end{split}
\end{equation}
(where $\{x_1,\dots ,x_\tau\}=\{1,\dots ,\tau\}$).

\par If $v=2$ they will be in the form:

\begin{equation}
\label{antestrateyeto3}
\begin{split}
&pcontr(\nabla_{i_{*}}Y\otimes {\nabla^{(A)}}^{f_{x_1}\dots
f_{x_{a-1}}} S_{*}R^{i f_{x_a} kl} \otimes {\nabla^{(B)}}^{f_{x_{a+1}}\dots f_{x_{\tau-1}}} {R^{i'f_{x_\tau} }}_{kl}
\\&\otimes\nabla_{f_1}\psi_1\otimes\dots\otimes\nabla_{f_\tau}
\phi_\tau\nabla_i\tilde{\phi}_1\otimes\nabla_{i'}\tilde{\phi}_2),
\end{split}
\end{equation}
(where $\{x_1,\dots ,x_\tau\}=\{1,\dots ,\tau\}$), or in the form:

\begin{equation}
\label{antestrateyeto4}
\begin{split}
&pcontr(\nabla_{q}Y\otimes {\nabla^{(A)}}^{f_{x_1}\dots
f_{x_{a-1}}}
S_{*}R^{i f_{x_a} ql} \otimes
{\nabla^{(B)}}^{f_{x_{a+1}}\dots f_{x_{\tau-1}}} {R^{i'f_{x_\tau}}}_{i_{*}l}
\\&\otimes\nabla_{y_1}\psi_1\otimes\dots\otimes\nabla_{y_\tau}
\phi_\tau\nabla_i\tilde{\phi}_1\otimes\nabla_{i'}\tilde{\phi}_2).
\end{split}
\end{equation}

\par If $v=1$ the equation (\ref{sygkino}) will hold with $P=\emptyset$:

\par Then, picking out the sublinear combination in
(\ref{sygkino}) which consists of terms with a factor $\nabla Y$
and replacing into our hypothesis, we derive that the coefficient
of each of the tensor fields indexed in $P$ must be zero. This
completes the proof of our claim. $\Box$

\subsection{Proof of Lemmas \ref{obote3}, \ref{vanderbi3}:}

{\it Proof of Lemma \ref{obote3}:}
 \newline

 The first claim follows
immediately, since each tensor field has a removable index (thus
each tensor field separately can be written as an $X_{*}div$).

The proof of the second claim essentially follows the ``manual''
construction of divergences, as in \cite{alexakis3}. By ``manually''
 constructing explicit divergences out of each $C^{h,i_1\dots
i_\alpha}_g(\Omega_1,\dots,\Omega_p,\phi_1,\dots,\phi_u)$, $h\in H_2$,
 we derive that we can write:

\begin{equation}
\label{amerwom}
\begin{split}
&\sum_{h\in H_2} a_h Xdiv_{i_1}\dots Xdiv_{i_\alpha}C^{h,i_1\dots
i_\alpha}_g(\Omega_1,\dots,\Omega_p,Y,\phi_1,\dots,\phi_u)=
\\&(Const)_1 Xdiv_{i_1}\dots Xdiv_{i_\xi}C^{1,i_1\dots
i_\xi}_g(\Omega_1,\dots,\Omega_p,Y,\phi_1,\dots,\phi_u)+
\\&(Const)_2 Xdiv_{i_1}\dots Xdiv_{i_\zeta}C^{2,i_1\dots
i_\zeta}_g(\Omega_1,\dots,\Omega_p,Y,\phi_1,\dots,\phi_u)+
\\&\sum_{q\in Q} a_q Xdiv_{i_1}\dots Xdiv_{i_\alpha}C^{q,i_1\dots
i_\alpha}_g(\Omega_1,\dots,\Omega_p,Y,\phi_1,\dots,\phi_u)+
\\&\sum_{j\in J} a_j C^{j}_g(\Omega_1,\dots,\Omega_p,Y,\phi_1,\dots,\phi_u),
\end{split}
\end{equation}
where the tensor fields indexed in $Q$ are as required by our
Lemma hypothesis, while the tensor fields $C^1,C^2$ are explicit
tensor fields which we will write out below (they depend on the
values $p,\sigma_1,\sigma_2$).\footnote{In some cases there will be no
tensor fields $C^1,C^2$ (in which case we will just say that in
(\ref{amerwom}) we have $(Const)_1=0$, $(Const)_2=0$).}

\par We will then show that in (\ref{amerwom}) we will have
$(Const)_1=(Const)_2=0$. That will complete the proof of Lemma
\ref{obote3}. We distinguish cases based on the value of $p$: Either $p=2$ or $p=1$ or $p=0$. 
\newline

 {\it The case $p=2$:} With no loss of
generality we assume that the factor $\nabla^{(A)}\Omega_1$ is
contracting against the factors $\nabla\phi_1,\dots,\nabla\phi_x$
and $\nabla^{(B)}\Omega_2$ is contracting against
$\nabla\phi_{x+1},\dots ,\nabla\phi_{x+t}$; we may also assume
wlog that $x\le t$. By  manually constructing divergences, it
follows that we can derive (\ref{amerwom}), where each of the
tensor fields $C^1,C^2$ will be in the forms, respectively:

\begin{equation}
\label{mitsakis} pcontr(\nabla_{i_{*}}Y\otimes
\nabla^{(A)}_{v_1\dots v_xi_1\dots
i_\gamma}\Omega_1\otimes\nabla^{(B)}_{y_1\dots y_t
i_{\gamma+1}\dots i_{\gamma+\delta}}\Omega_2
\otimes\nabla^{v_1}\phi_1\otimes\dots\otimes\nabla^{y_t}\phi_u),
\end{equation}
(where if $t\ge 2$ then $\delta=0$, otherwise $t+\delta=2$), or:

\begin{equation}
\label{mitsakis2} pcontr(\nabla_{q}Y\otimes
\nabla^q\nabla^{(A)}_{v_1\dots v_xi_1\dots
i_\gamma}\Omega_1\otimes\nabla^{(B)}_{y_1\dots y_t
i_{\gamma+1}\dots i_{\gamma+\delta}}\Omega_2
\otimes\nabla^{v_1}\phi_1\otimes\dots\otimes\nabla^{y_t}\phi_u),
\end{equation}
(where if $t\ge 2$ then $\delta=0$, otherwise $t+\delta=2$).

{\it The case $p=1$:}  We ``manually'' construct divergences to derive
(\ref{amerwom}), where  if $\sigma_1=1$ then there are no tensor
 fields $C^1,C^2$ (and hence (\ref{amerwom}) is our desired conclusion);
if $\sigma_1=0,\sigma_2=1$ then
there is only the tensor field $C^1$ in (\ref{amerwom}) and it is in the form:

\begin{equation}
\label{mitsakis3}\begin{split} &pcontr(\nabla^{q}Y\otimes
S_{*}\nabla^{(\nu)}_{v_2\dots v_xi_1\dots
i_\gamma}R_{ii_{\gamma+1}i_{\gamma+2}q}\otimes\nabla^{(B)}_{y_1\dots
y_t i_{\gamma+1}\dots i_{\gamma+\delta}}\Omega_2
\\&\otimes\nabla^i\tilde{\phi}_1\otimes\nabla^{v_1}\phi_2\otimes\dots\otimes\nabla^{y_t}\phi_u),
\end{split}
\end{equation}
where if $t\ge 2$ then $\delta=0$, otherwise $\delta=2-t$.
\newline

{\it The case $p=0$:} We have three subcases: Firstly
$\sigma_2=2$, secondly $(\sigma_2=1,\sigma_1=1)$, and thirdly
$\sigma_1=2$.

In the case $\sigma_2=2$, the tensor fields $C^1,C^2$ must be in the forms,
respectively:

\begin{equation}
\label{mitsakis3,5}
\begin{split} & pcontr(\nabla_{i_{*}}Y\otimes
S_{*}\nabla^{(\nu)}_{v_2\dots v_xi_1\dots
i_\gamma}R_{ii_{\gamma+1}i_{\gamma+2}l}\otimes
S_{*}\nabla^{(t-1)}_{y_1\dots
y_t}{R_{i'i_{\gamma+3}i_{\gamma+4}}}^l
\\&\otimes\nabla^i\tilde{\phi}_1\otimes\nabla^{i'}\tilde{\phi}_2
\otimes\nabla^{v_1}\phi_3\otimes\dots\otimes\nabla^{y_t}\phi_u),
\end{split}
\end{equation}

\begin{equation}
\label{mitsakis4}\begin{split} & pcontr(\nabla^{q}Y\otimes
S_{*}\nabla^{(\nu)}_{qv_2\dots v_xi_1\dots
i_\gamma}R_{ii_{\gamma+1}i_{\gamma+2}l}\otimes
S_{*}\nabla^{(t-1)}_{y_1\dots
y_t}{R_{i'i_{\gamma+3}i_{\gamma+4}}}^l
\\&\otimes\nabla^i\tilde{\phi}_1\otimes\nabla^{i'}\tilde{\phi}_2
\otimes\nabla^{v_1}\phi_3\otimes\dots\otimes\nabla^{y_t}\phi_u),
\end{split}
\end{equation}
(if $x=t=0$ then the tensor field $C^1$ above will not be
present).

In the case $\sigma_1=2$, the tensor fields $C^1,C^2$ must be in one of the
two forms:

\begin{equation}
\label{mitsakis5} pcontr(\nabla_{i_{*}}Y\otimes
\nabla^{(m_1)}_{v_1\dots v_xi_1\dots
i_\gamma}R_{ii_{\gamma+1}i_{\gamma+2}l}\otimes
\nabla^{(t-1)}_{y_1\dots y_t}{{R^i}_{i_{\gamma+3}i_{\gamma+4}}}^l
\otimes\nabla^{v_1}\phi_1\otimes\dots\otimes\nabla^{y_t}\phi_u),
\end{equation}

\begin{equation}
\label{mitsakis6} \begin{split} &pcontr(\nabla_{q}Y\otimes
\nabla^q\nabla^{(m_1)}_{v_1\dots v_xi_1\dots
i_\gamma}R_{ii_{\gamma+1}i_{\gamma+2}l}\otimes
\nabla^{(t-1)}_{y_1\dots y_t}{{R^i}_{i_{\gamma+3}i_{\gamma+4}}}^l
\\&\otimes\nabla^{v_1}\phi_1\otimes\dots\otimes\nabla^{y_t}\phi_u).
\end{split}\end{equation}

In the case $\sigma_1=1,\sigma_2=1$, there will be only one tensor
field $C^1$, in the form:

\begin{equation}
\label{mitsakis7} \begin{split}&pcontr(\nabla^qY\otimes
S_{*}\nabla^{(m_1)}_{v_1\dots v_xi_1\dots
i_\gamma}R_{ii_{\gamma+1}i_{\gamma+2}l}\otimes
\nabla^{(t-1)}_{y_1\dots y_t}{R_{qi_{\gamma+3}i_{\gamma+4}}}^l
\otimes\nabla^i\tilde{\phi}_1\\&\otimes\nabla^{v_1}\phi_2\otimes\dots\otimes\nabla^{y_t}\phi_u).
\end{split}\end{equation}

\par We then derive that $(Const)_1=(Const)_2=0$ as in \cite{alexakis3}
 (by picking out the sublinear combination in
(\ref{amerwom}) that consists of  complete contractions with a
factor $\nabla Y$--differentiated only once).
\newline

{\it Proof of Lemma \ref{vanderbi3}:}
\newline

\par We again ``manually'' construct explicit $Xdiv$' to write:

\begin{equation}
\label{amerwom'}
\begin{split}
&\sum_{h\in H_2} a_h Xdiv_{i_1}\dots Xdiv_{i_\alpha}C^{h,i_1\dots
i_\alpha}_g(\Omega_1,\dots,\Omega_p,Y,\phi_1,\dots,\phi_u)=
\\&(Const)_1 Xdiv_{i_1}\dots Xdiv_{i_\xi}C^{1,i_1\dots
i_\xi}_g(\Omega_1,\dots,\Omega_p,Y,\phi_1,\dots,\phi_u)+
\\&(Const)_2 Xdiv_{i_1}\dots Xdiv_{i_\zeta}C^{2,i_1\dots
i_\zeta}_g(\Omega_1,\dots,\Omega_p,Y,\phi_1,\dots,\phi_u)+
\\&\sum_{q\in Q} a_q Xdiv_{i_1}\dots Xdiv_{i_\alpha}C^{q,i_1\dots
i_\alpha}_g(\Omega_1,\dots,\Omega_p,Y,\phi_1,\dots,\phi_u)+
\\&\sum_{j\in J} a_j C^{j}_g(\Omega_1,\dots,\Omega_p,Y,\phi_1,\dots,\phi_u),
\end{split}
\end{equation}
where the tensor fields indexed in $Q$ are as required by our
Lemma hypothesis, while the tensor fields $C^1,C^2$ are explicit
tensor fields which we will write out below (they depend on the
values $p,\sigma_1,\sigma_2$). In some cases there will be no
tensor fields $C^1,C^2$ (in which case we will just say that in
(\ref{amerwom}) we have $(Const)_1=0$, $(Const)_2=0$).

 {\it The case $p=2$:} With no loss of
generality we assume that the factor $\nabla^{(A)}\Omega_1$ is
contracting against the factors $\nabla\phi_1,\dots,\nabla\phi_x$
and $\nabla^{(B)}\Omega_2$ is contracting against
$\nabla\phi_{x+1},\dots ,\nabla_{\phi_{x+t}}$; we may also assume
wlog that $x\le t$. By manual construction of divergences, it
follows that we can derive (\ref{amerwom}), where there is only the tensor field $C^1$
and it is in the form:

\begin{equation}
\label{bmitsakis}
pcontr(\nabla_{[i_{*}}\chi_1\otimes\nabla^{q]}\chi_2\otimes
\nabla^{(A)}_{v_1\dots v_xi_1\dots
i_\gamma}\Omega_1\otimes\nabla^{(B)}_{qy_1\dots y_t
i_{\gamma+1}\dots i_{\gamma+\delta}}\Omega_2
\otimes\nabla^{v_1}\phi_1\otimes\dots\otimes\nabla^{y_t}\phi_u),
\end{equation}
(where if $t\ge 1$ then $\delta=0$, otherwise $\delta=1$).
\newline

{\it The case $p=1$:}  We ``manually'' construct divergences to derive
(\ref{amerwom'}), where:  if $\sigma_1=1$ then there are no tensor
fields $C^1,C^2$ in the RHS of (\ref{amerwom'}) (and this is our desired conclusion);
 if $\sigma_1=0,\sigma_2=1$ then
there is only the tensor field $C^1$ in the RHS of (\ref{amerwom'}) and it
is of the form:

\begin{equation}
\label{bmitsakis3}
\begin{split}
&pcontr(\nabla_{[i_{*}}\omega_1\otimes\nabla^{q]}\omega_2\otimes
S_{*}\nabla^{(\nu)}_{v_2\dots v_xi_1\dots
i_\gamma}R_{ii_{\gamma+1}i_{\gamma+2}q}\otimes\nabla^{(B)}_{y_1\dots
y_t i_{\gamma+1}\dots i_{\gamma+\delta}}\Omega_2
\otimes\nabla^i\tilde{\phi}_1
\\&\otimes\nabla^{v_1}\phi_2\otimes\dots\otimes\nabla^{y_t}\phi_u),
\end{split}
\end{equation}
where if $t\ge 2$ then $\delta=0$, otherwise $\delta=2-t$.

{\it The case $p=0$:} We have three subcases: Firstly
$\sigma_2=2$, secondly  $(\sigma_2=1,\sigma_1=1)$, and thirdly
$\sigma_1=2$.

In the case $\sigma_2=2$, the tensor fields $C^1,C^2$ in the RHS of (\ref{amerwom})
will be in the two forms, respectively:

\begin{equation}
\label{bmitsakis4}
\begin{split}
&pcontr(\nabla_{[i_{*}}\omega_1\otimes\nabla^{q]}\omega_2\otimes
S_{*}\nabla^{(\nu)}_{qv_2\dots v_xi_1\dots
i_\gamma}R_{ii_{\gamma+1}i_{\gamma+2}l}\otimes
\\&S_{*}\nabla^{(t-1)}_{y_1\dots
y_t}{R_{i'i_{\gamma+3}i_{\gamma+4}}}^l
\otimes\nabla^i\tilde{\phi}_1\otimes\nabla^{i'}\tilde{\phi}_2
\otimes\nabla^{v_1}\phi_3\otimes\dots\otimes\nabla^{y_t}\phi_u),
\end{split}
\end{equation}

\begin{equation}
\label{bmitsakis5}
\begin{split}
&pcontr(\nabla_{[p}\omega_1\otimes\nabla^{q]}\omega_2\otimes
S_{*}\nabla^{(\nu)}_{qv_2\dots v_xi_1\dots
i_\gamma}R_{ii_{\gamma+1}i_{\gamma+2}p}\otimes
\\&S_{*}\nabla^{(t-1)}_{y_1\dots
y_t}{R_{i'i_{\gamma+3}i_{\gamma+4}q}}
\otimes\nabla^i\tilde{\phi}_1\otimes\nabla^{i'}\tilde{\phi}_2
\otimes\nabla^{v_1}\phi_3\otimes\dots\otimes\nabla^{y_t}\phi_u).
\end{split}
\end{equation}

In the case $\sigma_1=2$, the tensor fields $C^1,C^2$ will be  the
forms, respectively:

\begin{equation}
\label{bmitsakis5}
\begin{split}
&pcontr(\nabla_{[i_{*}}\omega_1\otimes\nabla^{q]}\otimes
\nabla^{(m_1)}_{v_1\dots v_xi_1\dots
i_\gamma}R_{ii_{\gamma+1}i_{\gamma+2}l}\otimes
\nabla^{(t-1)}_{qy_1\dots y_t}{{R^i}_{i_{\gamma+3}i_{\gamma+4}}}^l
\\&\otimes\nabla^{v_1}\phi_1\otimes\dots\otimes\nabla^{y_t}\phi_u),
\end{split}
\end{equation}

\begin{equation}
\label{bmitsakis6}
\begin{split}
&pcontr(\nabla_{[p}\omega_1\otimes\nabla^{q]}\omega_2\otimes
\nabla^{(m_1)}_{v_1\dots v_xi_1\dots
i_\gamma}R_{ii_{\gamma+1}i_{\gamma+2}p}\otimes
\nabla^{(t-1)}_{y_1\dots y_t}{{R^i}_{i_{\gamma+3}i_{\gamma+4}q}}
\\&\otimes\nabla^{v_1}\phi_1\otimes\dots\otimes\nabla^{y_t}\phi_u),
\end{split}
\end{equation}
(if at least one of the two factors $\nabla^{(m)}R_{ijkl}$ is
contracting against a factor $\nabla\phi_h$. Otherwise, we can prove (\ref{amerwom'})
 with no tensor fields $C^1,C^2$ on the RHS).

In the case $\sigma_1=1,\sigma_2=1$, the tensor fields $C^1,C^2$ must be in
 the forms, respectively:

\begin{equation}
\label{bmitsakis7}
\begin{split}
&pcontr(\nabla_{[i_{*}}\omega_1\otimes\nabla^{q]}\omega_2\otimes
S_{*}\nabla^{(\nu)}_{v_1\dots v_xi_1\dots
i_\gamma}R_{ii_{\gamma+1}i_{\gamma+2}l}\otimes
\nabla^{(t-1)}_{y_1\dots y_t}{R_{qi_{\gamma+3}i_{\gamma+4}}}^l
\\&\otimes\nabla^i\tilde{\phi}_1\otimes\nabla^{v_1}\phi_2
\otimes\dots\otimes\nabla^{y_t}\phi_u),
\end{split}
\end{equation}

\begin{equation}
\label{bmitsakis8}
\begin{split}
&pcontr(\nabla_{[p}\omega_1\otimes\nabla^{q]}\omega_2\otimes
S_{*}\nabla^{(m_1)}_{v_1\dots v_xi_1\dots
i_\gamma}R_{ii_{\gamma+1}i_{\gamma+2}l}\otimes
\nabla^{(t-1)}_{y_1\dots y_t}{R_{pqi_{\gamma+3}}}^l
\\&\otimes\nabla^i\tilde{\phi}_1\otimes\nabla^{v_1}
\phi_2\otimes\dots\otimes\nabla^{y_t}\phi_u).
\end{split}
\end{equation}

\par We then derive that $(Const)_1=(Const)_2=0$ by picking out the sublinear combination in
(\ref{amerwom'}) that consists of  complete contractions with two
factors $\nabla Y,\nabla\omega_2$--each factor differentiated only once). $\Box$

\section{The proof of Proposition \ref{giade} in the special cases:}
\label{specialcases}

\subsection{The direct proof  of Proposition \ref{giade} (in case II)
in the ``special cases''.}

We now prove Proposition \ref{giade} directly in the special subcases of case II. 
We recall the setting of the special subcases of Proposition \ref{giade} 
in case II  are as follows: 
In subcase IIA
for each $\mu$-tensor field of maximal refined double character,
 $C^{l,i_1\dots i_\mu}_g$ there is a unique factor 
in the form $T=\nabla^{(m)}R_{ijkl}$ for which two internal indices are free, and each 
derivative index is either free or contracting against a 
factor $\nabla\phi_h$. For subcase IIB there is a unique factor  
in the form $T=\nabla^{(m)}R_{ijkl}$ for which one internal index is free, and each 
derivative index is either free or contracting against a 
factor $\nabla\phi_h$. In both sucases IIA, IIB there is at 
least one free derivative index in the factor $T$. 

Moreover, both in subcases IIA, IIB, all {\it other} real factors
in $C^{l,i_1\dots i_\mu}_g$ are either  
in the form $S_{*}R_{ijkl}$ or $\nabla^{(2)}\Omega_h$, or they are in the form $\nabla^{(m)}R_{ijkl}$,
where all the $m$ derivative indices contract against 
factors $\nabla\phi_h$.\footnote{For the rest of this subsection, we will 
slightly abuse notation and {\it not} write out the derivative indices that contract 
against factors $\nabla\phi_h$--we will thus refer to factors $R_{ijkl}$, setting $m=0$.}
\newline

In order to prove Proposition \ref{giade} directly in the special subcases of subcases IIA, IIB 
we will rely on a new Lemma:

Our new Lemma deals with two different settings, which we will 
label setting A and setting B below.

In setting A, we let
$$ \sum_{l\in \overline{L}} a_l C^{l,i_1\dots i_\mu}_g
(\Omega_1,\dots,\Omega_p,\phi_1,\dots,\phi_u)$$ 
stand for a linear combination of $\mu$-tensor fields with one factor $\nabla^{(m)}R_{ijkl}$
containing  $\alpha\ge 2$ free indices, distributed according to the pattern 
\\$\nabla^{(m)}_{(free)\dots (free)}R_{(free)j(free)l}$, and all other other factors  
being  all in one of the forms $R_{ijkl}$,$S_{*}R_{ijkl}$,$\nabla^{(2)}\Omega_h$. 
(I.e.~they have no removable indices).

In setting B we let
$$\sum_{l\in \overline{L}} a_l C^{l,i_1\dots i_\mu}_g
(\Omega_1,\dots,\Omega_p,\phi_1,\dots,\phi_u)$$ 
stand for a linear combination of $\mu$-tensor fields with one factor $\nabla^{(m)}R_{ijkl}$
containing  $\alpha\ge 2$ free indices, distributed according to the pattern 
\\$\nabla^{(m)}_{(free)\dots (free)}R_{(free)j(free)l}$, 
and all {\it but one} of the other factors  
being   in one of the forms $R_{ijkl}$,$S_{*}R_{ijkl}$,
$\nabla^{(2)}\Omega_h$; one of the other factors (which we label $T'$) 
will be in the form 
$\nabla R_{ijkl}$,$S_{*}\nabla R_{ijkl}$, $\nabla^{(3)}\Omega_h$. We will 
call this other factor ``the factor with the extra derivative''. 
Moreover, in setting B we impose the additional
 restriction that if both the indices ${}_j,{}_l$ 
in the factor $\nabla^{(m)}_{(free)\dots (free)}R_{(free)j(free)l}$  
contract against the same other factor $T'$, then either $T'$ is {\it not} the 
factor with the extra derivative, or if it is, then $T'$ is in the form $\nabla_s R_{abcd}$,
and furthermore  the indices ${}_j,{}_l$ contract against the 
indices ${}_b,{}_c$ {\it and} we assume that the indices ${}_s,{}_a,{}_c$ 
are symmetrized over.\footnote{In other words, in that case the factors $T,T'$
contract according to the pattern: $\nabla^{(m)}_{(free)\dots (free)}R_{(free)j(free)l}\nabla_{(s}{{R_a}^{jk}}_{d)}$,
 where the indices ${}_s,{}_a,{}_d$ are symmetrized over.}

\begin{lemma}
\label{specialii}
Let $\sum_{l\in \overline{L}} a_l C^{l,i_1\dots i_\mu}_g$ 
be a linear combination of $\mu$-tensor fields as described above. 
We assume the following special case of (\ref{hypothese2}):

\begin{equation}
\label{gates} 
\begin{split}
&\sum_{l\in \overline{L}\bigcup L'} a_l Xdiv_{i_1}\dots Xdiv_{i_\mu} 
C^{l,i_1\dots i_\mu}_g(\Omega_1,\dots,\Omega_p,\phi_1,\dots,\phi_u)+
\\&\sum_{h\in H} a_h Xdiv_{i_1}\dots Xdiv_{i_{\beta}}
C^{h,i_1\dots i_\beta}_g(\Omega_1,\dots,\Omega_p,\phi_1,\dots,\phi_u)
\\&+\sum_{j\in J} a_j C^j_g(\Omega_1,\dots,\Omega_p,\phi_1,\dots,\phi_u);
\end{split}
\end{equation}
here, in both cases A and B the terms indexed in $\overline{L}$ 
will be as described above; the $\mu$-tensor fields indexed in 
$L'$ will have fewer
than $\alpha$ free indices in any given factor of the form $\nabla^{(m)}R_{ijkl}$. 
The tensor felds indexed in $H$ each have rank $>\mu$ 
and also each of them has fewer than $\alpha$ free indices
in any given factor of the form $\nabla^{(m)}R_{ijkl}$. 
Finally, the terms indexed in $J$
are simply subsequent to $\vec{\kappa}_{simp}$. 
\newline
We claim that:
\begin{equation}
\label{tasidera}
\sum_{l\in \overline{L}} a_l C^{l,i_1\dots i_\mu}_g\nabla_{i_1}\upsilon\dots \nabla_{i_\mu}\upsilon=0.
\end{equation}
\end{lemma}

\par We will prove this Lemma shortly. Let us now, however, note how the above 
Lemma directly implies Proposition \ref{giade} in the 
special subcases IIA (directly) and IIB (after some manipulation).
 \newline

{\it Lemma \ref{specialii} implies Proposition \ref{giade} 
in the special subcases of case II:}
 
We first start with subcase IIA: 
 Consider the sublinear combination of $\mu$-tensor fields of 
 maximal refined double character in (\ref{hypothese2}). 
 Denote their index set by $L_{Max}\subset L$. Recall that since 
we are considering the subcase where (\ref{hypothese2}) falls under the special case of 
Proposition \ref{giade} in case IIA, it follows that for each $C^{l,i_1\dots i_\mu}_g$
 there is a unique factor 
in the form $\nabla^{(m)}R_{ijkl}$ for which two internal indices are free, and each 
derivative index is either free or contracting against a 
factor $\nabla\phi_h$; denote by $M+2$ the number of 
free indices in that factor.\footnote{So we set $\alpha=M+2$.}

Now,  by weight considerations (since we are in a special subcase 
of Proposition \ref{giade} in case IIA), any tensor field of rank $>\mu$ in (\ref{hypothese2}) must have
 strictly fewer than $M+2$ free indices in any given factor $\nabla^{(m)}R_{ijkl}$. 
Therefore in subcase IA, (\ref{hypothese2})
 is of the form (\ref{gates}), with $L_{Max}\subset \overline{L}$. 
Therefore, we apply Lemma \ref{specialii} to (\ref{hypothese2}) and pick out the sublinear 
combination of terms with a refined double character $Doub(\vec{L}^z), z\in Z'_{Max}$\footnote{Recall 
that $\vec{L}^z$, $z\in Z'_{Max}$ is the collection of maximal 
refined double characters that Proposition \ref{giade} deals with.}
 we thus obtain a {\it new} true equation, 
since (\ref{tasidera}) holds formally, and the double character is 
 invariant under the formal permutations of indices that make (\ref{tasidera}) formally zero. 
This proves our claim in subcase IIA. 
 \newline
 
 Now we deal with subcase IIB:

We consider the $\mu$-tensor fields of maximal refined double character 
in (\ref{hypothese2}). By definition (since we now fall under a special case), 
they will each have a factor in the form 
$\nabla^{(m)}_{(free)\dots (free)}R_{(free)jkl}$, 
with a total of $M+1>1$ free indices.\footnote{So, we set $\alpha=M+1$.} 
Each of the other factors will be in the form $R_{ijkl}$ 
or be simple factors in the form 
$S_{*}R_{ijkl}$, or in the form $\nabla^{(2)}\Omega_h$.

We denote by $\overline{L}\subset L$ the index set of $\mu$-tensor fields 
with $M+1$ free indices in a factor 
$\nabla^{(m)}R_{ijkl}$. It follows by weight considerations 
that the factor in question will be unique for each $C^{l,i_1\dots i_\mu}_g, l\in \overline{L}$.
We then start out with some explicit manipulation of the terms indexed in $\overline{L}$:

We will prove that there exists a linear combination of $\mu+1$-tensor 
fields, $\sum_{h\in H} a_h C^{h,i_1\dots i_{\mu+1}}_g$, as allowed 
in the statement of Proposition \ref{giade}, so that:

\begin{equation}
\label{injury}
\begin{split}
& \sum_{l\in \overline{L}} a_l 
C^{l,i_1\dots i_\mu}_g\nabla_{i_1}\upsilon\dots\nabla_{i_\mu}\upsilon=
\sum_{h\in H} a_h Xdiv_{i_{\mu+1}}  
C^{h,i_1\dots i_{\mu+1}}_g\nabla_{i_1}\upsilon\dots\nabla_{i_\mu}\upsilon
\\&+\sum_{l\in \overline{L}_{new}} a_l 
C^{l,i_1\dots i_\mu}_g\nabla_{i_1}\upsilon\dots\nabla_{i_\mu}\upsilon
\sum_{j\in J} a_j C^{l,i_1\dots i_\mu}_g\nabla_{i_1}\upsilon\dots\nabla_{i_\mu}\upsilon.
\end{split}
\end{equation}
Here the $\mu$-tensor fields indexed in $\overline{L}_{new}$ have a factor 
\\$T=\nabla^{(M-1)}_{(free)\dots (free)}R_{(free)j(free)l}$, 
and {\it one other factor} $T'$ 
has an extra derivative (meaning that $T'$ is either in the form 
$\nabla R_{ijkl}$ or $S_{*}\nabla R_{ijkl}$, or $\nabla^{(3)}\Omega_h$). 
Moreover if both indices ${}_j,{}_l$ in $T$ contract against indices ${}^j,{}^l$ in 
the same factor $T''$ and at least one of ${}^j,{}^l$ is removable, 
then $T'\ne T''$. Clearly, (\ref{injury}) in conjunction with Lemma
 \ref{specialii} implies Proposition \ref{giade} in the ``special cases'' 
of case II. So, matters are reduced to showing (\ref{injury}) 
(and then deriving Lemma \ref{specialii}). 
\newline

{\it Proof of (\ref{injury}):} 
We first apply the second Bianchi identity to the factor $T$ 
to move one of the derivative free indices ito the position ${}_{kl}$ 
in the factor $\nabla^{(M-1)}_{(free)\dots (free)}R_{(free)j(free)l}$. 
Thus, we derive that modulo terms of length $\ge\sigma+u+1$: 

$$C^{l,i_1\dots i_\mu}_g=-C^{l,1,i_1\dots ,i_\mu}_g +C^{l,2,i_1\dots ,i_\mu}_g,$$
where the partial contractions $C^{l,1,i_1\dots i_\mu}_g, C^{l,2,i_1\dots i_\mu}_g$ have 
the factor $T$ replaced by a factor in the form:
 $\nabla^{(m)}_{k(free)\dots (free)} R_{(free)j(free)l}$,
$\nabla^{(m)}_{l(free)\dots (free)} R_{(free)jk(free)}$,
respectively. We then erase the indices ${}_k,{}_l$ in these two factors 
(thus creating a new tensor field 
$C^{l,1,i_1\dots ,i_\mu i_{\mu+1}}_g, C^{l,2,i_1\dots ,i_\mu i_{\mu+1}}_g$)
by creating a free index ${}_{i_{\mu+1}}$),
 and subtract the $Xdiv_{i_{\mu+1}}[\dots]$ of the corresponding 
$(\mu+1)$-tensor field. We then derive an equation: 

\begin{equation}
\label{surgery}
 C^{l,1,i_1\dots ,i_\mu}_g =Xdiv_{i_{\mu+1}}C^{l,1,i_1\dots i_{\mu+1}}_g+
\sum_{l\in L_{new}} C^{l,i_1\dots i_\mu}_g,
\end{equation}
where all the tensor fields indexed in $L_{new}$ satisfy 
the required property of Lemma \ref{specialii}, except
for the fact that one could have both indices ${}_j,{}_l$ 
in the factor $\nabla^{(M-1)}_{(free)\dots (free)}R_{(free)j(free)l}$
contracting against indices ${}^j,{}^l$ in a factor $T'$ which has 
an additionnal derivative index. If $C^{l,i_1\dots i_\mu}_g$, $l\in L_{new}$  
is not in the form allowed in the claim of Lemma \ref{specialii}, then 
(after possibly applying the second Bianchi identity and possibly introducing 
simply subsequent complete contractions) we may arrange  that one of the 
indices ${}^j,{}^l$ is a derivative index.

In that case we construct another 
$(\mu+1)$-tensor field by erasing the derivative index ${}^j$ or ${}^l$ and making 
the index ${}_j$ or ${}_l$ in a free index ${}_{i_{\mu+1}}$. 
Then, subtracting the corresponding $Xdiv_{i_{\mu+1}}$ of this new $(\mu+1)$-tensor 
field, we derive our claim.  $\Box$

Therefore, matters are reduced to proving Lemma \ref{specialii}. 
\newline

{\it Proof of Lemma \ref{specialii}:}
\newline

Let us start with some notational conventions

\par Recall the first variation law of the curvature tensor under
variations by a symmetric 2-tensor by $v_{ij}$: For any complete
or partial contraction $T(g_{ij})$ (which is a function of the
metric $g_{ij}$), we define:
$Image^1_{v_{ij}}=\frac{d}{dt}|_{t=0}[T(g_{ij}+tv_{ij})]$.
 (We write $Image^1_{v_{ij}}[\dots]$
 or  $Image^1_{v_{ab}}[\dots]$ below to stress
 that we are varying by a 2-tensor, rather than just by a scalar).

We consider the equation $Image^1_{v_{ij}}[L_g]=0$ (which
corresponds to the first {\it metric} variation of our Lemma
hypothesis (i.e.~of (\ref{hypothese2})). This equation holds modulo
 complete contractions with at least $\sigma+u+1$ factors.

\par Thus, we derive a new local equation:

\begin{equation}
\label{elalam}
\begin{split}
&\sum_{l\in L_\mu} a_l Xdiv_{i_1}\dots Xdiv_{i_\mu}
Image^{1}_{v_{ab}}[C^{l,i_1\dots i_\mu}_g]\\&+\sum_{l\in L\setminus L_\mu} a_l
Xdiv_{i_1}\dots Xdiv_{i_a}Image^{1}_{v_{ab}}[C^{l,i_1\dots i_a}_g]
\\&=\sum_{j\in J} a_j Image^{1}_{v_{ab}}[C^j_g],
\end{split}
\end{equation}
which holds modulo terms of length $\ge\sigma+u+1$. 

\par Now, we wish to pass from the local equation 
above to an integral equation, and then to apply 
the {\it silly divergence formula} from \cite{a:dgciI} to that integral equation 
(thus deriving a {\it new} local equation). 

In order to do this, we start by introducing some more notation:
Let us write out:
$$Image^{1}_{v_{ab}}[C^{l,i_1\dots i_\mu}_g]=\sum_{t\in T^l} a_t C^{t,i_1\dots i_a}_g$$
where each $C^{t,i_1\dots i_a}_g$ is in the form:

\begin{equation}
\begin{split}&pcontr(\nabla^{(A+2)}_{r_1\dots r_{A+2}}v_{ab}\otimes
\nabla^{(m_1)}R_{ijkl}\otimes\dots\otimes \nabla^{(m_{\sigma-1})}R\otimes \nabla^{(b_1)}\Omega_1
\otimes\dots\otimes\nabla^{(b_p)}\Omega_p\\&\otimes \nabla\phi_1\otimes\dots\otimes\nabla\phi_u).
\end{split}
\end{equation}

\par For our next technical tool we introduce some notation: For each tensor field
$C^{l,i_1\dots i_a}_{g}$ in the form above, we denote by $C^l_{g}$ the complete
contraction that arises by hitting each factor $T_i$ ($i=1,2,3$)
by $m$ derivative indices $\nabla^{u_1\dots u_m}$, where
${}_{u_1},\dots ,{}_{u_m}$ are the free indices that belong to
$T_i$ in $C^{l,i_1\dots i_a}_g$ (thus we obtain a factor with $m$
 internal contraction, each involving a derivative index).
Notice there is a one-to-one correspondence between the tensor fields and
 the complete contractions we are constructing.
 We can then easily observe that there are two linear combinations
$\Sigma_{r\in R_1} a_r C^r_{g}(\Omega_1,\dots
\Omega_p,\phi_1,\dots ,\phi_u)$, $\Sigma_{r\in R_2} a_r
C^r_{g}(\Omega_1,\dots \Omega_p,\phi_1,\dots ,\phi_u)$ where each
$C^r_g,r\in R_1$ has at least $\sigma+u+1$ factors, while each
$C^r_g, r\in R_2$ has $\sigma +u$ factors but at least one factor
$\nabla^{(p)}\phi_h\ne\Delta\phi_h$ with $p\ge 2$, so that for any
compact orientable $(M,g)$:

\begin{equation}
\begin{split}
\label{tejada}
&\int_M\sum_{l\in L} a_l \sum_{t\in T^l} a_t C^{t,*}_g(v_{ab})+
\sum_{r\in R_1} a_r C^r_g(v_{ab})+\sum_{r\in R_2} a_r C^r_g(v_{ab})dV_g=0
\end{split}
\end{equation}
(denote the integrand of the above by $Z_g(v_{ab})$).
Here again each $C^j_{g}$ has $\sigma +u$ factors and all
factors $\nabla\phi_h$ have only one derivative but its simple
character is subsequent to $\vec{\kappa}$. We call this technique
 (of going from the local equation (\ref{elalam}) to the
  integral equation (\ref{tejada})) the ```inverse integration by parts''.
\newline

\par Now, we derive a ``silly divergence formula'' from the
above by performing integrations by parts with respect to
the factor $\nabla^{(B)}v_{ab}$ (until we are left
with a factor $v_{ab}$--without derivatives). This produces a new local equation
which we denote by $silly[Z_g(v_{ab})]=0$. 
We will be using this equation in our derivation of Lemma \ref{specialii}.

Now, for each $C^{l,i_1\dots i_\mu}_g$, $l\in \overline{L}$, we 
consider the factor \\$T=\nabla^{(M)}_{(free)\dots (free)}R_{(free)j(free)l}$ 
with the $M+2$ free indices. We define $T^j$ to be the 
factor in $C^{l,i_1\dots i_\mu}_g$
that contracts against the index ${}_j$ in $T$ and by $T^l$ to be the factor in $C^{l,i_1\dots i_\mu}_g$
that contracts against the index ${}_l$ in $T$. We define $\overline{L}_{same}\subset \overline{L}$
to be the index set of tensor fields for which $T^j=T^l$; 
we define $\overline{L}_{not.same}\subset \overline{L}$
to be the index set of tensor fields for which $T^j\ne T^l$.
We will then prove (\ref{tasidera}) separately for the two sublinear combinations indexed in  
$\overline{L}_{same},\overline{L}_{not.same}$.
\newline

{\it Proof of (\ref{tasidera}) for the index set $\overline{L}_{same}$:} 

We first prove our claim for $\sigma>3$ and then note how to prove it when $\sigma=3$. 

Consider $silly[L_{g}(\Omega_1,\dots,\Omega_p,\phi_1,\dots,\phi_u,v_{ab})]=0$.
Pick out the sublinear combination $silly_+[L_{g}(\Omega_1,\dots,\Omega_p,\phi_1,\dots,\phi_u,v_{ab})]=0$ 
with $\mu-M-2$ internal contractions, and with the indices in the factor 
$v_{ab}$ contracting against a factor $T'$ which either has {\it no} extra derivative indices,
or if it does, then the contraction is according to the pattern $v^{ab}\otimes\nabla_sR_{ajbl}$;
 we also require that the two factors $T'',T'''$ with an extra $M+2$ extra derivatives each.
This sublinear combination must vanish separately, hence we derive: 
\begin{equation}
 \label{name}
silly_+[Z_{g}(\Omega_1,\dots,\Omega_p,\phi_1,\dots,\phi_u,v_{ab})]=0.
\end{equation}
We also observe that this sublinear combination can only arise (in the process of
 passing from the equation $L_g=0$ to deriving $silly_+[Z_g(v_{ab})]=0$) by 
 replacing the factor $\nabla^{(M)}_{(free)\dots (free)}R_{(free)j(free)l}$
 by $\nabla^{(M)}_{(free)\dots (free)}v_{jl}$ and then (in the inverse integration by parts) replacing 
all $\mu$ free indices by internal contractions,\footnote{(Thus the factor 
  $\nabla^{(M)}_{(free)\dots (free)}v_{jl}$ gets replaced by $\Delta^{M+2}v_{ij}$).}
  and finally integrating by parts the $M+2$ pairs of derivative 
  indices $(\nabla^a,\nabla_a)$ and forcing all upper indices hit a factor $T''\ne T'$ and the lower indices
to hit a factor $T'''\ne T', T'''\ne T''$.\footnote{The fact that 
$\sigma>3$ ensures the existence of two such factors.} 

Thus, we can prove our claim by starting from the equation (\ref{name}) and 
appling $Sub_\upsilon$ $\mu-M-2$ times,\footnote{See the
 Appendix in \cite{alexakis1} and just set $\omega=\upsilon$.}  just applying the eraser
to the extra $M+2$ pairs of contracting derivatives,\footnote{This can be done by just 
repeating the proof of the ``Eraser'' Lemma in the Appendix of \cite{alexakis1}.} 
and then replacing the factor $v_{ab}$ by 
\\$\nabla^{(M)}_{r_1\dots r_M}R_{iajb}\nabla^{r_1}\upsilon\dots\nabla^{r_M}\upsilon\nabla^a\upsilon\nabla^b\upsilon$. 
Finally we just divide by the combinatiorial constant ${\sigma-3}\choose{2}$.

\par Let us now consider the case $\sigma=3$: In those case the terms of 
maximal refined double character can only arise in 
the subcase IIA,\footnote{This follows by virtue of the
 symmetry of the indices ${}_s,{}_a,{}_d$ in any factor $\nabla_sR_{abcd}$ as discussed above.}
 and can only be in one of the  forms: 
$\nabla^{(M)}_{(free)\dots (free)}R_{(free)j(free)l}\otimes R^{ijkl}\otimes\nabla^{(2)}_{ik}\Omega_1)$,
$\nabla^{(M)}_{(free)\dots (free)}R_{(free)j(free)l}\otimes R^{ijkl}\otimes{{{R_{(free)}}^j}_{(free)}}^l)$.
 Thus, in that case we define $silly_+[Z_g(v_{ab})]$ to stand for the terms  
 $(v_{jl}\otimes\nabla^{(M+2)}_{t_1\dots t_{M+2}} R^{ijkl}\otimes\nabla^{(M+4)}_{t_1\dots t_{M+4}ik}\Omega_1)$,
 $(v_{jl}\otimes \nabla^{(M+2)}_{t_1\dots t_{M+2}}R^{ijkl}
(\nabla^{(M+2)})^{t_1\dots t_{M+2}}\otimes{{{R_{(free)}}^j}_{(free)}}^l)$ respectively,
and then repeat the argument above. 
\newline

{\it Proof of (\ref{tasidera}) for the index set $\overline{L}_{not.same}$:} 

We prove our claim in steps: We first denote by 
$\overline{L}_{not.same}^{**}\subset \overline{L}_{not.same}$
 the index set of tensor fields in $\overline{L}_{not.same}$ 
for which both indices ${}_j,{}_l$ in the factor 
 $T=\nabla^{(M)}_{(free)\dots (free)}R_{(free)j(free)l}$ contract
 against special indices in factors $T^j,T^l$ of the form $S_{*}R_{ijkl}$. 
 We will firstly prove that: 
 
 \begin{equation}
 \label{thisplace}
\sum_{l\in \overline{L}_{not.same}^{**}} a_l C^{l,i_1\dots i_\mu}_g \nabla_{i_1}\upsilon\dots\nabla_{i_\mu}\upsilon=
\sum_{l\in L'} a_l   C^{l,i_1\dots i_\mu}_g \nabla_{i_1}\upsilon\dots\nabla_{i_\mu}\upsilon.
 \end{equation}
 Here the terms in the RHS have all the features of the terms 
in $\overline{L}_{not.same}$, but in addition at most
 one of the indices in the factor  $T=\nabla^{(M)}_{(free)\dots (free)}R_{(free)j(free)l}$ 
 contract against a special index in a factor of the form $S_{*}R_{ijkl}$. 
 Thus, if we can prove (\ref{thisplace}), we are reduced to proving our claim 
 under the additional assumption that $\overline{L}_{not.same}^{**}=\emptyset$.
 
 For our next claim, we denote by  $\overline{L}_{not.same}^{*}\subset \overline{L}_{not.same}$
 the index set of tensor ields in $\overline{L}_{not.same}$ for which 
 one of the indices ${}_j,{}_l$ in the factor 
 $T=\nabla^{(M)}_{(free)\dots (free)}R_{(free)j(free)l}$ contracts 
 against a special index in factors $T^j,T^l$ of the form $S_{*}R_{ijkl}$.

 We will then prove that: 
 
 \begin{equation}
 \label{thisplace2}
\sum_{l\in \overline{L}_{not.same}^{*}} a_l C^{l,i_1\dots i_\mu}_g \nabla_{i_1}\upsilon\dots\nabla_{i_\mu}\upsilon=
\sum_{l\in L''} a_l   C^{l,i_1\dots i_\mu}_g \nabla_{i_1}\upsilon\dots\nabla_{i_\mu}\upsilon.
 \end{equation}
 Here the terms in the RHS have all the features of the terms in $\overline{L}_{not.same}$, but in addition 
 none of the indices in the factor  $T=\nabla^{(M)}_{(free)\dots (free)}R_{(free)j(free)l}$ 
 contracting against a special index in a factors of the form $S_{*}R_{ijkl}$. 
 Thus, if we can prove (\ref{thisplace}), we are reduced to proving our claim 
 under the additionnal assumption that for eac $C^{l,i_1\dots i_\mu}_g$, $l\in \overline{L}$
 the two indices ${}_j,{}_l$ in the factor 
$T=\nabla^{(M)}_{(free)\dots (free)}R_{(free)j(free)l}$  contract
 against two different factors and none of the indices ${}^j,{}^l$ are 
 special indices in a factor  of the form $S_{*}R_{ijkl}$. 
  
  In our third step, we prove (\ref{tasidera}) under this additonnal assumption.

    We now present our proof of the third step. We will indicate in the end how this
 proof can be easily modified to derive the first two steps. 
  \newline

  For each $l\in \overline{L}_{not.same}$, let us denote by $link(l)$ 
  the number of particular contractions betwen the 
  factors $T^j,T^l$ in the tensor fields $C^{l,i_1\dots i_\mu}_g$. 
(Note that by weight considerations $0\le link(l)\le 3$). 
  Let $B$ be the maximum value of $link(l), l\in \overline{L}_{not.same}$, and by 
  $\overline{L}_{not.same}^B\subset \overline{L}_{not.same}$ the corresponding 
index set. We will then prove our claim 
  for the tensor fields indexed in $\overline{L}_{not.same}^B$. 
  By repeating this step at most four times, we will derive our third claim.

Consider $silly[L_{g}(\Omega_1,\dots,\Omega_p,\phi_1,\dots,\phi_u,v_{ab})]=0$.
Pick out the sublinear combination $silly_*[L_{g}(\Omega_1,\dots,\Omega_p,\phi_1,\dots,\phi_u,v_{ab})]=0$ 
with $\mu-M-2$ internal contractions, and with an extra $M+2$ 
derivatives on the factors $T^j,T^l$ against which the two indices of the 
factor $v_{ab}$ contract, {\it and} with $M+2+B$ particual contractions betwen the
factors $T^j,T^l$. This sublinear combination must vanish separately: 
$$silly_*[L_{g}(\Omega_1,\dots,\Omega_p,\phi_1,\dots,\phi_u,v_{ab})]=0.$$
Moreover, we observe by following the ``inverse integration 
by parts'' and the silly 
divergence formula obtained from $\int_{M^n} Z_g(v_{ab})dV_g=0$, 
that the LHS of the above can be desrcibed as follows: 

 For each $C^{l,i_1\dots i_\mu}_g, l\in \overline{L}^B_{not.same}$, we denote by $\tilde{C}_g^l(v_{ab})$ 
 the complete contraction that arises by replacing 
 the factor $T=\nabla^{(M)}_{(free)\dots (free)}R_{(free)j(free)l}$ 
 by $\nabla^{(M+2)}_{(free)\dots (free)}v_{jl}$, and then
 replacing each free index that {\it does not} 
 belong to the factor $T$ by an internal contraction. We then denote by $\hat{C}^l_g(v_{ab})$
 the complete contraction that rises from $\tilde{C}_g^l(v_{ab})$ by hitting the factor $T^j$ 
 (against which the index ${}_j$ in $v_{jl}$ contracts) by $(M+2)$ derivative indices 
 $\nabla_{t_1},\dots ,\nabla_{t_{M+2}}$ and hitting the factor $T^l$ (against 
which the index ${}_l$ in $v_{jl}$ contracts) by derivatives 
  $\nabla^{t_1},\dots ,\nabla^{t_{M+2}}$.\footnote{These 
derivatives contract against the indices
  $\nabla_{t_1},\dots ,\nabla_{t_{M+2}}$ that have hit $T^j$.}
  It follows that:
  
  $$(0=)silly_*[L_{g}(\Omega_1,\dots,\Omega_p,\phi_1,\dots,\phi_u,v_{ab})]=
  \sum_{ \overline{L}^B_{not.same}} a_l 2^{M+1}[\hat{C}^l_g(v_{ab})].$$
Now, to derive our claim, we introduce a formal operation $Op[\dots]$ 
which acts on the terms above by applying $Sub_\upsilon$ to each of the 
$\mu-M-2$ internal contractions,\footnote{See the Appendix of \cite{alexakis1}
for the definition of this operation.} erasing $M+2$ particular 
contractions between the factors $T^j,T^l$ and then replacing the factor $v_{jl}$ by 
$\nabla^{(m)}_{r_1\dots r_M}R_{ijkl}\nabla^{r_1}\upsilon\dots\nabla^{r_M}\upsilon\nabla^i\upsilon\nabla^k\upsilon$. 
This operation produces a new true equation; after we divide this new true 
equation by $2^{M+1}$, we derive our claim. 
\newline

{\it Note on the derivation of (\ref{thisplace}), (\ref{thisplace2}):} 
The equations can be derived by a straighttforward modification of the ideas above:
The only extra feature we must add is that in the silly divergence formula we 
must pick out the terms for which (both/one of the) indices ${}_j,{}_l$ in $v_{jl}$ 
contract against a special index in a factor 
$S_{*}\nabla^{(M+2)}R_{abcd}\nabla^a\tilde{\phi}_h$.
This linear combination will vanish,
 modulo terms where one/none of the indices ${}_j,{}_l$ 
in $v_{jl}$ contract against a special index 
in the factor $S_{*}R_{ijkl}$: This follows by the same 
argument that is used in \cite{alexakis4} to derive that Lemma 3.1 in \cite{alexakis4} 
implies Proposition \ref{giade} in case I: We 
firstly replace the factor $v_{jl}$ by an expression $y_{(j}y_{l)}$. 
We then just replace both/one of the expressions 
$\nabla_i\tilde{\phi}_h, y_j$ by $g_{ij}$ and 
apply $Ricto\Omega$ twice/once.\footnote{Recall that this 
operation has been defined in the Appendix of \cite{alexakis1}
and produces a true equation.} The only terms that survive this true equation 
are the ones indexed in $\overline{L}_{not.same}$, for which
 the expression(s) $S_*\nabla^{(\nu)}_{r_1\dots r_\nu}R_{ijkl}\nabla^i\tilde{\phi}_h\nabla^ky$
are replaced by $\nabla^{(\nu+2)}_{r_1\dots r_\nu jl}Y_f$. 
We then proceed as above, deriving that the 
sublinear combination of terms indexed in $\overline{L}_{not.same}$ 
 must vanish, {\it after} we replace two/one expressions 
 $S_*\nabla^{(\nu)}_{r_1\dots r_\nu}R_{ijkl}\nabla^i\tilde{\phi}_h\nabla^ky$ 
 by $\nabla^{(\nu+2)}_{r_1\dots r_\nu jl}Y_f$. Then, repeating the 
 permutations applied to any factors $\nabla^{(\nu+2)}_{r_1\dots r_\nu jl}Y_f$,
  to $S_*\nabla^{(\nu)}_{r_1\dots r_\nu}R_{ijkl}\nabla^i\tilde{\phi}_h\nabla^ky$  we derive our claim. $\Box$

\subsection{The remaining cases of Proposition \ref{giade} in case III:}

\par We recall that there are remaining cases only when $\sigma=3$.
In that case we have the remaining cases when $p=3$ and $n-2u-2\mu\le 2$,
 or when $p=2,\sigma_2=1$ and $n=2u+2\mu$.

{\it The case $p=3$:} Let us start with the subcase $n-2u-2\mu=0$. In this case, all
tensor fields in (\ref{hypothese2}) will be in the form:

\begin{equation}
\label{thmbeki} \begin{split} &pcontr(\nabla^{(A)}_{i_1\dots
i_aj_1\dots j_b}\Omega_1\otimes \nabla^{(B)}_{i_{a+1}\dots
i_{a+a'}j_{b+1}\dots j_{b+b'}}\Omega_2\otimes
\\&\nabla^{(C)}_{i_{a+a'+1}\dots i_{a+a'+a''}j_{b+b'+1}\dots
j_{b+b'+b''}}\Omega_3
\otimes\nabla^{j_{x_1}}\phi_1\dots\otimes\nabla^{x_{j+j'+j''}}\phi_u)
\end{split}
\end{equation}
where we are making the following conventions: Each of the indices
${}_{i_f}$ is free; also, each of the indices ${}_{j_f}$ is
contracting against some factor $\nabla\phi_h$, and also $A,B,C\ge
2$.

\par Thus, we observe that is this subcase $\mu$ is also the {\it maximum}
rank among the tensor fields appearing in
 (\ref{hypothese2}). Now, assume that the $\mu$-tensor fields
in (\ref{hypothese2}) of {\it maximal refined double character}
have $a=\alpha,a'=\alpha',a''=\alpha''$. With no
loss of generality (only up to renaming the factors
 $\Omega_1,\Omega_2,\Omega_3,\phi_1,\dots,\phi_u$)
 we may assume that $\alpha\ge\alpha'\ge\alpha''$ and that only the functions
  $\nabla\phi_1,\dots,\nabla\phi_{u_1}$ contract
   against $\nabla^{(A)}\Omega_1$ in $\vec{\kappa}_{simp}$.
We will then show that the coefficient $a_{\alpha,\alpha',\alpha''}$
of this tensor field must be zero. This will prove Proposition \ref{giade} in this subcase.

\par We prove that  $a_{\alpha,\alpha',\alpha''}=0$ by considering
the global equation $\int Z_g dV_g=0$ and considering the silly
divergence formula $silly[Z_g]=0$. We then consider the sublinear
 combination $silly_{+}[Z_g]$ consisting of terms
with $\alpha',\alpha''$ internal contractions in the factors
$\nabla^{(D)}\Omega_2,\nabla^{(E)}\Omega_3$, with $\alpha$ particular
contractions between those factors and with all factors
$\nabla\phi_h$ that contracted against $\nabla^{(A)}\Omega_1$
in $\vec{\kappa}_{simp}$ being replaced by $\Delta\phi_h$, while all
factors $\nabla\phi_h$ that contracted against
$\nabla^{(B)}\Omega_2,\nabla^{(C)}\Omega_3$ still do so. We easily observe that
$silly_{+}[Z_g]=0$, and furthermore $silly_{+}[Z_g]$ consists of the complete contraction:

\begin{equation}
\begin{split} 
&contr(\Omega_1\otimes
{\nabla^{f_1\dots f_{\alpha}}}_{j_{b+1}\dots
j_{b+b'}}\Delta^{\alpha'}\Omega_2\otimes \nabla_{f_1\dots
f_{\alpha}j_{b+b'+1}\dots j_{b+b'+b''}}\Delta^{\alpha''}
\Omega_3\otimes\Delta\phi_1\\&\dots
\nabla^{j_{x_{b+b'+b''}}}\phi_u)
\end{split}
\end{equation} 
times the constant
$(-1)^{u_1}2^{\alpha}a_{\alpha,\alpha',\alpha''}$. Thus, we derive
that $a_{\alpha,\alpha',\alpha''}=0$.
\newline

{\it The second subcase:} We now consider the setting where
$\sigma=p=3$, $n-2u-2\mu=2$. In this setting, the maximum rank of
the tensor fields appearing in (\ref{hypothese2}) is $\mu+1$. In
this case, all $(\mu+1)$-tensor fields in (\ref{hypothese2}) will
be in the form (\ref{thmbeki}) (with
$\alpha+\alpha'+\alpha''=\mu+1$, while all the $\mu$-tensor fields
will be in the form (\ref{thmbeki}) but with
$\alpha+\alpha'+\alpha''=\mu$, {\it and} with
 one particular contraction ${}_c,{}^c$ between two of the factors
 $\nabla^{(A)}\Omega_1,\nabla^{(B)}\Omega_2,\nabla^{(C)}\Omega_3$.

\par Now, if both the indices ${}_c,{}^c$ described above are removable,
we can explicitly express $C^{l,i_1\dots i_\mu}_g$ as an $Xdiv$ 
of an acceptable $(\mu+1)$-tensor field. Therefore, we are reduced
to showing our claim in this setting where for each $\mu$-tensor
field in (\ref{hypothese2}) at least one of the indices
${}_c,{}^c$ is not removable. Now, let $z\in Z_{Max}$ stand for
one of the index sets for which the sublinear combination
$\sum_{l\in L^z} a_l C^{l,i_1\dots i_\mu}_g$ in (\ref{hypothese2})
indexes tensor fields of maximal refined double character.
 We assume with no loss of generality that for each $l\in L^z$ the
 factors $\nabla^{(A)}\Omega_1$, $\nabla^{(B)}\Omega_2$,
  $\nabla^{(C)}\Omega_3$ have $\alpha\ge\alpha'\ge\alpha''$ free
  indices respectively.\footnote{Recall that by our hypothesis $\alpha'\ge 2$.}
  Therefore, the tensor fields indexed in $L^z$
   can be in one of the following two forms:

\begin{equation}
\label{thmberkia}
\begin{split}
&pcontr(\nabla^c\nabla^{(A)}_{i_1\dots i_\alpha j_1\dots
j_b}\Omega_1\otimes \nabla^{(B)}_{i_{\alpha+1}\dots
i_{\alpha+\alpha'}j_{b+1}\dots j_{b+b'}}\Omega_2\otimes
\\&\nabla^{(2)}_{ci_{\alpha+\alpha'+1}\dots
i_{\alpha+\alpha'+\alpha''}j_{b+b'+1}\dots j_{b+b'+b''}}
\Omega_3\otimes\nabla^{j_{x_1}}\phi_1\dots\otimes\nabla^{x_{j+j'+j''}}\phi_u),
\end{split}
\end{equation}

\begin{equation}
\label{thmberkib}
\begin{split}
&pcontr(\nabla^{(A)}_{i_1\dots i_\alpha j_1\dots
j_b}\Omega_1\otimes \nabla^c\nabla^{(B)}_{i_{\alpha+1}\dots
i_{\alpha+\alpha'}j_{b+1}\dots j_{b+b'}}\Omega_2\otimes
\\&\nabla^{(2)}_{ci_{\alpha+\alpha'+1}\dots
i_{\alpha+\alpha'+\alpha''}j_{b+b'+1}\dots
j_{b+b'+b''}}\Omega_3
\otimes\nabla^{j_{x_1}}\phi_1\dots\otimes\nabla^{x_{j+j'+j''}}\phi_u),
\end{split}
\end{equation}
(where $A,B\ge 3$).

\par Now, by ``manually subtracting'' $Xdiv$'s from these $\mu$-tensor
 fields, we can assume wlog that the tensor fields indexed
  in our chosen $L^z$ are in the from (\ref{thmberkib}).

\par With that extra assumption, we can show that the coefficient
of the tensor field (\ref{thmberkib}) is zero. We see this by
considering the (global) equation $\int_M Z_g dV_g=0$
 and using the silly divergence formula $silly[Z_g]=0$
(which arises by integrations by parts w.r.t. to the factor
$\nabla^{(A)}\Omega_1$). Picking out the sublinear combination
$silly_{+}[Z_g]$ which consists of the complete contraction:

\begin{equation}\begin{split}&contr(\Omega_1\otimes
\nabla^c\nabla^{f_1\dots f_{\alpha}}_{j_{b+1}\dots
j_{b+b'}}\Delta^{\alpha'}\Omega_2\otimes \nabla_c\nabla_{f_1\dots
f_{\alpha}j_{b+b'+1}\dots j_{b+b'+b''}}\Delta^{\alpha''}
\Omega_3\otimes\\&\Delta\phi_1\dots
\nabla^{j_{x_{b+b'+b''}}}\phi_u)
\end{split}
\end{equation} (notice that
$silly_{+}[Z_g]=0$), we derive that the coefficient of
(\ref{thmberkib}) must vanish. Thus, we have shown our claim in
this second subcase also. $\Box$
\newline

{\it The case $p=2$, $\sigma_2=1$:} Recall that in this case
we fall under the special case when $n=2u+2\mu$. In this setting, we will have
that in each index set $L^z,z\in Z'_{Max}$ (see the statement of
 Lemma 3.5 in \cite{alexakis4}) there is
a {\it unique} $\mu$-tensor field of maximal refined
double character in (\ref{hypothese2}), where the two indices
 ${}_k,{}_l$ in the factor $S_{*}\nabla^{(\nu)}R_{ijkl}$ will
 be contracting against one of the factors
 $\nabla^{(A)}\Omega_1,\nabla^{(B)}\Omega_2$ (wlog we may
  assume that they are contracting against {\it different factors}).
But now, recall that since we are considering case $A$ 
of Lemma 3.5 in \cite{alexakis4}, one of the factors
$\nabla^{(A)}\Omega_1,\nabla^{(B)}\Omega_2$ will have at least two free indices.
 Hence, in at least one of the factors
 $\nabla^{(A)}\Omega_1,\nabla^{(B)}\Omega_2$,
  the index ${}^k,{}^l$ is removable
(meaning that it can be erased, and
we will be left with an acceptable tensor field). We denote
 by $C^{l,i_1\dots i_\mu i_{\mu+1}}_g$ the tensor field that arises from
$C^{l,i_1\dots i_\mu}_g$ by erasing the aforementioned
 ${}^k,{}^l$ and making ${}_k$ or ${}_l$ into a free index, we then observe that:

 \begin{equation}
\label{kyriakos}
C^{l,i_1\dots i_\mu}_g-Xdiv_{i_{\mu+1}}C^{l,i_1\dots i_\mu i_{\mu+1}}_g=0
 \end{equation}
(modulo complete contractions of length $\ge\sigma+u+1$).
This completes the proof of our claim. $\Box$

\end{document}